\documentclass[11pt]{amsart}
\usepackage{graphicx,indentfirst,pict2e,amsmath,amsthm,amsfonts,amssymb,mathrsfs}
\newtheorem{theorem}{Theorem}[section]
\newtheorem{theorem*}{Theorem}
\newtheorem{definition}[theorem]{Definition}

\newtheorem{lemma}[theorem]{Lemma}

\newcommand{\C}{\mathbb{C}}

\newcommand{\Q}{\mathbb{Q}}
\newcommand{\incidenceLocus}{I}

\newcommand{\canonicalPower}{\nu}
\newcommand{\dA}{d}
\newcommand{\NA}{N}
\newcommand{\Ki}{K_{i}}

\newcommand{\wkwi}{W(k,w,i)}

\newcommand{\Pro}{\mathbf{P}}
\newcommand{\Hilb}{\mbox{\it{Hilb}}}
\newcommand{\Cliff}{\mbox{\scriptsize{\it{Cliff}}}}
\newcommand{\RR}{\mbox{\scriptsize{\it{RR}}}}
\newcommand{\vir}{\protect \mbox{\scriptsize{\textmd{{vir}}}}}
\newcommand{\indz}{\mathbf{1}_{\zeta \neq 0}}
\newcommand{\indzki}{\mathbf{1}_{\zeta_{k,i} \neq 0}}
\newcommand{\dblq}{/\!/}

\newcommand{\Mgn}{\overline{{M}}_{g,n}}
\newcommand{\MgA}{\overline{{M}}_{g,\mathcal{A}}}

\newcommand{\codim}{\operatorname{codim}}
\newcommand{\Span}{\operatorname{Span}}
\newcommand{\wt}{\operatorname{wt}}
\newcommand{\Sym}{\operatorname{Sym}}

\title{GIT stability of weighted pointed curves}
\author{David Swinarski}
\address{Department of Mathematics, Columbia University, New York NY 10027, USA}
\email{swinarski@math.columbia.edu}

\begin{document}

\pagenumbering{arabic}
\maketitle


\setlength{\baselineskip}{15pt}

In the late 1970s Mumford established Chow stability of smooth unpointed genus $g$ curves embedded by complete linear systems of degree $d \geq 2g+1$,  and at about the same time Gieseker established asymptotic Hilbert stability (that is, stability of $m^{th}$ Hilbert points for some large values of $m$) under the same hypotheses.  Both of them then use an indirect argument to show that nodal Deligne-Mumford stable curves are GIT stable.  The case of marked points lay untouched until 2006, when Elizabeth Baldwin proved that pointed Deligne-Mumford stable curves are asymptotically Hilbert stable.  (Actually, she proved this for stable maps, which includes stable curves as a special case.)   Her argument is a delicate induction on $g$ and the number of marked points $n$; elliptic tails are glued to the marked points one by one, ultimately relating stability of an $n$-pointed genus $g$ curve to Gieseker's result for genus $g+n$ unpointed curves.

There are three ways one might wish to improve upon Baldwin's results.  First, one might wish to construct moduli spaces of weighted pointed curves or maps; it appears that Baldwin's proof can accommodate some, but not all, sets of weights.   Second, one might wish to study Hilbert stability for small values of $m$; since Baldwin's proof uses Gieseker's proof as the base case, it is not easy to see how it could be modified to yield an approach for small $m$.  Finally, the Minimal Model Program for moduli spaces of curves has generated interest in GIT for 2, 3, or 4-canonical linear systems; due to its use of elliptic tails, Baldwin's proof cannot be used to study these, as elliptic tails are known to be GIT unstable in these cases.

In this paper I give a direct proof that smooth curves with distinct weighted marked points are asymptotically Hilbert stable with respect to a wide range of parameter spaces and linearizations.  Some of these yield the (coarse) moduli space of Deligne-Mumford stable pointed curves $\Mgn$ and Hassett's moduli spaces of weighted pointed curves $\MgA$, while other linearizations may give other quotients which are birational to these and which may admit interpretations as moduli spaces.   The full construction of the moduli spaces is not contained in this paper, only the proof that smooth curves with distinct weighted marked points are stable, which is the key new result needed for the construction.  For this I follow Gieseker's approach to reduce to the GIT problem to a combinatorial problem, though the solution is very different.

\section*{Introduction}
Let $(C, P_{1},\ldots, P_{n}, \mathcal{A})$ be a weighted pointed stable curve.  That is, 
\begin{itemize}
\item $C$ is a reduced connected projective algebraic curve with at worst nodes as singularities, 
\item the points $P_{i}$ lie on $C$ and are ordered (note we do not require that they be distinct, nor that they be smooth points of $C$),
\item $\mathcal{A} = (a_{1},\ldots,a_{n})$, where the $a_{i}$ are rational numbers between 0 and 1 inclusive,
\item  $a_{i} = 0$ if $P_{i}$ is a node, 
\item a subset of the points is allowed to collide if the sum of their weights does not exceed 1, and
\item the $\Q$-line bundle $\omega(\sum a_{i}P_{i})$ is ample on $C$.
\end{itemize}
Hassett introduced weighted pointed stable curves in \cite{Hass}; the theory is extended to stable maps by several people (\cite{BM}, \cite{AG}, \cite{MM}).  

The goal of this paper is to describe linearizations for which the points of an appropriate space parametrizing embedded weighted pointed stable curves $(C \subset \Pro^{N},  P_{1},\ldots,P_{n}, \mathcal{A})$ are GIT stable.  The main result of this paper, Theorem \ref{stabilitytheorem}, does not say exactly this.  Instead, for most of this paper, we do the following:
\begin{itemize}
\item We ignore the set of weights $\mathcal{A}$ and just study embedded pointed curves $(C \subset \Pro^{N},  P_{1},\ldots,P_{n})$.  
\item We assume that the curve $C$  is smooth.
\item We assume that the points $\{ P_{i} \}$ are distinct.
\end{itemize}
Theorem \ref{stabilitytheorem} asserts that smooth pointed curves with distinct marked points are GIT stable with respect to certain linearizations.  Armed with this result, one may proceed to show that all weighted pointed stable curves are GIT stable for certain linearizations, justifying the title of this paper.  This is not fully written out here, but it is discussed in Section \ref{moduli section}.

So, let $x$ be a point parametrizing an embedded smooth pointed curve $(C \subset \Pro^{N}, P_{1},\ldots,P_{n} )$.  Following Gieseker, the  numerical criterion is reformulated in a way that permits a more combinatorial approach.  A 1-PS $\lambda$ of $SL(N+1)$ induces a weighted filtration of $H^{0}(C,\mathcal{O}(1))$ and a weighted filtration of $H^{0}(C,\mathcal{O}(m))$.  The value of Mumford's function $\mu^{L}(x,\lambda)$ may be interpreted as the ``minimum weight of a basis of $H^{0}(C, \mathcal{O}(m))$ compatible with this filtration plus a contribution from the marked points.''  (From now on, whenever we refer to a basis of $H^{0}(C,\mathcal{O}(m))$, we always implicitly mean one that is compatible with the weighted filtration.)   The numerical criterion says that if $\mu^{L}(x,\lambda)$ is sufficiently small, then $x$ is GIT stable with respect to $\lambda$.  Any basis therefore gives an upper bound for $\mu^{L}(x,\lambda)$, so the goal becomes: find a basis of sufficiently small weight.

Our main tool for computing (a bound for) the weight of a basis is something I call a {\it profile}.  This is a graph which may be associated to any filtration of a vector space such that the weight decreases at each stage.   Suppose $\tilde{F}_{\bullet}$ is such a filtration of $H^{0}(C, \mathcal{O}(m))$.   (I use tildes for filtrations of $H^{0}(C,\mathcal{O}(m))$; no tilde indicates a filtration of $H^{0}(C,\mathcal{O}(1))$.)  Suppose the weight on the $k^{th}$ stage of $\tilde{F}_{\bullet}$ is $\tilde{r}_{k}$.  Then the profile associated to $\tilde{F}_{\bullet}$ is just the decreasing step function in the first quadrant of the $(\mbox{codimension} \times \mbox{weight})$-plane whose value is $\tilde{r}_{k}$ over the interval $[\codim \tilde{F}_{k}, \codim \tilde{F}_{k+1} ).$   Given any profile, it is possible to choose a basis whose weight is less than the area under the profile.

There is a notion of an absolute weight filtration on  $H^{0}(C, \mathcal{O}(m))$ (see Section \ref{profilessection}); the area under its profile is the minimum weight of a basis.  This is perhaps the most natural filtration to consider, but it is too difficult to compute.  So, like Gieseker, we study other filtrations.

The action of a 1-PS $\lambda$ induces a filtration $V_{\bullet}$ of $H^{0}(C,\mathcal{O}(1))$.  By considering specific spaces of degree $m$ monomials in elements of $V$ diagonalizing the $\lambda$-action, Gieseker produces a very straightforward filtration $\tilde{V}_{\bullet}$ of $H^{0}(C,\mathcal{O}(m))$ as well as a second, slightly fancier filtration $\tilde{G}_{\bullet}$.  Gieseker is able to show that the weight (or area) associated to $\tilde{G}_{\bullet}$ is sufficiently small to establish $\lambda$-stability of smooth unpointed curves.  Unfortunately, as we show with a concrete example, the analogue of $\tilde{G}_{\bullet}$ is not sufficient to establish $\lambda$-stability when there are marked points.

One could try to improve $\tilde{G}_{\bullet}$, but it is too difficult (at least for me) to show that the sum of its area and the marked points contribution is sufficiently small.  Therefore I use $\tilde{V}_{\bullet}$ as a starting point to build a new filtration, $\tilde{X}_{\bullet}$, which is obtained by taking spans of carefully chosen spaces of monomials.   The recipe is given in terms of the combinatorics of the base loci of the stages of the filtration $V_{\bullet}$.  Although $\tilde{X}_{\bullet}$ is rather tedious to define, it has the virtue that we can bound the sum of its area and the marked points contribution sufficiently well to show that smooth curves with distinct marked points are stable.  The key new ingredients in my proof are the definition/choice of $\tilde{X}_{\bullet}$; an easy but important lemma (Lemma \ref{trace}) which allows us to compute spans of spaces of monomials in the $V_{j}$'s using multiplicities of points in the base loci; and the combinatorial argument (see the proof of Lemma \ref{creep}) which allows us to effectively bound the sum of the marked points contribution and the area of the profile associated to $\tilde{X}_{\bullet}$.  

Gieseker's proof establishes stability for smooth unpointed curves embedded by complete linear systems of degree $d \geq 2g+1$.  (There are some misleadingly placed hypotheses in \cite{G}, but one can check that everything works with the hypotheses just mentioned.)  At the present time it is necessary for me to make the hypotheses:
\begin{itemize}
\item If $n=0$, the parameter space satisfies $N \geq 2g-2$.
\item If $n \geq 1$, then either the parameter space satisfies $N \geq 2g-1$, or else the linearization satisfies the following condition (the notation is explained in Section \ref{numerical criterion}): $\gamma b > \frac{g-1}{N}$.
\end{itemize}
One might hope to do a little better (see Section \ref{lowerenvelopesarebetter}), but at least this includes the important case of bicanonically embedded pointed curves (i.e. pointed curves embedded by sections of \\ $(\omega(P_{1}+\cdots P_{n}))^{2}$).

Here is an outline of the paper: in Section \ref{GITsetup} I describe the GIT problem carefully, specifying the parameter spaces and linearizations we will consider, and reformulate the numerical criterion in the form we shall use it.  Profiles are also defined here.  In Section \ref{Giesreview} I review Gieseker's proof, with a few enhancements, to fix notation; a reader familiar with Gieseker's proof should be able to read it very quickly.  In Section \ref{keyexample}  I give an example showing why his proof does not suffice for marked points, and a hint illustrating how we will go about fixing it.

Throughout Sections \ref{GITsetup}--\ref{keyexample} we steadily extract combinatorial data from the algebro-geometric action of a 1-PS $\lambda$ acting on the Hilbert point of a weighted pointed stable curve.  The last result of this type is Lemma \ref{trace}, which allows us to compute codimensions of spans of monomial-type sublinear series of $H^{0}(C,\mathcal{O}(m))$ using only the multiplicities of points in the base loci.  After this, the problem becomes almost entirely combinatorial.   

In Section \ref{EWP}, I produce the filtration $\tilde{X}_{\bullet}$ on $H^{0}(C,\mathcal{O}(m))$ which is built using the filtration $\tilde{V}_{\bullet}$ as scaffolding.   The goal is now to show that the area under the profile for $\tilde{X}_{\bullet}$ plus the contribution from the marked points is less than the bound specified by the numerical criterion.

This is established in two steps: first, I describe a second, simpler graph called the {\it virtual profile} which is bounded above by the profile for $\tilde{X}_{\bullet}$.  Basically it is the graph of the piecewise linear function connecting the left endpoints of the steps in the weight profile.  (I'm oversimplifying things a little here---I'm glossing over some rounding errors.)  The virtual profile is not really the profile of any filtration, nor does it compute or bound the weight of a basis; the most rigorous interpretation I have for it is on the level of graphs.  Again, while it is easy to compute the area of the profile (it's a step function, after all!), when it is time to add the contribution from the marked points, it is easier to do this with the virtual profile than with the profile.  In Section \ref{virvsactual} I bound the discrepancy between the areas of the two graphs and show that this is relatively small when $m$ is large.  Then in Section \ref{boundingTvir} I bound the sum of the area under the virtual profile for $\tilde{X}_{\bullet}$ and the weight from the marked points.  Everything comes together in Section \ref{putTogether} to show that smooth pointed curves with distinct marked points have GIT stable Hilbert points, and the application of this to construction of moduli spaces is stated but not proven.  Finally, this preprint concludes with a short section of additional remarks which are likely to be omitted from a published version.

Here is a picture illustrating the profile and virtual profile associated to an example that is explained in detail in Section \ref{illustration}.  Note that I will always fill in the graphs of step functions to obtain staircase figures.  

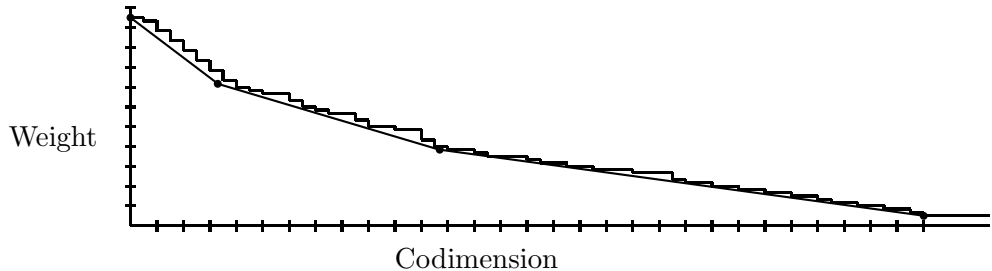
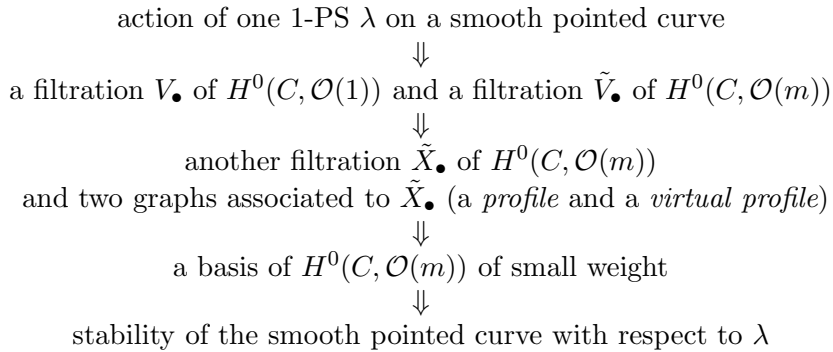
\setlength{\unitlength}{1pt}
\thicklines
\hspace{1in} \begin{picture}(300,90)
\put (0,78.75){\line(1,0){5}}
\put (5,78.75){\line(0,-1){1.25}}
\put (5,77.5){\line(1,0){5}}
\put (10,77.5){\line(0,-1){3.75}}
\put (10,73.75){\line(1,0){5}}
\put (15,73.75){\line(0,-1){3.75}}
\put (15,70){\line(1,0){5}}
\put (20,70){\line(0,-1){3.75}}
\put (20,66.25){\line(1,0){5}}
\put (25,66.25){\line(0,-1){3.75}}
\put (25,62.5){\line(1,0){5}}
\put (30,62.5){\line(0,-1){3.75}}
\put (30,58.75){\line(1,0){5}}
\put (35,58.75){\line(0,-1){3.75}}
\put (35,55){\line(1,0){5}}
\put (40,55){\line(0,-1){2.5}}
\put (40,52.5){\line(1,0){5}}
\put (45,52.5){\line(0,-1){1.25}}
\put (45,51.25){\line(1,0){5}}
\put (50,51.25){\line(0,-1){1.25}}
\put (50,50){\line(1,0){10}}
\put (60,50){\line(0,-1){2.5}}
\put (60,47.5){\line(1,0){5}}
\put (65,47.5){\line(0,-1){2.5}}
\put (65,45){\line(1,0){5}}
\put (70,45){\line(0,-1){1.25}}
\put (70,43.75){\line(1,0){5}}
\put (75,43.75){\line(0,-1){1.25}}
\put (75,42.5){\line(1,0){10}}
\put (85,42.5){\line(0,-1){2.5}}
\put (85,40){\line(1,0){5}}
\put (90,40){\line(0,-1){2.5}}
\put (90,37.5){\line(1,0){10}}
\put (100,37.5){\line(0,-1){1.25}}
\put (100,36.25){\line(1,0){10}}
\put (110,36.25){\line(0,-1){3.75}}
\put (110,32.5){\line(1,0){5}}
\put (115,32.5){\line(0,-1){2.5}}
\put (115,30){\line(1,0){5}}
\put (120,30){\line(0,-1){1.25}}
\put (120,28.75){\line(1,0){10}}
\put (130,28.75){\line(0,-1){1.25}}
\put (130,27.5){\line(1,0){5}}
\put (135,27.5){\line(0,-1){1.25}}
\put (135,26.25){\line(1,0){15}}
\put (150,26.25){\line(0,-1){1.25}}
\put (150,25){\line(1,0){5}}
\put (155,25){\line(0,-1){1.25}}
\put (155,23.75){\line(1,0){10}}
\put (165,23.75){\line(0,-1){1.25}}
\put (165,22.5){\line(1,0){10}}
\put (175,22.5){\line(0,-1){1.25}}
\put (175,21.25){\line(1,0){15}}
\put (190,21.25){\line(0,-1){1.25}}
\put (190,20){\line(1,0){15}}
\put (205,20){\line(0,-1){2.5}}
\put (205,17.5){\line(1,0){5}}
\put (210,17.5){\line(0,-1){1.25}}
\put (210,16.25){\line(1,0){10}}
\put (220,16.25){\line(0,-1){1.25}}
\put (220,15){\line(1,0){10}}
\put (230,15){\line(0,-1){1.25}}
\put (230,13.75){\line(1,0){10}}
\put (240,13.75){\line(0,-1){1.25}}
\put (240,12.5){\line(1,0){10}}
\put (250,12.5){\line(0,-1){1.25}}
\put (250,11.25){\line(1,0){10}}
\put (260,11.25){\line(0,-1){1.25}}
\put (260,10){\line(1,0){5}}
\put (265,10){\line(0,-1){1.25}}
\put (265,8.75){\line(1,0){10}}
\put (275,8.75){\line(0,-1){1.25}}
\put (275,7.5){\line(1,0){10}}
\put (285,7.5){\line(0,-1){1.25}}
\put (285,6.25){\line(1,0){10}}
\put (295,6.25){\line(0,-1){1.25}}
\put (295,5){\line(1,0){5}}
\put (300,5){\line(0,-1){1.25}}
\put (300,3.75){\line(1,0){25}}
\put (0,78.75){\line(33,-25){33}}
\put (0,78.75){\circle*{3}}
\put (33,53.75){\line(84,-25){84}}
\put (33,53.75){\circle*{3}}
\put (117,28.75){\line(183,-25){183}}
\put (117,28.75){\circle*{3}}
\put (300,3.75){\circle*{3}}
\put (0,0){\line(1,0){325}}
\put (0,0){\line(0,1){82.5}}
\put (10,-2){\line(0,1){4}}
\put (20,-2){\line(0,1){4}}
\put (30,-2){\line(0,1){4}}
\put (40,-2){\line(0,1){4}}
\put (50,-2){\line(0,1){4}}
\put (60,-2){\line(0,1){4}}
\put (70,-2){\line(0,1){4}}
\put (80,-2){\line(0,1){4}}
\put (90,-2){\line(0,1){4}}
\put (100,-2){\line(0,1){4}}
\put (110,-2){\line(0,1){4}}
\put (120,-2){\line(0,1){4}}
\put (130,-2){\line(0,1){4}}
\put (140,-2){\line(0,1){4}}
\put (150,-2){\line(0,1){4}}
\put (160,-2){\line(0,1){4}}
\put (170,-2){\line(0,1){4}}
\put (180,-2){\line(0,1){4}}
\put (190,-2){\line(0,1){4}}
\put (200,-2){\line(0,1){4}}
\put (210,-2){\line(0,1){4}}
\put (220,-2){\line(0,1){4}}
\put (230,-2){\line(0,1){4}}
\put (240,-2){\line(0,1){4}}
\put (250,-2){\line(0,1){4}}
\put (260,-2){\line(0,1){4}}
\put (270,-2){\line(0,1){4}}
\put (280,-2){\line(0,1){4}}
\put (290,-2){\line(0,1){4}}
\put (300,-2){\line(0,1){4}}
\put (-2,7.5){\line(1,0){4}}
\put (-2,15){\line(1,0){4}}
\put (-2,22.5){\line(1,0){4}}
\put (-2,30){\line(1,0){4}}
\put (-2,37.5){\line(1,0){4}}
\put (-2,45){\line(1,0){4}}
\put (-2,52.5){\line(1,0){4}}
\put (-2,60){\line(1,0){4}}
\put (-2,67.5){\line(1,0){4}}
\put (-2,75){\line(1,0){4}}
\put (-2,82.5){\line(1,0){4}}
\put (100,-16){Codimension}
\put (-46,30){Weight}
\end{picture}

\vspace{.25in}
In summary: 

\begin{displaymath}
\begin{array}{c}
\mbox{action of one 1-PS $\lambda$ on a smooth pointed curve} \\
\Downarrow \\ 
\mbox{a filtration $V_{\bullet}$ of $H^{0}(C,\mathcal{O}(1))$ and a filtration $\tilde{V}_{\bullet}$ of $H^{0}(C,\mathcal{O}(m))$} \\
\Downarrow \\ 
\mbox{another filtration $\tilde{X}_{\bullet}$ of $H^{0}(C,\mathcal{O}(m))$} \\
\mbox{ and two graphs associated to $\tilde{X}_{\bullet}$ (a {\it profile} and a {\it virtual profile})} \\
\Downarrow \\ 
\mbox{a basis of $H^{0}(C,\mathcal{O}(m))$ of small weight} \\
\Downarrow \\
\mbox{stability of the smooth pointed curve with respect to $\lambda$}  
\end{array}
\end{displaymath}

Two remarks on notation here may reduce anxiety for those skimming the proof:

Note that from Section \ref{EWP} onward it may appear at times as though we are using rational numbers as exponents of monomials.  Although the resulting ``virtual'' spaces are usually nonsensical, in cases where they do make sense they are useful in motivating some definitions and calculations.  However, such spaces are never used to produce basis elements in $H^{0}(C,\mathcal{O}(m))$; to get basis elements, we always round exponents appropriately.

We will obtain two-dimensional arrays of integers $c_{j,i}$.  That is, $j$ indexes the row, and $i$ indexes the column, opposite the usual alphabetic convention.  There is nothing deep happening here; the reasons I made this choice are too silly to discuss further.

\subsection*{Acknowledgements}
It is a great pleasure to thank my advisors, Ian Morrison and Michael Thaddeus, for their help with this work.  I am also very grateful to Elizabeth Baldwin for sharing much of her early work with me, which got me interested in the problem and helped me get started.  Finally, I would like to thank Johan de Jong and Brendan Hassett for their technical help and encouragement.

\section{The GIT setup} \label{GITsetup}

\subsection{The parameter spaces and linearizations we use}
\label{numerical criterion}
In this chapter we investigate GIT stability for the following general setup.  Let $P(t)  := \dA t - g+1$ be a degree one polynomial.  We  form the incidence locus $\incidenceLocus \subset   \Hilb(\Pro^{N},P(t)) \times \prod^{n} \Pro^{N}$ where the points in the projective space factors lie on the curve in $\Pro^N$ parametrized by the point in the first factor.  We study the GIT stability of points of $\incidenceLocus$.  Note two things: no sets of weights $\mathcal{A}$ appear in this paragraph; we will see in Section \ref{moduli section} that considering weighted marked points influences the choice of $d$, but otherwise plays no role in the stability proof.  Also, we do not assume that $C \subset \Pro^N$ is pluricanonically embedded, or even that the degree of $C \subset \Pro^N$ matches the degree of the pluricanonical embedding--- we can investigate GIT stability for more general setups than just those which have an obvious application to construction of moduli spaces of curves.  All we need is that the embedding $C \subset \Pro^N$ is by a complete linear system, and some precise degree/dimension bounds in terms of the genus, which will be carefully stated at the end in Theorem \ref{stabilitytheorem}.  These will even allow some special embeddings.

To do GIT, one must specify a linearization on the $G$-space (here, $\incidenceLocus$).  Although not necessary, perhaps the easiest way to do this is to embed $ \Hilb(\Pro^{N},P(t)) \times \prod^{n} \Pro^{N}$ in a high-dimensional projective space and use its $\mathcal{O}(1)$.

Let $C \subset \Pro^{N}$ be a subscheme with Hilbert polynomial $P(t)$.  For sufficiently large $m, m_{i}'$, the maps 
\begin{displaymath}
\begin{array}{llll}
\mbox{ev}_{C}^{m}: & H^{0}(\Pro^{N},\mathcal{O}(m)) & \rightarrow & H^{0}(C,\mathcal{O}_{C}(m)) \\
\mbox{ev}_{P_{i}}^{m_{i}'}: & H^{0}(\Pro^{N},\mathcal{O}(m_{i}')) & \rightarrow & H^{0}(P_{i},\mathcal{O}_{P_{i}}(m_{i}')) \cong \C
\end{array}
\end{displaymath}
are surjective.   The first map gives rise to an embedding of the Hilbert scheme in a Grassmannian, which in turn embeds in a projective space by the Pl\"{u}cker embedding.  The maps in the second line correspond to $m_{i}'$-uple embeddings of $ \Pro^{N}$.  Finally, a Segre embedding of all these projective spaces yields an embedding of  $ \Hilb(\Pro^{N},P(t)) \times \prod^{n} \Pro^{N}$ into a very large projective space, as desired.

Now, to specify a linearization  on $\incidenceLocus \subset  \Hilb(\Pro^{N},P(t)) \times \prod^{n} \Pro^{N}$, it suffices to specify the ratios between $m$ and each $m_{i}'$.  I will do this as follows: let $\mathcal{B} = (b_{1},\ldots, b_{n}) \in \Q^{n} \cap [0,1]^{n}$ be a set of weights, which I call the \emph{linearizing weights}.  Then set $m_{i}' = \gamma b_{i}m^{2}$.  (The coefficient $\gamma$ will be specified later, at least for the moduli spaces $\MgA$, where it is approximately 1/2; see Section \ref{moduli section}.  Factoring $\gamma$ out of the ratios $m_{i}'/m^{2}$ like this now simplifies the statements of later results needed to construct the moduli spaces.)  Finally, write $b:= \sum_{i=1}^{n} b_{i}$.

\subsection{The numerical criterion for our setup}
By being a little more explicit, we obtain a useful reformulation of the numerical criterion.  

In Gieseker's paper and this paper we use Grothendieck's convention that if $V$ is a vector space, then $\Pro(V)$ is the collection of equivalence classes under scalar action of the nonzero elements of the dual space $V^\vee$.  One consequence of this convention is that the numerical criterion takes the opposite sign from how it appears in \cite{GIT}.

Let $X$ be a projective algebraic scheme with the action of a group $G$ linearized on a very ample line bundle $L$.  Let $\lambda:\mathbb{G}_m\rightarrow G$ be a 1-PS of $G$.  Choose a basis $\{e_0,\ldots ,e_{N}\}$ of $H^0(X,L)$ diagonalizing the $\lambda$ action and ordered so that the weights $r_0\leq\cdots\leq r_{N}\in\mathbb{Z}$ increase.  The weights on the dual basis then have the opposite signs: $-r_0,\ldots,-r_{N}$.

A point $x\in X$ is represented by some non-zero $\hat{x}=\sum_{i=0}^N x_i e_i^\vee\in H^0(X,L)^\vee$.  Define
\[
\mu^L(x,\lambda):=\min\{r_i|x_i\neq 0\}.
\]
Then, with our sign conventions, we have the following characterization of semistability:
\begin{theorem}[cf. \cite{GIT} Theorem 2.1]
\begin{eqnarray*}
x\in X^{ss}(L)	& \Longleftrightarrow	& \mu^L(x,\lambda)\leq 0 \mbox{ for all 1-PS $\lambda\neq 0$} \\
x\in X^{s}(L)	& \Longleftrightarrow	& \mu^L(x,\lambda)< 0 \mbox{ for all 1-PS $\lambda\neq 0$} .
\end{eqnarray*}
\end{theorem}

In our situation $X$ is the incidence scheme $I$, the point $x \in X$ parametrizes an embedded pointed curve $(C \subset \Pro^N,P_{1},\ldots,P_{n})$, the scheme $I$ is embedded in  \mbox{$\Pro(\bigwedge^{P(m)}\Sym^{m}V \otimes \bigotimes^{n} \Sym^{m_{i}'}V)$} where $V = H^{0}(\Pro^N,\mathcal{O}(1))$, and $L$ is the $\mathcal{O}(1)$ on this very large projective space.  Let $\lambda$ be a 1-PS of $SL(V)$.  One particularly nice basis of $ \bigwedge^{P(m)}\Sym^{m}V \otimes \bigotimes^{n} \Sym^{m_{i}'}V$ is given by elements of the form 
\begin{equation} \label{basisform} (M_{1} \wedge \cdots \wedge M_{P(m)}) \otimes (M'_{1}) \otimes \cdots \otimes (M'_{n}),
\end{equation}
where each $M_{j}$ is a monomial of degree $m$ and each $M'_{i}$ is a monomial of degree $m_{i}'$ in the basis elements of $V$ diagonalizing $\lambda$.

The numerical criterion may be translated as follows: a point of 
$I$ is stable with respect to $\lambda$ if and only if there is a basis element of the form (\ref{basisform}) such that 
 \begin{enumerate}
{}{\setlength{\leftmargin}{.75in} \setlength{\rightmargin}{.75in} \setlength{\itemsep}{-1.5mm} \setlength{\itemindent}{-.25in}}
\item the images of the $M_{\ell}$ under the evalution map form a basis of $H^{0}(C,\mathcal{O}_{C}(m))$, 
\item $M'_{i}$ does not vanish at $P_{i}$,
\item the $SL(\NA+1)$ weights satisfy
\begin{displaymath} \displaystyle{\sum_{\ell=1}^{P(m)}} \wt_{\lambda}({M}_{\ell}) + \sum^{n}\wt_{\lambda}(M_{i}') < 0
\end{displaymath}
\end{enumerate}

In fact, it will be convenient to renormalize the $\lambda$ weights so that they decrease to 0 and sum to 1.  If $s_{\NA},\ldots,s_{0}$ are the original weights, (so $s_{\NA} \geq \cdots \geq s_{0}$ and $\sum s_{j} = 0$), then the desired transformation is $r_{j} = (s_{\NA -j}-s_{0})/((\NA +1)|s_{0}|)$.  Also, we write
\begin{eqnarray*} A &:= & \sum_{\ell=1}^{P(m)} \wt_{\lambda}({M}_{\ell}) \\
 T  & := &\sum_{\ell=1}^{P(m)} \wt_{\lambda}({M}_{\ell}) + \sum^{n}\wt_{\lambda}(M_{i}') 
\end{eqnarray*}
for parts of the left hand side of condition $3.$ above.  We may rewrite condition $3.$ as follows.

\begin{lemma}  Condition $3.$ above with the unnormalized weights $s_{j}$ is equivalent to the following condition:
 \begin{enumerate}
\item[3$.'$] With the normalized weights $r_{j}$, the following inequality is satisfied:
{}{\setlength{\leftmargin}{.75in} \setlength{\rightmargin}{.75in} \setlength{\itemsep}{-1.5mm} \setlength{\itemindent}{-.25in}}
\begin{eqnarray} T:= \displaystyle{\sum_{\ell=1}^{P(m)}} \wt_{\lambda}({M}_{\ell}) + \sum^{n}\wt_{\lambda}(M_{i}') & < & \left( 1 + \frac{g-1}{\NA+1} \right) m^2 + \frac{1}{\NA+1} \sum^{n} m_{i}' - \frac{g-1}{\NA+1} m \nonumber \\
\label{numcrit} & = &  \left( 1+ \frac{g-1+\gamma b }{\NA+1} \right) m^2 - \frac{g-1}{\NA+1} m.
\end{eqnarray}
\end{enumerate}
\end{lemma}

{\it Proof.}  Suppose that we have the required collection of monomials satisfying 
\[\displaystyle{\sum_{\ell=1}^{P(m)}} \wt_{\lambda}({M}_{\ell}) + \sum^{n}\wt_{\lambda}(M_{i}') < 0 \]
with the weights $s_{j}$.  Let $w_{0},\ldots,w_{N}$ be a basis of $H^{0}(C,\mathcal{O}(1))$ diagonalizing the $\lambda$ action.  If $M_{\ell} = w_{0}^{f_{\ell,0}} \cdots w_{N}^{f_{\ell,N}}$, then $\wt_{\lambda}({M}_{\ell}) = \sum_{j=0}^{N} f_{\ell,j}s_{j}$.  

Let $j(i)$ be the function whose value for each $i = 1,\ldots,n$ is the largest index (hence giving the smallest weight) such that the section $w_{j(i)}$ does not vanish at $P_{i}$. Then $\wt_{\lambda}(M_{i}') = m_{i}' s_{j(i)}$.  Thus condition $3.$ may be rewritten
\begin{displaymath}
 \sum_{\ell=1}^{P(m)} \sum_{j=0}^{N} f_{\ell,j} s_{j} + \sum_{i=1}^{n} m_{i}' s_{j(i)} < 0
\end{displaymath}
\begin{displaymath} \\
\Leftrightarrow \sum_{\ell=1}^{P(m)} \sum_{j=0}^{N} f_{\ell,N-j} ((N+1) |s_{0}| r_{j}+s_{0}) + \sum_{i=1}^{n} m_{i}' ((N+1) |s_{0}| r_{N-j(i)}+s_{0}) <  0.
\end{displaymath}
  We proceed to divide by $|s_{0}|$.  Note that our conventions imply that $s_{0}<0$:
\begin{eqnarray}
\sum_{\ell=1}^{P(m)} \sum_{j=0}^{N} f_{\ell,N-j} ((N+1)  r_{j} -1 ) + \sum_{i=1}^{n} m_{i}' ((N+1) r_{N-j(i)} -1) < 0 \nonumber &&\\
\Leftrightarrow (N+1)\sum_{\ell=1}^{P(m)} \sum_{j=0}^{N} f_{\ell,N-j} r_{j} - \sum_{\ell=1}^{P(m)} \sum_{j=0}^{N} f_{\ell,N-j}  + (N+1) \sum_{i=1}^{n} m_{i}' r_{N-j(i)} -  \sum_{i=1}^{n} m_{i}'  <  0  &&\nonumber \\
\Leftrightarrow \sum_{\ell=1}^{P(m)} \sum_{j=0}^{N} f_{\ell,N-j} r_{j} + (N+1) \sum_{i=1}^{n} m_{i}' r_{N-j(i)}  <  \frac{1}{N+1} ( \sum_{\ell=1}^{P(m)} \sum_{j=0}^{N} f_{\ell,N-j}  +   \sum_{i=1}^{n} m_{i}') \nonumber &&
\end{eqnarray}
But we have $ \sum_{j=0}^{N} f_{\ell,N-j} = m$ since each $M_{\ell}$ is a monomial of degree $m$.  Hence we obtain 
\begin{eqnarray}  \sum_{\ell=1}^{P(m)} \sum_{j=0}^{N} f_{\ell,N-j} r_{j} + (N+1) \sum_{i=1}^{n} m_{i}' r_{N-j(i)} & < & \frac{dm-g+1  +   \sum_{i=1}^{n} m_{i}'}{N+1} \nonumber
\end{eqnarray}
Finally, we apply the relation $m_{i} = \gamma b_{i} m^{2}$ associated to the linearization and use $b = \sum b_{i}$:
\begin{eqnarray} \label{to interpret} \sum_{\ell=1}^{P(m)} \sum_{j=0}^{N} f_{\ell,N-j} r_{j} + (N+1) \sum_{i=1}^{n} m_{i}' r_{N-j(i)} & < & \frac{dm-g+1  +   \gamma b m^{2}}{N+1}
\end{eqnarray}

Now, if we let $v_{j} = w_{N-j}$, then the term $ \sum_{j=0}^{N} f_{\ell,N-j} r_{j}$ is the weight of the monomial $v_{0}^{f_{\ell,0}} \cdots v_{N}^{f_{\ell,N}}$.  Also,  $v_{N-j(i)}$ is the smallest weight section among the $v_{j}$'s which does not vanish at $P_{i}$.  Thus we may interpret the left hand side of (\ref{to interpret}) as: the $r$-weight of a collection of monomials restricting to the basis of $H^{0}(C, \mathcal{O}(m))$ plus the $r$-weight of a collection of degree $m_{i}'$ monomials which do not vanish at $P_{i}$.  

This argument can be run in reverse, so given a collection of monomials satisfying $3.'$ we can produce a collection of monomials satisfying $3$.
$\hfill \Box$

Note that property $1.$ above requires a set of \emph{monomials} in $H^{0}(\Pro(V),\mathcal{O}(m))$ which map to a basis of $H^{0}(C,\mathcal{O}(m))$ of small weight.  We want to turn things around, and instead start on the curve in $H^{0}(C,\mathcal{O}(m))$ and work our way back to $H^{0}(\Pro(V),\mathcal{O}(m))$.  The action of a 1-PS $\lambda$ of $SL(V)$ on the Hilbert point of a curve induces a weights on elements of $H^{0}(C, \mathcal{O}_{C}(m))$ (cf. \cite{HM} p. 208).  Briefly, take a basis of $H^{0}(\Pro^{N}, \mathcal{O}(1))$ diagonalizing the $\lambda$ action.  There is an obvious way to define the weight of any degree $m$ monomial, the weight of any degree $m$ homogeneous polynomial is defined to be the maximum weight of its constituent monomials,  and the weight of an element of $H^{0}(C, \mathcal{O}_{C}(m))$ is the minimum of the weights of its preimages in $H^{0}(\Pro^{N}, \mathcal{O}(m))$.

The next proposition says that to establish GIT stability, it is enough to show that there exists \emph{any} basis of $H^{0}(C,\mathcal{O}(m))$ of small weight.

\begin{lemma} \label{curve basis lemma} If there exist a basis of $H^{0}(C,\mathcal{O}(m))$ of $\lambda$-weight $W$, and monomials $M_{1}',\ldots,M_{n}'$ satisfying condition (2) above, and together these satisfy
\begin{displaymath} W + \sum \wt_{\lambda}M_{i}'  \leq  \left( 1+ \frac{g-1+\gamma b }{\NA+1} \right) m^2 - \frac{g-1}{\NA+1} m,
\end{displaymath}
then there are monomials $M_{{1}},\ldots , M_{P(m)}$  which together with $M_{1}',\ldots,M_{n}'$ satisfy conditions 1, 2, and 3' of the numerical criterion.  
\end{lemma}
{\it Proof.}  Let $q_{1},\ldots,q_{P(m)}$ be a basis of $H^{0}(C,\mathcal{O}(m))$satisfying 
\begin{displaymath} W + \sum \wt_{\lambda}M_{i}'  \leq  \left( 1+ \frac{g-1+\gamma b }{\NA+1} \right) m^2 - \frac{g-1}{\NA+1} m.
\end{displaymath}
We may assume that the $q$'s are in order of decreasing weight.  Let $p_{1},\ldots,p_{P(m)}$ be a set of preimages of the $q$'s of minimal weight (that is, $\wt p_{i} = \wt q_{i}$ for each $i$).  Let $\{ M_{i,j} \}$ be the monomials constituting $p_{i}$, so that $p_{i} = \sum_{j=1}^{j_{i}} \alpha_{i,j}M_{i,j}$.  

Write the list of monomials $\{ M_{i,j} \}$ in order of decreasing weight.  If there are ties, choose any order on the tied entries.  Write $y = \# \{ M_{i,j} \}$.  Form the $(P(m)  \! \times \!  y)$-matrix whose entry in row $i$ and the column labelled by $M_{i,j}$ is the coefficient of $M_{i,j}$ in $p_{i}$.  Each row has a leading monomial (the monomial corresponding to the leftmost column with a nonzero entry in that row).  Row reduce this matrix to upper triangular form; this can only lower the leading weight in each row.  Now choose the leading monomials in each row.  Either these map to a basis of $H^{0}(C,\mathcal{O}(m))$ having weight less than or equal to the weight of the basis given by $q_{1}, \ldots, q_{P(m)}$, or else there is a relation between these terms after restriction to the curve.  If this happens, delete the column corresponding to the leftmost monomial appearing in the relation, and begin again (row reduce to upper triangular form, check whether the leading terms in each row give a basis...).  Eventually we must arrive at a set of monomials which give a basis for $H^{0}(C,\mathcal{O}(m))$ (since $\{ \rho(p_{i}) \}$ is a basis of $H^{0}(C,\mathcal{O}(m))$) and the weight of this set of monomials is less than or equal to the weight of the basis given by $q_{1}, \ldots, q_{P(m)}$.   $\hfill \Box$

\subsection{Generalities on profiles} \label{profilessection}
As mentioned in the introduction, the main tool for computing the weight of a basis is something I call a {\it profile}.  (Gieseker uses profiles in his proof, but he doesn't use the word ``profile.'')  We define this abstractly now.  

Let $V$ be a vector space such that every element of $V$ has a weight associated to it.  Let ${F}_{\bullet}$ be a decreasing weighted filtration on ${W}$.  That is, $V = F_{0} \supset F_{1} \supset \cdots \supset F_{N} = 0$, and there is a (finite) decreasing sequence of weights $r_{0} > r_{1} > \cdots > r_{N} = 0 $ such that all the elements of $F_{h}$ have weight less than or equal to $r_{h}$.  

\begin{definition}  The \emph{profile} of a decreasing weighted filtration $F_{\bullet}$ as described above is the graph of the decreasing step function in the $(\mbox{codimension} \times \mbox{weight})$-plane whose value is $r_{h}$ over the interval $[\codim {F}_{h}, \codim {F}_{h+1} ).$
\end{definition}

This is like a distribution function bounding how many linearly independent elements have at most a given weight.  Indeed, given a profile, it is possible to choose a basis whose weight is no greater than the area under the profile.  We will sometimes speak of the ``weight of a filtration'' or ``weight of a profile''; of course what we mean by this is the area underneath the profile, which is a bound for the weight of a basis adapted to this filtration.

Now, there is a notion of an {\it absolute weight filtration}.  It may be described as follows: For each possible weight ${r}_{h}$, form 
\begin{displaymath} {\Omega}({r}_{h}) := \Span \{ v : v \in V, \wt(v) \leq {r}_{h} \}.
\end{displaymath}
Then the profile associated to $\Omega_{\bullet}$ can be used to choose a basis of minimum weight, as it tells exactly how many elements of high weight must be added to the basis before elements of lower weight may be added.

In this paper, we will encounter filtrations of $H^{0}(C,\mathcal{O}(1))$ and $H^{0}(C,\mathcal{O}(m))$.  To help keep track of the ambient vector space of the filtration, we will use tildes for filtrations of $H^{0}(C,\mathcal{O}(m))$.  The filtration of greatest importance for us, $\tilde{X}_{\bullet}$ (to be defined in Section \ref{EWP}), is of this type.

\section{A review of Gieseker's proof}
\label{Giesreview}
Let us quickly review Gieseker's proof from \cite{G}, viewing it as the $n=0$ case of the above setup.  We have recast the numerical criterion to say: the $m$-th Hilbert point of a smooth curve is GIT stable if and only if there exists a basis of $H^{0}(C, \mathcal{O}_{C}(m))$ such that the sum of its weights is less than $(1+\epsilon) m^2$.   

As discussed before Lemma \ref{curve basis lemma}, the action of a 1-PS $\lambda$ of $SL(N+1)$ on the Hilbert point of a curve induces a weights on elements of $H^{0}(C, \mathcal{O}_{C}(m))$ (cf. \cite{HM} p. 208).  Now, it is probably most natural to consider the absolute weight filtration on $H^{0}(C,\mathcal{O}(m))$.  If one could compute its profile, then one could compute Mumford's function $\mu^{L}(x,\lambda)$ on the nose.  However, this is too difficult to compute, so Gieseker considers another filtration instead.

Here is a brief and slightly simplified description of the weighted filtration $\tilde{G}_{\bullet}$ Gieseker uses and its profile.  Given: a curve and a 1-PS $\lambda$.  As before, renormalize the $\lambda$-weights so that they are decreasing and sum to 1.  Let $\{ w_{i} \}$ be a basis of $H^{0}(C,\mathcal{O}_{C}(1)) \cong H^{0}(\Pro^N, \mathcal{O}(1))$ diagonalizing the $\lambda$ action (and compatible with the order of the $r_{i}$).  Let $V_{i} := \mbox{span} ( \{ w_{j} | j \geq i \} ) \subseteq V$.  The normalization ensures that all the points $(im,r_{i}m)$ lie in the first quadrant.  Form the lower envelope of these points, and let $0 = i_{0}, i_{1}, \ldots, $ index the subsequence of points lying on the lower envelope.  Then in $H^{0}(\Pro(V), \mathcal{O}_{\Pro(V)}(m)) \cong \Sym^{m}V$ we have the following filtration:
\begin{equation}
\begin{array}{ccccccccccc}
\Sym^{m} V  = 
V_{i_{0}}^{m}V_{i_{1}}^{0} &  \supset &  V_{i_{0}}^{m-1}V_{i_{1}}^{1} &  \supset & \cdots & \supset &  V_{i_{0}}^{m-p}V_{i_{1}}^{p} & \supset &  \cdots & \supset & V_{i_{0}}^{0}V_{i_{1}}^{m} \\
V_{i_{1}}^{m}V_{i_{2}}^{0} & \supset &  V_{i_{1}}^{m-1}V_{i_{2}}^{1} & \supset & \cdots & \supset & V_{i_{1}}^{m-p}V_{i_{2}}^{p} & \supset &  \cdots & \supset& V_{i_{1}}^{0}V_{i_{2}}^{m} \\
\mbox{etc.}
\end{array}
\end{equation}

The image of this filtration under restriction to the curve gives a filtration $\tilde{G}_{\bullet}$ of $H^{0}(C,\mathcal{O}_{C}(m))$.  We can compute the dimension of each stage of the filtration in $H^{0}(C,\mathcal{O}_{C}(m))$, and we know the weight of each stage, so this is the data of a profile.  The profile is the graph of a step function; its left endpoints lie on the lower envelope of the set of points $\{ (im,r_{i}m ) \}$.   Here is a picture:

\hspace{1in}
\setlength{\unitlength}{18pt}
\thicklines
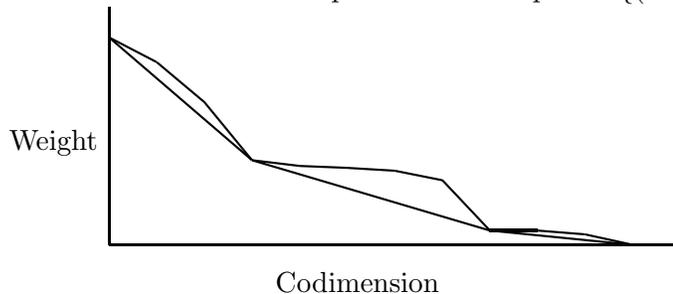
\begin{picture}(11,5)
\put (0,4.36){\line(25,-13){1}}
\put (1,3.84){\line(25,-21){1}}
\put (2,3){\line(50,-61){1}}
\put (3,1.78){\line(25,-3){1}}
\put (4,1.66){\line(25,-1){1}}
\put (5,1.62){\line(50,-3){1}}
\put (6,1.56){\line(5,-1){1}}
\put (7,1.36){\line(50,-53){1}}
\put (8,.3){\line(1,0){1}}
\put (9,.3){\line(25,-2){1}}
\put (10,.22){\line(50,-11){1}}
\put (0,4.36){\line(50,-43){3}}
\put (3,1.78){\line(125,-37){5}}
\put (8,.3){\line(10,-1){3}}
\put (0,0){\line(1,0){12}}
\put (0,0){\line(0,1){5}}
\put (3.5,-1){Codimension}
\put (-2.1,2){Weight}	
\end{picture}
\vspace{0.5in}

\noindent (Looking ahead, the lower envelope here is the inspiration for what I will later call the {\it virtual profile}.)

Any basis adapted to this filtration will establish stability, as the area $A$ under the profile is very close to the area under the lower envelope, and the area under the lower envelope is less than $1 m^{2}$, by a combinatorial lemma due to Morrison (\cite{Morr}, Section 4).

\subsection{The weighted filtration on $H^{0}(C,\mathcal{O}(1))$}
For speed, the previous subsection oversimplified some details of Gieseker's proof.   We will now take the opportunity to begin building up the definitions and notation we need; I have grouped these in this section with his proof, because most of the ideas here are extracted from his proof or follow easily from it.

As we have observed already, the action of the 1-PS $\lambda$ induces most fundamentally a weighted filtration on $H^{0}(C,\mathcal{O}(1))$, but to establish stability we need to find a basis of $H^{0}(C,\mathcal{O}(m))$ of small weight.  We will be going back and forth between these two vector spaces for the rest of the proof.  We begin with $H^{0}(C,\mathcal{O}(1))$, and see what our knowledge of this filtration tells us about filtrations on $H^{0}(C,\mathcal{O}(m))$.  Once we find formulas for the area under the profile for a certain filtration on $H^{0}(C,\mathcal{O}(m))$, we will ultimately bound the weight of the basis by relating quantities back to their counterparts in $H^{0}(C,\mathcal{O}(1))$.

Let $V_{\bullet}$ be the weighted filtration on $H^{0}(C,\mathcal{O}(1))$ induced by the action of the 1-PS $\lambda$.  That is, the stages of the filtration are distinguished by decreasing weight.  Let $z_{j}$ be the size of the $j^{th}$ stage of the filtration, so $z_{j} = \codim V_{j+1} - \codim V_{j}$, and let $r_{j}$ be the weight.  Assume that the weights $r_{j}$ have been normalized so that they are decreasing to zero and sum to 1 (that is, $r_{N} = 0$ and $\sum z_{j} r_{j} = 1$).  Let $D_{j}$ be the base locus of the sublinear series $V_{j}$, and let $d_{j} = \deg D_{j}$.  Let $Q_{1}, ..., Q_{q}$ be the points in $\mbox{Supp} D_{N}$.  (There will be a natural way to order them, but the order is immaterial.)  The marked points $P_{i}$ may or may not show up among the $Q_{i}$; set 
\begin{equation} B_{i} = \left\{ \begin{array}{ll} b_{k}, & Q_{i} = P_{k} \mbox{ for some $k$} \\
0, & Q_{i} \neq P_{k} \mbox{ for any $k$.} \end{array} \right.
\end{equation}
(Note I am already assuming that the marked points are distinct, so $Q_{i}$ can only equal $P_{k}$ for at most one $k$.)  Let $c_{j,i}$ be the multiplicity of $Q_{i}$ in $D_{j}$.  (Note that the indices are not in alphabetic order, opposite the usual convention.  The reasons I have made this choice are too silly to discuss.)   In general  $V_{j}$ is contained in but not equal to $H^{0}(C,\mathcal{O}(1)(-D_{j}))$.  My experience with this problem leads me to conjecture that the maximum of Mumford's $\mu^{L}(x,\lambda)$ function occurs for 1-PS where equality holds at every stage.

\subsection{Relating codegrees and codimensions in $H^{0}(C,\mathcal{O}(1))$}
We have one obvious bound on the weights: $\sum z_{j} r_{j} =1$.  We will need to relate codegrees $d_{j} =  \sum_{i=1}^{n}c_{j,i}$ and codimensions $\sum_{\tau=0}^{j-1}z_{\tau}$.

Near the top of the weighted filtrations, the base loci have low degree, so $\mathcal{O}(1)(-D_{j})$ has high degree, and the dimension/codimension of $H^{0}(C,\mathcal{O}(1)(-D_{j}))$ may be computed using Riemann-Roch.  More precisely: if $\deg D_{j} > \dA-2g+1$, then $\codim V_{j} > \NA - g$.  So if  $\codim V_{j} \leq \NA - g$, then $\deg D_{j} \leq \dA-2g+1$, so $\deg \mathcal{O}(1)(-D_{j}) > 2g-2 $, so $h^{1}( \mathcal{O}(1)(-D_{j}) ) = 0$.  Since $V_{j} \subseteq H^{0}(C,\mathcal{O}(1)(-D_{j}))$, we get a bound: the codegree of $\mathcal{O}(1)(-D_{j})$ cannot exceed the codimension of $V_{j}$.  Recall from the definition of the $z_{j}$'s that $\codim V_{j} = \sum_{\tau=0}^{j-1} z_{\tau}$.  Writing $D_{j} =  \sum_{i=1}^{q} c_{j,i}Q_{i} $, we have: $\deg D_{j} = \sum_{i=1}^{q} c_{j,i}$.  We thus obtain: 
\begin{equation} \label{RR region} \mbox{if $\sum_{\tau=0}^{j-1} z_{\tau} \leq \NA - g$, then $\sum_{i=1}^{q} c_{j,i} \leq \sum_{\tau=0}^{j-1} z_{\tau}$. }
\end{equation}
I call this the {\it Riemann-Roch region} of the filtration.  Write $j_{\RR}$ for the largest index $j$ which satisfies  $\sum_{\tau=0}^{j-1} z_{\tau} \leq \NA - g$.

On the other hand, if $\mathcal{O}(1)$ itself is special, or for stages of the filtration of high codimension (that is, near the bottom), the line bundles $\mathcal{O}(-D_{j})$ have low degree, and we might have $h^{1}(\mathcal{O}(-D_{j})) \neq 0$.  Here we can use Clifford's Theorem to get the following bound:   
\begin{equation} \label{Clifford region} \mbox{if $\sum_{\tau=0}^{j-1} z_{\tau} > \NA - g$, then $\sum_{i=1}^{q} c_{j,i} \leq \sum_{\tau=0}^{j-1} z_{\tau} + \left(\sum_{\tau=0}^{j-1} z_{\tau} - (\NA -g) \right) - h^{1}(C,\mathcal{O}(1))$. }
\end{equation}
I call this the {\it Clifford region} of the filtration and write $j_{\Cliff}$ for the smallest index $j$ which satisfies $\sum_{\tau=0}^{j-1} z_{\tau} > \NA - g$.  (So of course $j_{\Cliff} = j_{\RR}+1$.)

Note that in the case of principal interest (when $d = \nu(2g-2+a)$ and $\nu$ is large, so that $N$ is also large), the Riemann-Roch region accounts for the lion's share of the filtration.

\subsection{Passing to $H^{0}(C,\mathcal{O}(m))$}

We want to use the base loci $D_{j}$ to control how multiples of the $V_{j}$ intersect, and this would work best if  $V_{j} = H^{0}(C,\mathcal{O}(1)(-D_{j}))$.   Gieseker observed that if we pass from $H^{0}(C,\mathcal{O}(1))$ to $H^{0}(C,\mathcal{O}(m))$ (which is where we ultimately need to produce a basis anyway), then we will be able to treat an arbitrary 1-PS $\lambda$ as if it were of this form.  Most of the proof of Lemma \ref{multup lemma} below comes from pages 54--55 of \cite{GCime}.  However, I want to add a few comments to Gieseker's proof, so I will run through the argument here.

Let $(V_{s}^{u-w}V_{t}^{w}V_{0})^{v}$ denote the subspace of $H^{0}(C,\mathcal{O}((u+1)v))$ generated by expressions of the form $x_{1} \cdots x_{v(u-w)} y_{1} \cdots y_{vw} z_{1} \cdots z_{v}$ where the $x$'s come from $V_{s}$, the $y$'s come from $V_{t}$, and the $z$'s come from $V_{0}$.  

\begin{lemma} \label{multup lemma} Let $u,v,w$ be nonnegative integers with $0 \leq w \leq u$ and $v \geq 1$.  Suppose $C$ is an arbitrary subscheme of $\Pro^{N}$ with Hilbert polynomial $dt-g+1$ and 
\[ v \geq \frac{d^2 (u+1)^2 - d(u+1)}{2} -g+1.
\]
Then 
\begin{displaymath} (V_{s}^{u-w}V_{t}^{w}V_{0})^{v} = H^{0}(C, \mathcal{O}((u+1)v)(-(u-w)D_{s} - wD_{t}))
\end{displaymath}
\end{lemma}

{\it Remark.}   Note that the bound on $v$ depends on $u$ and the Hilbert polynomial $P(z) = \dA z-g+1$, but not on the curve $C$ or the line bundle $\mathcal{O}_{C}(1)$ embedding $C$ into $\Pro^{N}$.    

{\it Proof.}   Let $L_{s}$ and  $L_{t}$ be the line bundles generated by the sections in $V_{s}$ and $V_{t}$.  Here is the first comment to add to Gieseker's proof: then $L_{s} = \mathcal{O}_{C}(1)(-D_{s})$.  We have 
\[ (V_{s}^{u-w} V_{t}^{w} V_{0})^{v} \subset H^{0}(C,(L_{s}^{u-w} L_{t}^{w} L_{0})^{v})
= H^{0}(C,\mathcal{O}((u+1)v)(-(u-w)D_{s} - wD_{t})). \]

Now, since sections in $V_{s}^{u-w} V_{t}^{w}$ generate $L_{s}^{u-w} L_{t}^{w}$, and  $V_{0}$ is very ample, we have that $V_{s}^{u-w} V_{t}^{w} V_{0}$ is very ample, and hence determines an embedding $C \hookrightarrow \Pro^{M}$.  We have a short exact sequence
\[ 0 \rightarrow \mathcal{I}(v) \rightarrow \mathcal{O}_{\Pro^{M}}(v) \rightarrow \mathcal{O}_{C}(v) \rightarrow 0.
\]
(We now have two $\mathcal{O}_{C}(1)$'s in this proof, corresponding to the embeddings in $\Pro^{N}$ and $\Pro^{M}$, but it is not difficult to tell them apart.)  Write $d_{s} = \deg D_{s}$, respectively for $t$; then $\deg L_{s} = d-d_{s}$ and $ \deg L_{t} = d-d_{t}$.  Then the Hilbert polynomial for $C \subset \Pro^{M}$ is  
\[ P(z) = ( (d-d_{s})(u-w) + (d-d_{t})(w) + d )z -g+1.
\]
The Gotzmann number for this Hilbert polynomial is 
\[ m_{0} = \frac{( (d-d_{s})(u-w) + (d-d_{t})(w) + d )^{2} -  ((d-d_{s})(u-w) + (d-d_{t})(w) + d) }{2} -g+1;
\]
recall that the Gotzmann number for a Hilbert polynomial has the property that it is the maximum regularity for any sheaf with that Hilbert polynomial (\cite{Gotz} Lemma 2.9).  Hence, $H^{1}(\mathcal{I}(v))=0$ since $v$ is larger than the Gotzmann number.  But then  
\[ H^{0} ( \Pro^{M}, \mathcal{O}(v)) \rightarrow H^{0}(C,(L_{s}^{u-w} L_{t}^{w} L_{0})^{v}
\]
is surjective.

Comparing this to the definition of $(V_{s}^{u-w} V_{t}^{w} V_{0})^{v}$, this says that 
\[ (V_{s}^{u-w} V_{t}^{w} V_{0})^{v} = H^{0}(C,(L_{s}^{u-w} L_{t}^{w} L_{0})^{v}) = H^{0}(C,\mathcal{O}((u+1)v)(-(u-w)D_{s} - wD_{t}))
\]
as desired.  

Finally note that $d-d_{s}$ and $d-d_{t}$ are no larger than $d$; hence taking
\[ v \geq \frac{d^2 (u+1)^2 - d(u+1)}{2} -g+1.
\]
ensures that $v$ is greater than or equal to the Gotzmann number for any $V_{s}$ and $V_{t}$.   

$\hfill$ $\Box$

{\it Remark.}  We will be applying this result when $C$ is a smooth curve in $\Pro^{N}$; for this application, the Gotzmann number is really much larger than we should need.   I hope to improve this result significantly, which should be helpful (if not necessary) when studying stability for small values of $m$.  

\vspace{.25in}
Let $m = (u+1)v$.  Then there is a filtration $\tilde{V}_{\bullet}$ of $H^{0}(C,\mathcal{O}(m))$ by the subspaces $(V_{j}^{u}V_{0})^{v}$.

Note however that if in the original filtration, there are two successive stages where the base locus does not increase, now, after passing to $H^{0}(C,\mathcal{O}(m))$, the second of these stages has risen up to replace the first of these two stages.  Thus, in $H^{0}(C,\mathcal{O}(m))$, we need only record the subsequence of the $j$'s where the degree of the base locus increases.  I will index these by the letter $k$.

The filtration $\tilde{V}_{\bullet}$ may be further refined by using spaces of the form $(V_{k}^{u-w}V_{k+1})^{w}V_{0})^{v}$.  We will abuse notation and write $\tilde{V}_{\bullet}$ for this refinement also.   Thus, the index of the filtration $\tilde{V}_{\bullet}$ may be the single index $k$, or a pair $(k,w)$.

I will use tildes for quantities associated to $\tilde{V}_{\bullet}$. We have $\tilde{V}_{k} = H^{0}(C,\mathcal{O}(m)(-\tilde{D}_{k}))$, where $\tilde{D}_{k} = uv D_{j_{k}}$.  We write $\tilde{d}_{k}:= uv d_{j_{k}} $ and $\tilde{c}_{k,i} := uv c_{j_{k},i}$.  Then 
\[ \tilde{V}_{k} = H^{0}(C,\mathcal{O}(m)(-\tilde{c}_{k,1}Q_{1} - \cdots - \tilde{c}_{k,q}Q_{q} ))
\]
and elements of this space have weight $\leq \tilde{r}_{k} := uv r_{j_{k}} +vr_{0}$.

 Define $\tilde{N}$ to be the smallest index giving the $vr_{0}$-weight space.  We have: 
\begin{equation} \label{multupfilt}
\begin{array}{lc} \mbox{Space} &  \mbox{Weight} \\ [.0625in]
\tilde{V}_{0}  =  H^{0}(C,\mathcal{O}(m)) &  \tilde{r}_{0} \\ [.0625in]
\tilde{V}_{1}  =  H^{0}(C,\mathcal{O}(m)(-\tilde{c}_{1,1}Q_{1}
- \cdots - \tilde{c}_{1,q}Q_{q})) &  \tilde{r}_{1} \\ [.0625in]
\tilde{V}_{2}  =  H^{0}(C,\mathcal{O}(m)(-\tilde{c}_{2,1}Q_{1}
 - \cdots - \tilde{c}_{2,q}Q_{q})) &   \tilde{r}_{2} \\ [.0625in]
 \vdots &   \vdots \\ [.0625in]  
\tilde{V}_{\tilde{N}}  =  H^{0}(C,\mathcal{O}(m)(-\tilde{c}_{\tilde{N},1}Q_{1} - \cdots - \tilde{c}_{\tilde{N},q}Q_{q}))  &   \tilde{r}_{\tilde{N}} = vr_{0} \\ [.0625in]
\end{array}
\end{equation}

We may extract the multiplicities of the points in the base loci in the weighted filtration $\tilde{V}_{\bullet}$ and the weights to obtain an $(\tilde{N}+1) \times (q+1)$ array:

\begin{equation} \label{M}
\left( \begin{array}{cccc} 
\tilde{c}_{0,1} & \cdots & \tilde{c}_{0,q} & \tilde{r}_{0} \\
\tilde{c}_{1,1} & \cdots & \tilde{c}_{1,q} & \tilde{r}_{1} \\
\vdots         & \vdots & \vdots & \vdots \\
\tilde{c}_{\tilde{N},1} & \cdots & \tilde{c}_{\tilde{N},q} & \tilde{r}_{\tilde{N}} = vr_{0}
\end{array} \right)
\end{equation}

This array has the following properties:  the $\tilde{c}_{k,i}$'s are all nonnegative integers; the $\tilde{r}_{i}$'s are rational numbers weakly decreasing to $vr_{0}$; and in the first row the $\tilde{c}_{0,i}$'s are all zero.  Furthermore we see that the sum of the entries in row $k$ is governed by either a Riemann-Roch bound (\ref{RR region}) or a Clifford bound (\ref{Clifford region}).

\section{Why Gieseker's proof doesn't cover marked points}
\label{keyexample}

To my knowledge, Elizabeth Baldwin first wrote down the straightforward generalization of Gieseker's result to $\Mgn$ (unpublished), and it is not difficult to see that the analogue of Gieseker's filtration does not suffice to establish stability in cases where $b_{i}$ is more than a little larger than 0.  Here is a counterexample:

\subsection{Example 1} \label{lindecnoskips}
{\it Purpose:} to show that the profile associated to $\tilde{G}_{\bullet}$ (which equals $\tilde{V}_{\bullet}$ in this example) does not suffice to establish asymptotic Hilbert stability when there are marked points.

Suppose $n \geq 3$.  Consider the 1-PS $\lambda$ which acts with linearly decreasing weights on the marked points.  That is, $\lambda$ induces the following weighted filtration:

\begin{displaymath}
\begin{array}{lclc}\mbox{Space} & & & \mbox{Weight} \\ [.0625in]
 V_{0} & = & H^{0}(C,\mathcal{O}(1)) & \hspace{.125in} \frac{1}{2} \\ [.0625in]
V_{1} & = & H^{0}(C,\mathcal{O}(1)(-P_{1})) &  \hspace{.125in} \frac{1}{3} \\ [.0625in]
V_{2} & = & H^{0}(C,\mathcal{O}(1)(-P_{1}-P_{2})) &  \hspace{.125in} \frac{1}{6} \\ [.0625in]
V_{3} & = & H^{0}(C,\mathcal{O}(1)(-P_{1}-P_{2}-P_{3})) &  \hspace{.125in} 0
\end{array}
\end{displaymath}  

The points $(im,r_{i}m)$ all lie on their lower envelope.  Also, we have $r_{0} + r_{1} + r_{2} = 1$.  Using $\gamma b_{i} = 1/2$, we have $T \approx 1m^2 - \frac{1}{4}m^2 + \gamma b m^2 = 5/4m^2 > (1+ \epsilon) m^2$.  

So the straightforward adaptation of Gieseker's proof is not enough to establish the stability of smooth pointed curves with respect to the linearizations we have specified.

\subsection{The key observation} \label{keyobservation}

In fact it is not difficult to show that the 1-PS of Example 1 is not destabilizing.  

We use the following easy linear algebra lemma:

\begin{lemma} Let $V_{1},\ldots, V_{n}$ be subspaces of a vector space $V$.  Write $V_{ij} := V_{i} \cap V_{j}$, $V_{ijk} := V_{i} \cap V_{j} \cap V_{k}$, etc.  Then 

\begin{displaymath} \codim \,  \Span \{V_{1},\ldots,V_{n} \} = \sum \codim V_{i} - \sum_{i<j} \codim V_{ij}+ \sum_{i<j<k} \codim V_{ijk} - \cdots + (-1)^{n-1} \codim V_{123 \cdots n}.
\end{displaymath}
\end{lemma}

Gieseker's proof proceeds as follows: $H^{0}(C,\mathcal{O}(m))$ contains the following spaces with the following codimensions and weights:

\begin{displaymath}
\begin{array}{ccl} \mbox{weight} & \mbox{codimension} & \mbox{space} \\
\frac{1}{2}m & 0 &  H^{0}(C,\mathcal{O}(m)) \\
\frac{1}{2}m -\frac{1}{6} & 1 & H^{0}(C,\mathcal{O}(m)(-P_{1})) \\
\frac{1}{2}m -\frac{2}{6} & 2 & H^{0}(C,\mathcal{O}(m)(-2P_{1})) \\
\frac{1}{2}m -\frac{3}{6}& 3 & H^{0}(C,\mathcal{O}(m)(-3P_{1})) \\
\frac{1}{2}m -\frac{4}{6}& 4 & H^{0}(C,\mathcal{O}(m)(-4P_{1})) \\
\vdots & \vdots & \vdots
\end{array}
\end{displaymath} 

As discussed above, if one basis element is chosen from each of these spaces, then $A$ is approximately $1-\frac{1}{n+1} = \frac{3}{4}$.  

The key observation is that we know more subspaces corresponding to each weight in the left column.  For instance, elements of the spaces $ H^{0}(C,\mathcal{O}(m)(-2P_{1}))$ and $ H^{0}(C,\mathcal{O}(m)(-P_{1}-P_{2}))$ each have weight $\frac{1}{2}m -\frac{2}{6}$.  These spaces each have codimension 2, and their intersection  $ H^{0}(C,\mathcal{O}(m)(-2P_{1}-P_{2}))$ has codimension 3, so their span has codimension $2 + 2 - 3 = 1$.

Now we try again to choose a basis of $H^{0}(C,\mathcal{O}(m))$ of lowest weight.  For the first basis element, we may be obliged to choose one element of top weight $ \frac{1}{2}m$.  But for the second basis element, we now know that we may bypass the elements of weights $\frac{1}{2}m -\frac{1}{6}$ and  $\frac{1}{2}m -\frac{2}{6}$ and instead choose an element of weight $\frac{1}{2}m -\frac{3}{6}$.  My proof repeatedly uses this trick, suggested by Ian Morrison, to establish a filtration and profile giving a basis of lower weight than Gieseker's.

\subsection{Minimizing multiplicities}   

Soon we are going to put a lot of effort into minimizing multiplicities.  The following lemma shows that this makes easy work of computing spans of spaces of the form we have encountered.  

\begin{lemma} \label{trace} Suppose we are given $q$ subspaces $E_{1},\ldots,E_{q}$ of $H^{0}(C,\mathcal{O}(m))$ of the form:
\begin{displaymath}
\begin{array}{lcl}
E_{1} & = & H^{0}(C,\mathcal{O}(m)(-d_{1,1}Q_{1} - \cdots -  d_{1,q}Q_{q}) \\
E_{2} & = & H^{0}(C,\mathcal{O}(m)(-d_{2,1}Q_{1} - \cdots  - d_{2,q}Q_{q})) \\
 &\vdots & \\
E_{q} & = & H^{0}(C,\mathcal{O}(m)(-d_{q,1}Q_{1} - \cdots  - d_{q,q}Q_{q})) 
\end{array}
\end{displaymath}
The $E_{i}$ need not be distinct, and though the notation looks a little similar to that of filtrations above, we do not mean in any way to imply that the $E_{i}$ form a filtration---in the applications we have in mind, they do not.

Suppose that $E_{i}$ minimizes the multiplicity of $Q_{i}$---that is, the minimum in each column appears along the diagonal.  Suppose also that 
\begin{displaymath} \sum_{i=1}^{q} \max_{j} d_{j,i} < dm - 2g.
\end{displaymath}
Then 
\begin{displaymath}\Span (E_{1},\ldots,E_{q}) = H^{0}(C,\mathcal{O}(m)(-\sum_{i=1}^{q} d_{i,i} Q_{i}))
\end{displaymath}
and
\begin{displaymath}
\codim \Span (E_{1},\ldots,E_{q}) = d_{1,1} + d_{2,2} + \cdots + d_{q,q}.  
\end{displaymath}
\end{lemma} 

{\it Proof.}  The condition 
\begin{displaymath} \sum_{i=1}^{q} \max_{j} d_{j,i} Q_i < dm - 2g.
\end{displaymath}
ensures that the codimension of the intersection of any subset of these $q$ spaces may be computed using Riemann-Roch.  Thus, for each subset $I \subseteq \{ 1,\ldots,q \}$, say $I = \{ i_{1},\ldots,i_{k} \} $ we have
\begin{displaymath} \codim E_{i_{1} \cdots i_{k}} = \max(d_{i_{1},1},\ldots,d_{i_{k},1}) + \max(d_{i_{1},2},\ldots, d_{i_{k},2}) + \cdots + \max(d_{i_{1},q},\ldots,d_{i_{k},q}). 
\end{displaymath}  

Suppose $j \not\in I$.  Then the term $\max(d_{i_{1},j},\ldots,d_{i_{k},j})$ is cancelled by a term coming from $I \cup \{ j \} $.  Being a subset of cardinality one greater, its codimension gets opposite sign from that of $I$.  And since by hypothesis $d_{j,j}$ is the smallest term in column $j$, it drops out of $\max(d_{i_{1},j},\ldots,d_{i_{k},j},d_{j,j})$, giving us exactly the cancellation we claimed.  Given $I$, every $j \in \{ 1,\ldots,q \}$ is either in $I$ or not in $I$, so it is clear whether the term $\max(d_{i_{1},j},\ldots,d_{i_{k},j})$ is cancelling or being cancelled.  The only terms surviving are the $d_{i,i}$ since there are no double intersections of the form $E_{ii}$ in our setup to cancel them.  

Finally, the base locus of $\Span (E_{1},\ldots,E_{q})$ must be $\sum_{i=1}^{q} d_{i,i} Q_{i}$ (since we can find sections that vanish to each $Q_{i}$ to exactly order $d_{i,i}$).  This gives   
\begin{equation} \label{spansubset} \Span (E_{1},\ldots,E_{q}) \subset H^{0}(C,\mathcal{O}(m)(-\sum_{i=1}^{q} d_{i,i} Q_{i})).
\end{equation}
But the codimensions of the two spaces in line (\ref{spansubset}) are the same, so we must actually have equality:
\begin{displaymath}\Span (E_{1},\ldots,E_{q}) = H^{0}(C,\mathcal{O}(m)(-\sum_{i=1}^{q} d_{i,i} Q_{i})).
\end{displaymath}
$\hfill \Box $

\section{The filtration $\tilde{X}_{\bullet}$ and its profile}
\label{EWP}

\subsection{Subscript conventions}
\label{subscriptsI}
In the course of the proof we will need to keep track of a set of subsequences of a subsequence of a sequence.   My first attempt, using several layers of subscripts, proved unsatisfactory; I know of no good convention for this kind of accounting, so I will use the following notation and conventions.

\subsubsection{Tildes}
Recall that $k$ indexes a subset of the rows $j$ of the original filtration $V_{\bullet}$.   Quantities associated to $\tilde{V}_{\bullet}$ (like the multiplicities $\tilde{c}$ and weights $\tilde{r}$ are written with tildes and indexed by $k$'s; quantities associated to $V_{\bullet}$ (such as $c$ and $r$) have no tildes and are indexed by $j$'s.  When I want to refer to a subsequence of $c$ or $r$, rather than using nested subscripts and writing for instance $r_{j_{k}}$ I will simply write $r_{k}$; this should cause no confusion, since the presence or absence of a tilde indicates whether a layer has been suppressed.  

\subsubsection{Cases I-IV and the functions $s(k,i)$ and $t(k,i)$} \label{Cases, s, and t}
  It is useful to define two functions $s$ and $t$ in some (but not all) situations.  We will take the time now to define four cases, which will be referred to in this section and in Section \ref{virvsactual}. 
\begin{enumerate}
\item[I.]   We have $\tilde{c}_{k,i} < \tilde{c}_{k+1,i} < \tilde{c}_{k+2,i}$.  That is, the multiplicity of the point $Q_{i}$ jumps at row $k$ and again at row $k+1$.  In this case we do not define $s(k,i)$ and $t(k,i).$
  
\item[II.]  We have $\tilde{c}_{k,i} = \tilde{c}_{k+1,i} = \tilde{c}_{k+2,i}$.  That is, the multiplicity of $Q_{i}$ does not jump at row $k$ or at row $k+1$.   Define $s(k,i)$ to be the last row where this multiplicity jumped, and let $t(k,i)$ be the next row where it jumps, or else $\tilde{N}$ if $\tilde{c}_{k,i} = \tilde{c}_{\tilde{N},i}$.  In symbols, $s(k,i)$ is the largest index strictly (in Case II) less than $k$ such that $\tilde{c}_{s(k,i),i} < \tilde{c}_{s(k,i)+1,i}$, and $t(k,i)$ is the smallest index strictly (in Case II) greater than $k$ such that $\tilde{c}_{t(k,i),i} < \tilde{c}_{t(k,i)+1,i}$ if this exists, or else $\tilde{N}$.  Finally, the reader will see after reading Case III and Case IV that in Case II we have $s(k,i) = s(k+1,i)$ and $t(k,i) = t(k+1,i)$.
 
\item[III.]  We have $\tilde{c}_{k,i} = \tilde{c}_{k+1,i} < \tilde{c}_{k+2,i}$.  That is, the multiplicity of $Q_{i}$ does not jump at row $k$ but jumps at row $k+1$.  Then as in Case II we define $s(k,i)$ to be the last row where this multiplicity jumped, and we define $t(k,i) = k+1$.  
  
\item[IV.] We have $\tilde{c}_{k,i} < \tilde{c}_{k+1,i} = \tilde{c}_{k+2,i}$.  That is, the multiplicity of $Q_{i}$ jumps at row $k$ but not at row $k+1$.  We define $s(k,i) = k$, and as in Case II let $t(k,i)$ be the next row where this multiplicity jumps, or else $\tilde{N}$ if $\tilde{c}_{k,i} = \tilde{c}_{\tilde{N},i}$.
\end{enumerate}

Defining $s$ and $t$ differently in Cases II-IV as we have done permits us to treat these cases simultaneously in Section \ref{Computing bounds for A_k,i}, which more than makes up for the extra work involved here.  There are two reasons why Case I is treated separately from the other cases. First, there is an easy way to deal with Case I that is not available in Cases II-IV.  Second, if one tries to analyze Case I the way we analyze Cases II-IV, one obtains a coefficient which I can bound in Case II-IV which I have not figured out how to bound in Case I.  So, it is desirable to treat Case I separately.

\subsubsection{Eliminating redundancies}
Rather than printing $i$ redundantly in subscripts, whenever I can I will leave it off the second time.  For example I will simply write $\tilde{c}_{s(k,i)}$ for $\tilde{c}_{s(k,i),i}$.

\subsubsection{The functions $j(i,\ell)$ and $k(i,\ell)$}
We will also want to keep track of the subset of $j$'s or $k$'s where the multiplicity of the point $Q_{i}$ in the base locus increases.  I will do this as follows:

Say the multiplicity of $Q_{i}$ jumps $\Ki$ times between the top of the filtration and the bottom.  We start counting from zero, so these stages of the filtration are the $0^{th}$ jump up through the $(\Ki-1)^{th}$ jump.  As a convention, we append $\bar{N}$ (the index of the last row of the filtration $V_{\bullet}$) or $\tilde{N}$ (the index of the last row of the filtration $\tilde{V}_{\bullet}$) as the $\Ki^{th}$ element of this sequence.  We write two increasing set functions
\[ j(i,\bullet): \{ 0, \ldots, \Ki \} \rightarrow \{ 0, \ldots, \bar{N} \}
\]
and
\[ k(i,\bullet): \{ 0, \ldots, \Ki \} \rightarrow \{ 0, \ldots, \tilde{N} \}
\]
and use these to index the rows where the multiplicity of the point $Q_{i}$ in the base locus increases.  That is, the function $j(i,\bullet)$ takes values in the $j$'s, and similarly $k(i,\bullet)$ takes values in the $k$'s.  Here is an example to give a little practice with this notation:  $j(i,0)$ means the index $j$ where the multiplicity of $Q_{i}$ jumps for the $0^{th}$ time.  This is the lowest row of the filtration where $Q_{i}$ is not in the base locus, so $r_{j(i,0)}$ is the least weight of a section not vanishing at $Q_{i}$.

As before, when $i$ appears more than once in a subscript, we will leave it off the second time.  Thus $c_{j(i,0),i} $ becomes $c_{j(i,0)} $  and we have $c_{j(i,0)} =0 $ while $c_{j(i,0)+1} = c_{j(i,1)} > 0$.

\subsubsection{A consequence of these conventions}
As a consequence, note that previously when going between the filtrations $V_{\bullet}$ and $\tilde{V}_{\bullet}$ we had $\tilde{c}_{k,i} = uv c_{j_{k},i}$.  But now with our new notation we can write $\tilde{c}_{k(i,\ell)} = uv c_{j(i,\ell)}$.  In this sense the definitions of $j(i,\ell)$ and $k(i,\ell)$ have eliminated some of the need for nested subscripts.

Finally we note that although the notations are similar in format, $j$ and $k$ are somewhat different in character from $s$ and $t$.  Briefly, $j$ and $k$ are ``lookup'' functions, whereas $s$ and $t$ are ``previous'' and ``next'' functions.

\subsection{The filtration $\tilde{X}_{\bullet}$ and its profile} \label{xtilde section}
Here we describe the filtration $\tilde{X}_{\bullet}$ of $H^{0}(C,\mathcal{O}(m))$ and its weight profile.  $\tilde{X}_{\bullet}$ is obtained from the filtration  $\tilde{V}_{\bullet}$ by taking spans of the stages of $\tilde{V}_{\bullet}$ with other cleverly chosen spaces.

The filtration $\tilde{X}_{\bullet}$ will have $\tilde{N} \times u +1$ stages.  

For each $k = 0, \ldots, \tilde{N}-1$, and for each $w = 0, \ldots, u-1$ we want to describe the space $\tilde{X}_{k, w}$.  Our starting point is the space $(V_{k}^{u-w} V_{k+1}^{w} V_{0})^{v}$.  Elements of this space have weight less than or equal to $v(u-w) r_{k} + vwr_{k+1} +vr_{0}$.

Our goal: for each $i$ from 1 to $q$, find subspaces of $H^{0}(C,\mathcal{O}(m))$ whose weight is less than or equal to $v(u-w) r_{k} + vwr_{k+1} +vr_{0}$, for which the multiplicity of $Q_{i}$ is less than the multiplicity in the base locus of $(V_{k}^{u-w} V_{k+1}^{w} V_{0})^{v}$.  We do this as described in the following definition.  Also, it is convenient to define certain quantities $\tilde{x}(k,i,w)$ at this time; their role will be explained soon.

\begin{definition}[The filtration $\tilde{X}_{\bullet}$ and its profile] \label{ProfileDefinition}
First, $\tilde{X}_{0,0} = H^{0}(C,\mathcal{O}(m)).$

For the remaining triples $(k,w,i)$ with $(k,w) \neq (0,0)$, where $k = 0, \ldots, \tilde{N}-1$, $w = 0, \ldots, u-1$,  and $i = 1, \ldots, q $, the contribution to the profile is found as follows:
\begin{itemize}
\item If the multiplicity of $Q_{i}$ is zero in row $k+1$ (and hence zero in row $k$ also), there is no contribution to $\tilde{X}_{k,w}$, and $\tilde{x}(k,i,w)=0$.
\item If the multiplicity of $Q_{i}$ is nonzero in row $k+1$ and we are in Case I as defined in Section \ref{Cases, s, and t}, so the multiplicity of $Q_{i}$ jumps at row $k$ and row $k+1$, then we add no new spaces to $\tilde{X}_{k,w}$, and the space $(V_{k}^{u-w} V_{k+1}^{w} V_{0})^{v}$ into $\tilde{X}_{k, w}$, and $\tilde{x}(k,i,w)$ is the multiplicity of $Q_{i}$ in $(V_{k}^{u-w} V_{k+1}^{w} V_{0})^{v}$; 
\item If the multiplicity of $Q_{i}$ is nonzero in row $k+1$ and we are in Case II, III, or IV as defined in Section \ref{Cases, s, and t}, so the multiplicity of $Q_{i}$ jumps at no more than one of the rows $k$ and $k+1$, let $s(k,i)$ and $t(k,i)$ be as defined there.  For each $w$ we find the smallest integer $W = W(u,v;k,w,i)$ such that $(V_{s(k,i)}^{u-W}V_{t(k,i)}^{W}V_{0})^{v}$ has weight less than $v(u-w)r_{k} + vwr_{k+1}+vr_{0}$.  Then $(V_{s(k,i)}^{u-W}V_{t(k,i)}^{W}V_{0})^{v}$ is added to $\tilde{X}_{k,w}$, and $\tilde{x}(k,i,w)$ is the multiplicity of $Q_{i}$ in the base locus of $(V_{s(k,i)}^{u-W}V_{t(k,i)}^{W}V_{0})^{v}$.
\end{itemize}
 Then
\begin{displaymath} \tilde{X}_{k, w} = \Span \{ (V_{k}^{u-w} V_{k+1}^{w} V_{0})^{v}, \mbox{spaces of type $(V_{s(k,i)}^{u-W}V_{t(k,i)}^{W}V_{0})^{v}$ if there are any} \}, 
\end{displaymath}
and let $\tilde{x}(k,w)$ be the codimension of $\tilde{X}_{k,w}$.

Note $\tilde{X}_{k, w}$ is the span of between $1$ and $q+1$ distinct spaces; there may be fewer than $q+1$ distinct spaces in the span, as there may be points $Q_{i}$, which make no contribution, and/or repeats may occur among the spaces of the form $(V_{s(k,i)}^{u-W}V_{t(k,i)}^{W}V_{0})^{v}$.   

Finally, for the last stage of the filtration, define $\tilde{X}_{\tilde{N}} := \tilde{V}_{\tilde{N}}$.

Thus, the profile associated to $\tilde{X}_{\bullet}$ is the graph of decreasing step function whose value over the intervals $[ \tilde{x}(k,w), \tilde{x}(k,w+1) )$ is $v(u-w) r_{k} + vwr_{k+1} +vr_{0}$, and whose value over the interval $[ \codim \tilde{X}_{\tilde{N}}, \dim H^{0}(C,\mathcal{O}(m)) ]$ is $v r_{0}$.
\end{definition}

Note that the spaces used to construct each $\tilde{X}_{k,w}$ satisfy the degree hypothesis of Lemma \ref{trace}: every space going into the span is either of the form $(V_{k}^{u-w}V_{k+1}^{w}V_{0})^{v}$ or $(V_{s(k,i)}^{u-\wkwi}V_{t(k,i)}^{\wkwi}V_{0})^{v}$.  But the base locus of any space of this form is bounded by the base locus of $(V_{\bar{N}}^{u}V_{0})^{v}$, which is $uvc_{\bar{N},1} + \cdots + uv c_{\bar{N},q}.$  That is, $\max_{j} \{ d_{j,i} \} \leq uvc_{\bar{N},i}$, so we have \begin{displaymath} \sum_{i=1}^{q} \max_{j} \{ d_{j,i} \} \leq  \sum_{i=1}^{q}  uvc_{\bar{N},i} \leq uv d < uv d + ud -2g = dm-2g.
\end{displaymath}

However, it is not always true that $(V_{k}^{u-w}V_{k+1}^{w}V_{0})^{v}$ or $(V_{s(k,i)}^{u-W}V_{t(k,i)}^{W}V_{0})^{v}$ always minimizes the multiplicity of $Q_{i}$ among these $q$ spaces.  (It is possible to find the minimum, but we will not do this now.  See Section \ref{lowerenvelopesarebetter} for a little more discussion.)  Therefore, we cannot apply Lemma \ref{trace} to conclude that $\tilde{x}(k,w) = \sum_{i=1}^{q} \tilde{x}(k,w,i)$.  However, we may use Lemma \ref{trace} to conclude that $\tilde{x}(k,w) \leq \sum_{i=1}^{q} \tilde{x}(k,w,i)$, since the minimum multiplicity for the point $Q_{i}$ must be smaller than $\tilde{x}(k,w,i)$.  Of course, this is not enough to bound $ \tilde{x}(k,w+1) -  \tilde{x}(k,w)$.  But since the $\tilde{r}_{k}$'s are decreasing, the weight $A$ of this profile will only decrease if some  $\tilde{x}(k,w) < \sum_{i=1}^{q} \tilde{x}(k,w,i)$.  So computing using equality at every stage gives the following upper bound for $A$:  \begin{equation} \label{firstAbound} A \leq \sum_{k=0}^{\tilde{N}-1} \sum_{w=0}^{u-1} (v(u-w)r_{k} + vwr_{k+1}+vr_{0})  (\tilde{x}(k,w+1) -  \tilde{x}(k,w))+ (\dim \tilde{X}_{\tilde{N}})vr_{0}.
\end{equation}

We have $\tilde{X}_{\tilde{N}} = H^{0}(C,\mathcal{O}(m)(-uv D_{\bar{N}}))$, and so we may compute $\dim \tilde{X}_{\tilde{N}} = dm - uv d_{\bar{N}} -g+1 = (d-d_{\bar{N}})uv + dv -g+1$.  Substituting this into (\ref{firstAbound}), we obtain 
\begin{equation} \label{secondAbound} A \leq \sum_{k=0}^{\tilde{N}-1} \sum_{w=0}^{u-1} (v(u-w)r_{k} + vwr_{k+1}+vr_{0})  (\tilde{x}(k,w+1) -  \tilde{x}(k,w))+((d-d_{\bar{N}})uv + dv -g+1)vr_{0}.
\end{equation}

Rather than trying to bound the right hand side of (\ref{secondAbound}), we will follow a different approach.  We will define a ``virtual'' profile whose graph has area $A^{\vir}$ nearly the same as the area of the graph $A$ of the actual profile, but which is computationally a little easier to work with.   Let $\Delta = A  - A^{\vir}$ be the discrepancy.  Also, for each $i$ between $1$ and $q$, recall that $r_{j(i,0)}$ is the $r_{j}$ such that $c_{j,i} = 0$ and $c_{j+1,i} > 0$.  Then 
\begin{equation} \label{firstTbound} T \leq A^{\vir} + \Delta + \sum_{i=1}^{n} \gamma B_{i} r_{j(i,0)} (u+1)^2 v^2.
\end{equation}

We use the rest of this section to define the virtual profile.  In the next section we bound $\Delta$, and in Section \ref{boundingTvir} we bound $A^{\vir} + \sum_{i=1}^{n} \gamma B_{i} r_{j(i,0)}(u+1)^2 v^2$.  Putting this all together with (\ref{firstTbound}), we will get a bound for $T$.

\subsection{The virtual profile}

\label{virtualProfile}

The virtual profile simplifies the graph of the profile in three ways:
\begin{itemize}
\item In the profile, we form a span of $q$ spaces for all $k$ and for all $w$, so the step function is defined over $ \tilde{N} \times u +1$ intervals; in the virtual profile, we only partition the domain (the codimension axis) into $\tilde{N}+1$ intervals.
\item In the profile, we round so that $W = W(u,v;k,w,i)$ is always an integer, so exponents, multiplicities, and codimensions are integers; in the virtual profile, their counterparts are rational numbers.   
\item In particular the quantity $\tilde{f}(k)$ (defined below) is the virtual counterpart to $\tilde{x}(k,0)$.  The profile is a step function, so the two points $(\tilde{x}(k,0), uv r_{k} + vr_{0})$ and $(\tilde{x}(k+1,0), uv r_{k+1} + vr_{0})$ are connected by a staircase; but in the virtual profile, we connect the two points $(\tilde{f}(k), \tilde{r}_{k})$ and $(\tilde{f}(k+1), \tilde{r}_{k+1})$ by straight line segments. 
\end{itemize}
  We will call the figure so obtained the {\it virtual profile} and use $A^{\vir}$, the area under the virtual profile, to approximate $A$.

\begin{definition}[The virtual profile] \label{virtual profile definition} For each $k = 0, \ldots, \tilde{N}-1$, we define $\tilde{f}(k)$ as follows.  We begin by defining $\tilde{f}_{i}(k)$ for each $i$.  Fix $i$.  Graph the set of points $\{ (\tilde{r}_{k(i,\ell)}, \tilde{c}_{k(i,\ell)}) : \ell = 0, \ldots, \Ki \}$ and connect these by straight line segments.  Then $\tilde{f}_{i}(k)$ is the piecewise linear function whose value at $k$ is the second coordinate of the point on this graph lying over $\tilde{r}_{k}$.  

The picture described above translates into the following rules.  We refer to Cases I-IV as defined in \ref{Cases, s, and t}:  
\begin{enumerate}
\item[0.]  If $\tilde{c}_{k+1,i}=0$, then $\tilde{f}_{i}(k)=0$.  
\item[I.]  In Case I, we have $\tilde{c}_{k+1,i}\neq 0$ and the multiplicity $\tilde{c}_{k,i}$ of $Q_{i}$ jumps at row $k$ (that is, $\tilde{c}_{k,i} < \tilde{c}_{k+1,i}$).  Then $\tilde{f}_{i}(k) = \tilde{c}_{k,i}$.  
\item[II,III,IV.]  Otherwise, let $s(k,i)$ and $t(k,i)$ be as defined in Section  \ref{Cases, s, and t}.  Then
\begin{displaymath} \tilde{f}_{i}(k)   =     \left(  \frac{\tilde{r}_{k}-\tilde{r}_{t(k,i)}}{\tilde{r}_{s(k,i)}-\tilde{r}_{t(k,i)}} \tilde{c}_{s(k,i)} + (1-\frac{\tilde{r}_{k}-\tilde{r}_{t(k,i)}}{\tilde{r}_{s(k,i)}-\tilde{r}_{t(k,i)}}) \tilde{c}_{t(k,i)} \right).
\end{displaymath} 
Note that in Case IV the formula above just gives $\tilde{f}_{i}(k) = \tilde{c}_{k,i}$, since $s(k,i) = k$ in Case IV.  
\end{enumerate}

Finally,
\begin{displaymath}
\tilde{f}(k)  :=  \sum_{i=1}^{q} \tilde{f}_{i}(k).
\end{displaymath}

The virtual profile is the graph of the piecewise linear function connecting the points $\{ (\tilde{f}(k), \tilde{r}_{k}) \}$.
\end{definition}  

Note the switch in the order of the coordinates that takes place: $\tilde{f}_{i}(k)$ is defined by a graph in the $(\mbox{weight} \times \mbox{multiplicity of $Q_{i}$})$-plane, whereas the virtual profile is graphed along with the profile in the $(\mbox{codimension} \times \mbox{weight})$-plane.

The quantity $\tilde{f}(k)$ is an approximate upper bound for the codimension of the $\tilde{r}_{k}$-weight space in $H^{0}(C,\mathcal{O}_{C}(m))$.  We have: 
\begin{eqnarray} A^{\vir}  & = & \sum_{k=0}^{\tilde{N}-1} \frac{1}{2}(\tilde{f}(k+1)-\tilde{f}(k))(\tilde{r}_{k+1}+\tilde{r}_{k}) + (\dim \tilde{V}_{\tilde{N}})vr_{0} \nonumber \\
& = & \sum_{k=0}^{\tilde{N}-1} \frac{1}{2}(\tilde{f}(k+1)-\tilde{f}(k))(\tilde{r}_{k+1}+\tilde{r}_{k}) + (d-d_{\bar{N}})uv + dv -g+1.
\end{eqnarray}

Also, for each $i$ between $1$ and $q$, recall that $r_{j(i,0)}$ is the $r_{j}$ such that $c_{j,i} = 0$ and $c_{j+1,i} > 0$.  Let $T^{\vir} = A^{\vir} + (u+1)^{2} v^{2} \gamma \sum_{i=1}^{q} B_{i} r_{j(i,0)}$ denote the approximation to $T$ obtained by approximating $A$ by $A^{\vir}$.    We have the following upper bound for $T^{\vir}$:
\begin{equation} \label{Tbound} T^{\vir} \leq \sum_{k=0}^{\tilde{N}-1} \frac{1}{2}(\tilde{f}(k+1)-\tilde{f}(k))(\tilde{r}_{k+1}+\tilde{r}_{k}) + ((d-d_{\bar{N}})uv + dv -g+1)vr_{0} + (u+1)^{2} v^{2} \gamma \sum^{n} B_{i} r_{j(i,0)}.
\end{equation}

Before we proceed, I will illustrate the ideas described above by applying them to Example 1.

\subsection{Illustration:  the profile and virtual profile for $\tilde{X}_{\bullet}$ in  Example 1} \label{illustration}
Recall that Example 1 concerns the 1-PS with $q=3$ which induces the following weight filtration:

\begin{displaymath}
\begin{array}{lclc}\mbox{Space} & & & \mbox{Weight} \\ [.0625in]
 V_{0} & = & H^{0}(C,\mathcal{O}(1)) & \hspace{.125in} \frac{1}{2} \\ [.0625in]
V_{1} & = & H^{0}(C,\mathcal{O}(1)(-P_{1})) &  \hspace{.125in} \frac{1}{3} \\ [.0625in]
V_{2} & = & H^{0}(C,\mathcal{O}(1)(-P_{1}-P_{2})) &  \hspace{.125in} \frac{1}{6} \\ [.0625in]
V_{3} & = & H^{0}(C,\mathcal{O}(1)(-P_{1}-P_{2}-P_{3})) &  \hspace{.125in} 0
\end{array}
\end{displaymath}

After passing to $H^{0}(C,\mathcal{O}(m))$ we obtain:

\begin{displaymath}
\begin{array}{lclc}\mbox{Space} & & & \mbox{Weight $\tilde{r}$} \\ [.0625in]
 \tilde{V}_{0} & = & H^{0}(C,\mathcal{O}(m)) & \hspace{.125in} \frac{1}{2}uv + \frac{1}{2}v \\ [.0625in]
\tilde{V}_{1} & = & H^{0}(C,\mathcal{O}(m)(-uvP_{1})) &  \hspace{.125in} \frac{1}{3}uv + \frac{1}{2}v \\ [.0625in]
\tilde{V}_{2} & = & H^{0}(C,\mathcal{O}(m)(-uvP_{1}-uvP_{2})) &  \hspace{.125in} \frac{1}{6}uv + \frac{1}{2}v \\ [.0625in]
\tilde{V}_{3} & = & H^{0}(C,\mathcal{O}(m)(-uvP_{1}-uvP_{2}-uvP_{3})) &  \hspace{.125in}  \frac{1}{2}v
\end{array}
\end{displaymath}  

\subsubsection{The virtual profile for Example 1}
Let us compute the virtual profile first, as this requires fewer calculations than computing $\tilde{X}_{\bullet}$ and the profile.  We can compute the virtual profile for an arbitrary $u,v$:

For $k=0$ there is nothing to compute.  

For $k=1$, the multiplicity of $P_{1}$ does not jump from row 1 to row 2.  We are in Case II.  Looking at where the multiplicity $P_{1}$ jumps, we have $s(1,1) = 0$ and $t(1,1)=3$, and we find that $\tilde{f}_{1}(1) = \frac{1}{3}uv$.   The multiplicity of $P_{2}$ jumps between row 1 and row 2; we are in Case IV, and we have $\tilde{f}_{2}(1) = \tilde{c}_{1,2} = 0$.    Finally, since the multiplicity of $P_{3}$ is zero in both row 1 and row 2, $\tilde{f}_{3}(1) = 0$.  Then $\tilde{f}(1) = \frac{1}{3}uv$.  Also, $\tilde{r}_{1} = \frac{1}{3}uv+ \frac{1}{2}v$.  

For $k=2$, the multiplicity of $P_{1}$ does not jump from row 2 to row 3.  We are in Case II,  $s(2,1) = 0$ and $t(2,1)=3$, and $\tilde{f}_{1}(2) = \frac{2}{3}uv$.   The multiplicity of $P_{2}$ does not jump between row 2 and row 3; we are in Case II, and $s(2,2) = 1$ and $t(2,2) = 3$, giving  $\tilde{f}_{2}(2) = \frac{1}{2}uv$.  Finally, the multiplicity of $P_{3}$ jumps at row 2; we are in Case IV, so $\tilde{f}_{3}(2) = \tilde{c}_{2,3} = 0$.  Then $\tilde{f}(2) = \frac{7}{6}uv$.  Also, $\tilde{r}_{2} = \frac{1}{6}uv+ \frac{1}{2}v$. 

Finally, for $k=\tilde{N}=3$ there is also nothing to compute.

The area of the region under the graph connecting the points $(0uv,\frac{1}{2}uv+ \frac{1}{2}v)$, $(\frac{1}{3}uv,\frac{1}{3}uv+ \frac{1}{2}v)$, $(\frac{7}{6}uv,\frac{1}{6}uv+ \frac{1}{2}v)$ and $(3uv, \frac{1}{2}v)$ is $\frac{1}{2}u^2 v^2 + \frac{3}{2}uv^2$.  To this we add the weight of the $vr_{0}$ region, which is $(\dim \tilde{V}_{\tilde{N}})vr_{0} = ((d-3)uv+dv-g+1)(\frac{1}{2}v)$.  We have:

\[ A^{\vir} = \frac{1}{2}u^2 v^2 + \frac{1}{2}d uv^2 + \frac{1}{2}dv^2 - \frac{1}{2}(g-1)v.
\]   

Using $\gamma B_{i} = \frac{1}{2}$, the contribution from the marked points is $\frac{1}{2}(u^{2}v^2 + 2uv^2 + v^2) $.  We have:

\[ T^{\vir} = 1u^2 v^2 + (\frac{1}{2}d+1) uv^2 + (\frac{1}{2}d+1) v^2 - \frac{1}{2}(g-1)v.
\]

\subsubsection{Interpreting the vertices of the virtual profile}
If we suppose that the integer $uv$ is divisible by 6, we can give a little more meaning to the calculations above.  

For $k=1$ we can begin with the space $\tilde{V}_{1}$, which gives us the point  $(1uv,\frac{1}{3}uv+ \frac{1}{2}v)$.  To this we add the space $V_{0}^{\frac{2}{3}uv}V_{3}^{\frac{1}{3}uv}V_{0}^{v}$ to minimize the multiplicity of $P_{1}$.  Similarly we add $V_{1}^{uv}V_{0}^{v}$ to minimize the multiplicity of $P_{2}$.  The multiplicity of $P_{3}$ is zero in all the spaces of this weight.  The codimension of $V_{1}^{uv}V_{0}^{v}$  is $uv$, and the codimension of $V_{0}^{\frac{2}{3}uv}V_{3}^{\frac{1}{3}uv}V_{0}^{v}$ is also $uv$.  However, using Lemma \ref{trace}, the codimension of their span is $\frac{1}{3}uv$.  In other words, the point $(1uv,\frac{1}{3}uv+ \frac{1}{2}v)$ in the profile of $\tilde{V}_{\bullet}$ slides left to $(\frac{1}{3}uv,\frac{1}{3}uv+ \frac{1}{2}v)$ in the virtual profile for $\tilde{X}_{\bullet}$. 

A similar analysis for $k=2$ yields the list of spaces $V_{0}^{\frac{1}{3}uv}V_{3}^{\frac{2}{3}uv}V_{0}^{v}$, $V_{1}^{\frac{1}{2}uv}V_{3}^{\frac{1}{2}uv}V_{0}^{v}$, and $V_{2}^{uv}V_{0}^{v}$ minimizing the multiplicities of $P_{1}$, $P_{2}$, and $P_{3}$ respectively.  The codimension of their span is $\frac{7}{6}uv$, so the point $(2uv,\frac{1}{6}uv+ \frac{1}{2}v)$ in the profile of $\tilde{V}_{\bullet}$ slides left to $(\frac{7}{6}uv,\frac{1}{6}uv+ \frac{1}{2}v)$ in the virtual profile for $\tilde{X}_{\bullet}$.  

It seems that for any fixed 1-PS $\lambda$ we could choose $uv$ sufficiently divisible to clear any denominators which may arise.  However, we cannot do this across all 1-PS, so we will consider this interpretation of the vertices of the virtual profile as motivational, not part of the rigorous proof.  Also, even when we have such divisibility, so that the virtual profile's vertices have this interpretation, I see no rigorous way to interpret the straight line segments connecting the vertices.  So, it seems best to regard the virtual profile merely as a graph and not an algebro-geometric object of any kind. 

\subsubsection{The filtration $\tilde{X}_{\bullet}$ and its profile for Example 1}
Now we compute the filtration $\tilde{X}_{\bullet}$ and its profile.  For this, we ought to specify $u,v$ first.  We choose $u=3$ and $v=5$.  Of course, this value of $v$ is really too small to use with Lemma \ref{multup lemma}, but let us ignore this in the interest of presenting a reasonably sized example.  Also, in this example, we will always have $\tilde{x}(k,w) = \sum_{i=1}^{q} \tilde{x}(k,w,i)$.  (How do I know this?  See Section \ref{lowerenvelopesarebetter} for a hint.)

The filtration $\tilde{X}_{\bullet}$ has ten stages.  The first and the last are easy to compute---we have $\tilde{X}_{0,0} = H^{0}(C,\mathcal{O}(m))$ and $\tilde{X}_{3} = (V_{3}^{3}V_{0})^{5}$.  Let's compute one of the middle stages, $\tilde{X}_{1,1}$, as an example:  The multiplicity of $P_{1}$ does not increase from row 1 to row 2 to row 3, so we are in Case II, and $s(1,1)=0$ and $t(1,1)=3$.  We find $W = 2$.  (Here $W$ may be computed from its defining properties, or by skipping ahead and using Formula (\ref{Wusefulformula}) derived in Section \ref{virvsactual}.)  Thus the contribution to $\tilde{X}_{1,1}$ from $P_{1}$ is $(V_{0}^{1}V_{3}^{2}V_{0})^{5}$, and $\tilde{x}(1,1,1) = 10$.  The multiplicity of $P_{2}$ increases from row 1 to row 2, but not from row 2 to row 3, so we are in Case IV, and $s(1,2) = 1$ and $t(1,2) = 3$.  Here $W = 1$, and the contribution from $P_{2}$ to $\tilde{X}_{1,1}$ is $(V_{1}^{2}V_{3}^{1}V_{0})^{5}$, and $\tilde{x}(1,1,2) = 5$.  The multiplicity of $P_{3}$ is zero in both row 1 and row 2, so $P_{3}$ does not contribute to $\tilde{X}_{1,1}.$  We have: $\tilde{X}_{1,1} = \Span \{  (V_{1}^{2}V_{2}^{1}V_{0})^{5}, (V_{0}^{1}V_{3}^{2}V_{0})^{5}, (V_{1}^{2}V_{3}^{1}V_{0})^{5}\}$, and $\tilde{x}(1,1) = 15.$

Here is the filtration $\tilde{X}_{\bullet}$.  I have left the spans unsimplified.
\begin{displaymath}
\begin{array}{lclcc} \mbox{Stage} & & \mbox{Space} & \mbox{Codim} & \mbox{Wt} \\
\tilde{X}_{0,0} & = & H^{0}(C,\mathcal{O}(m)) & 0 & 10 \\
\tilde{X}_{0,1} & = & \Span \{ (V_{0}^{2}V_{1}^{1}V_{0})^{5}, (V_{0}^{2}V_{3}^{1}V_{0})^{5} \}  & 5 & 55/6 \\
\tilde{X}_{0,2} & = & \Span \{ (V_{0}^{1}V_{1}^{2}V_{0})^{5}, (V_{0}^{2}V_{3}^{1}V_{0})^{5} \}  & 5 & 50/6 \\
\tilde{X}_{1,0} & = & \Span \{ (V_{1}^{3}V_{2}^{0}V_{0})^{5}, (V_{0}^{2}V_{3}^{1}V_{0})^{5}, (V_{1}^{3}V_{3}^{0}V_{0})^{5}\}  & 5 & 45/6 \\
\tilde{X}_{1,1} & = & \Span \{ (V_{1}^{2}V_{2}^{1}V_{0})^{5}, (V_{0}^{1}V_{3}^{2}V_{0})^{5}, (V_{1}^{2}V_{3}^{1}V_{0})^{5} \}  & 15 & 40/6 \\
\tilde{X}_{1,2} & = & \Span \{ (V_{1}^{1}V_{2}^{2}V_{0})^{5}, (V_{0}^{1}V_{3}^{2}V_{0})^{5}, (V_{1}^{2}V_{3}^{1}V_{0})^{5} \}  & 15 & 35/6 \\
\tilde{X}_{2,0} & = & \Span \{ (V_{2}^{3}V_{3}^{0}V_{0})^{5}, (V_{0}^{1}V_{3}^{2}V_{0})^{5}, (V_{1}^{1}V_{3}^{2}V_{0})^{5}, (V_{2}^{3}V_{3}^{0}V_{0})^{5} \}  & 20 & 5 \\
\tilde{X}_{2,1} & = & \Span \{ (V_{2}^{2}V_{3}^{1}V_{0})^{5}, (V_{0}^{0}V_{3}^{3}V_{0})^{5}, (V_{1}^{1}V_{3}^{2}V_{0})^{5}, (V_{2}^{2}V_{3}^{1}V_{0})^{5}\}  & 30 & 25/6 \\
\tilde{X}_{2,2} & = & \Span \{ (V_{2}^{1}V_{3}^{2}V_{0})^{5}, (V_{0}^{0}V_{3}^{3}V_{0})^{5}, (V_{1}^{0}V_{3}^{3}V_{0})^{5}, (V_{2}^{1}V_{3}^{2}V_{0})^{5} \}  & 40 & 20/6 \\
\tilde{X}_{3} & = & (V_{3}^{3}V_{0})^{5}  & 45 & 15/6 
\end{array}
\end{displaymath}

Notice that $\tilde{X}_{0,1} = \tilde{X}_{0,2} = \tilde{X}_{1,0}$, and $\tilde{X}_{1,1} = \tilde{X}_{1,2}$.  Nothing in our definitions prevents this, and it does not harm us either---all it means is that when we compute the area under the profile between these stages of the filtration, we will obtain a complicated expression for zero.


Here are the profile and virtual profile for $\tilde{X}_{\bullet}$ in Example 1 with $u=3$, $v=5$.  Tick marks on the horizontal axis show units of 5; tick marks on the vertical axis show units of 2.5.  \mbox{} \\

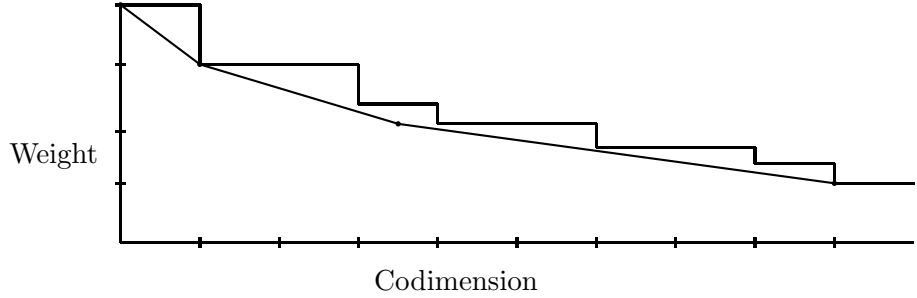
\setlength{\unitlength}{6pt}
\thicklines
\hspace{1in} \begin{picture}(50,15)
\put (0,15){\line(1,0){5}}
\put (5,15){\line(0,-1){3.75}}
\put (5,11.25){\line(1,0){10}}
\put (15,11.25){\line(0,-1){2.5}}
\put (15,8.75){\line(1,0){5}}
\put (20,8.75){\line(0,-1){1.25}}
\put (20,7.5){\line(1,0){10}}
\put (30,7.5){\line(0,-1){1.5}}
\put (30,6){\line(1,0){10}}
\put (40,6){\line(0,-1){1}}
\put (40,5){\line(1,0){5}}
\put (45,5){\line(0,-1){1.25}}
\put (45,3.75){\line(1,0){5}}
\put (0,15){\line(4,-3){5}}
\put (0,15){\circle*{.333}}
\put (5,11.25){\line(10,-3){12.5}}
\put (5,11.25){\circle*{.333}}
\put (17.5,7.5){\line(22,-3){27.5}}
\put (17.5,7.5){\circle*{.333}}
\put (45,3.75){\circle*{.333}}
\put (0,0){\line(1,0){50}}
\put (0,0){\line(0,1){15}}
\put (5,-.333){\line(0,1){.666}}
\put (10,-.333){\line(0,1){.666}}
\put (15,-.333){\line(0,1){.666}}
\put (20,-.333){\line(0,1){.666}}
\put (25,-.333){\line(0,1){.666}}
\put (30,-.333){\line(0,1){.666}}
\put (35,-.333){\line(0,1){.666}}
\put (40,-.333){\line(0,1){.666}}
\put (45,-.333){\line(0,1){.666}}
\put (-.333,3.75){\line(1,0){.666}}
\put (-.333,7){\line(1,0){.666}}
\put (-.333,11.25){\line(1,0){.666}}
\put (-.333,15){\line(1,0){.666}}
\put (16,-3){Codimension}
\put (-7,5){Weight}
\end{picture}
\vspace{0.5in} 

In this picture the area under the profile looks significantly larger than the area under the virtual profile, but for larger values of $u$ these areas become relatively closer.  This is made rigorous in the next section, but as an example, here are the profile and virtual profile for $\tilde{X}_{\bullet}$ in Example 1 with $u=20$, $v=5$.  Tick marks on the horizontal axis show units of 10; tick marks on the vertical axis show units of 5.  \mbox{}  \\ 

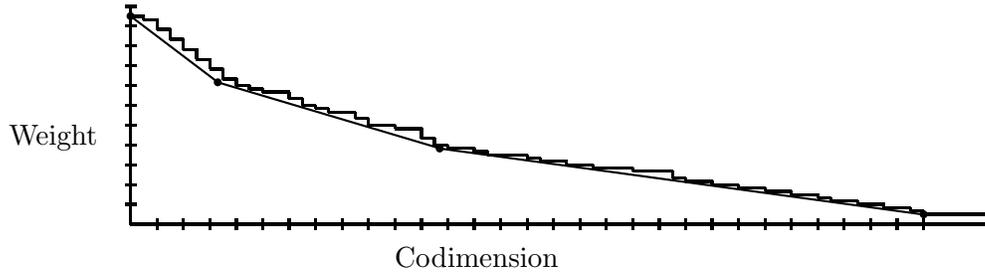
\setlength{\unitlength}{1pt}
\thicklines
\hspace{1in} \begin{picture}(300,90)
\put (0,78.75){\line(1,0){5}}
\put (5,78.75){\line(0,-1){1.25}}
\put (5,77.5){\line(1,0){5}}
\put (10,77.5){\line(0,-1){3.75}}
\put (10,73.75){\line(1,0){5}}
\put (15,73.75){\line(0,-1){3.75}}
\put (15,70){\line(1,0){5}}
\put (20,70){\line(0,-1){3.75}}
\put (20,66.25){\line(1,0){5}}
\put (25,66.25){\line(0,-1){3.75}}
\put (25,62.5){\line(1,0){5}}
\put (30,62.5){\line(0,-1){3.75}}
\put (30,58.75){\line(1,0){5}}
\put (35,58.75){\line(0,-1){3.75}}
\put (35,55){\line(1,0){5}}
\put (40,55){\line(0,-1){2.5}}
\put (40,52.5){\line(1,0){5}}
\put (45,52.5){\line(0,-1){1.25}}
\put (45,51.25){\line(1,0){5}}
\put (50,51.25){\line(0,-1){1.25}}
\put (50,50){\line(1,0){10}}
\put (60,50){\line(0,-1){2.5}}
\put (60,47.5){\line(1,0){5}}
\put (65,47.5){\line(0,-1){2.5}}
\put (65,45){\line(1,0){5}}
\put (70,45){\line(0,-1){1.25}}
\put (70,43.75){\line(1,0){5}}
\put (75,43.75){\line(0,-1){1.25}}
\put (75,42.5){\line(1,0){10}}
\put (85,42.5){\line(0,-1){2.5}}
\put (85,40){\line(1,0){5}}
\put (90,40){\line(0,-1){2.5}}
\put (90,37.5){\line(1,0){10}}
\put (100,37.5){\line(0,-1){1.25}}
\put (100,36.25){\line(1,0){10}}
\put (110,36.25){\line(0,-1){3.75}}
\put (110,32.5){\line(1,0){5}}
\put (115,32.5){\line(0,-1){2.5}}
\put (115,30){\line(1,0){5}}
\put (120,30){\line(0,-1){1.25}}
\put (120,28.75){\line(1,0){10}}
\put (130,28.75){\line(0,-1){1.25}}
\put (130,27.5){\line(1,0){5}}
\put (135,27.5){\line(0,-1){1.25}}
\put (135,26.25){\line(1,0){15}}
\put (150,26.25){\line(0,-1){1.25}}
\put (150,25){\line(1,0){5}}
\put (155,25){\line(0,-1){1.25}}
\put (155,23.75){\line(1,0){10}}
\put (165,23.75){\line(0,-1){1.25}}
\put (165,22.5){\line(1,0){10}}
\put (175,22.5){\line(0,-1){1.25}}
\put (175,21.25){\line(1,0){15}}
\put (190,21.25){\line(0,-1){1.25}}
\put (190,20){\line(1,0){15}}
\put (205,20){\line(0,-1){2.5}}
\put (205,17.5){\line(1,0){5}}
\put (210,17.5){\line(0,-1){1.25}}
\put (210,16.25){\line(1,0){10}}
\put (220,16.25){\line(0,-1){1.25}}
\put (220,15){\line(1,0){10}}
\put (230,15){\line(0,-1){1.25}}
\put (230,13.75){\line(1,0){10}}
\put (240,13.75){\line(0,-1){1.25}}
\put (240,12.5){\line(1,0){10}}
\put (250,12.5){\line(0,-1){1.25}}
\put (250,11.25){\line(1,0){10}}
\put (260,11.25){\line(0,-1){1.25}}
\put (260,10){\line(1,0){5}}
\put (265,10){\line(0,-1){1.25}}
\put (265,8.75){\line(1,0){10}}
\put (275,8.75){\line(0,-1){1.25}}
\put (275,7.5){\line(1,0){10}}
\put (285,7.5){\line(0,-1){1.25}}
\put (285,6.25){\line(1,0){10}}
\put (295,6.25){\line(0,-1){1.25}}
\put (295,5){\line(1,0){5}}
\put (300,5){\line(0,-1){1.25}}
\put (300,3.75){\line(1,0){25}}
\put (0,78.75){\line(33,-25){33}}
\put (0,78.75){\circle*{3}}
\put (33,53.75){\line(84,-25){84}}
\put (33,53.75){\circle*{3}}
\put (117,28.75){\line(183,-25){183}}
\put (117,28.75){\circle*{3}}
\put (300,3.75){\circle*{3}}
\put (0,0){\line(1,0){325}}
\put (0,0){\line(0,1){82.5}}
\put (10,-2){\line(0,1){4}}
\put (20,-2){\line(0,1){4}}
\put (30,-2){\line(0,1){4}}
\put (40,-2){\line(0,1){4}}
\put (50,-2){\line(0,1){4}}
\put (60,-2){\line(0,1){4}}
\put (70,-2){\line(0,1){4}}
\put (80,-2){\line(0,1){4}}
\put (90,-2){\line(0,1){4}}
\put (100,-2){\line(0,1){4}}
\put (110,-2){\line(0,1){4}}
\put (120,-2){\line(0,1){4}}
\put (130,-2){\line(0,1){4}}
\put (140,-2){\line(0,1){4}}
\put (150,-2){\line(0,1){4}}
\put (160,-2){\line(0,1){4}}
\put (170,-2){\line(0,1){4}}
\put (180,-2){\line(0,1){4}}
\put (190,-2){\line(0,1){4}}
\put (200,-2){\line(0,1){4}}
\put (210,-2){\line(0,1){4}}
\put (220,-2){\line(0,1){4}}
\put (230,-2){\line(0,1){4}}
\put (240,-2){\line(0,1){4}}
\put (250,-2){\line(0,1){4}}
\put (260,-2){\line(0,1){4}}
\put (270,-2){\line(0,1){4}}
\put (280,-2){\line(0,1){4}}
\put (290,-2){\line(0,1){4}}
\put (300,-2){\line(0,1){4}}
\put (-2,7.5){\line(1,0){4}}
\put (-2,15){\line(1,0){4}}
\put (-2,22.5){\line(1,0){4}}
\put (-2,30){\line(1,0){4}}
\put (-2,37.5){\line(1,0){4}}
\put (-2,45){\line(1,0){4}}
\put (-2,52.5){\line(1,0){4}}
\put (-2,60){\line(1,0){4}}
\put (-2,67.5){\line(1,0){4}}
\put (-2,75){\line(1,0){4}}
\put (-2,82.5){\line(1,0){4}}
\put (100,-16){Codimension}
\put (-46,30){Weight}
\end{picture}

\subsection*{Progress report}
We have at last defined all the key ingredients mentioned in the introduction: 
\begin{displaymath}
\begin{array}{c}
\mbox{one filtration $V_{\bullet}$ of $H^{0}(C,\mathcal{O}(1))$,} \\
\mbox{two filtrations $\tilde{V}_{\bullet}$ and $\tilde{X}_{\bullet}$ of $H^{0}(C,\mathcal{O}(m))$,} \\
\mbox{ and two graphs associated to $\tilde{X}_{\bullet}$} \\
\end{array}
\end{displaymath}
In Sections \ref{virvsactual}, \ref{boundingTvir}, and \ref{putTogether} it remains to study these filtrations and graphs more closely and show that they have the properties claimed.

\section{The discrepancy between the profile and virtual profile}
\label{virvsactual}
This section is devoted to showing that the areas of the profile and virtual profile are very close when $m$ is large.   That is, we bound the discrepancy  $\Delta := A - A^{\vir}$.  The strategy and methods of this section are extremely straightforward.

We will bound $\Delta$ by computing bounds for several terms which contribute to it.  Roughly speaking, we will compute the discrepancy $\Delta_{k,i}$ for each $k$ and $i$, but it takes a little care to say exactly what we mean by this, as the regions of the graph may be offset a little bit.  For instance, in the picture corresponding to Example 1 with $u=3$, $v=5$, we would partition the virtual profile at codimension 17.5 (a breakpoint of the piecewise linear function) but the corresponding partition for the profile occurs at codimension 20.

For the virtual profile this is straightforward.  The area under the graph of the virtual profile may be divided in an obvious way into $\tilde{N}$ trapezoids and one final rectangle.  Let us focus on the area $A_{k}^{\vir}$ of the $k^{th}$ trapezoid:
\begin{eqnarray*} A_{k}^{\vir} & = & \frac{1}{2}(\tilde{f}(k+1)-\tilde{f}(k))(\tilde{r}_{k+1}+\tilde{r}_{k}) \\
& = & \frac{1}{2}(\sum_{i=1}^{q} \tilde{f}_{i}(k+1)- \sum_{i=1}^{q} \tilde{f}_{i}(k))(\tilde{r}_{k+1}+\tilde{r}_{k}) \\
& = & \sum_{i=1}^{q} \frac{1}{2}( \tilde{f}_{i}(k+1)- \tilde{f}_{i}(k))(\tilde{r}_{k+1}+\tilde{r}_{k}).
\end{eqnarray*}
I will write $A_{k,i}^{\vir}$ for the $i^{th}$ summand:
\[ A_{k,i}^{\vir}  =  \frac{1}{2}(\tilde{f}_{i}(k+1)-\tilde{f}_{i}(k))(\tilde{r}_{k+1}+\tilde{r}_{k}).
\]
We compute $A_{k,i}^{\vir}$ now.    

\subsection{Computing $A_{k,i}^{\vir}$} \label{computeAkivir}
$A_{k,i}^{\vir}$ is the area of the trapezoid whose vertices are $(\tilde{f}_{i}(k),0)$, $(\tilde{f}_{i}(k+1),0)$, $(\tilde{f}_{i}(k+1),\tilde{r}_{k+1})$, and $(\tilde{f}_{i}(k),\tilde{r}_{k})$.  To compute $\tilde{f}_{i}(k+1) - \tilde{f}_{i}(k)$, recall the definition of $\tilde{f}_{i}(k)$ given in Definition \ref{virtual profile definition}.  We use the four cases defined in Section \ref{Cases, s, and t}.    
\begin{enumerate}
\item[I.] The multiplicity $\tilde{c}_{\bullet, i}$ jumps at row $k$ and again at row $k+1$.  Then the spaces contributing to the profile are $V_{k}^{uv}V_{0}^{v}$ at the $k^{th}$ vertex and $V_{k+1}^{uv}V_{0}^{v}$ at the $(k+1)^{th}$ vertex, and in between, spaces of the form $V_{k}^{(u-w)v}V_{k+1}^{wv}V_{0}^{v}$ are used.  Thus in the virtual profile we are calculating as if spaces of the form $V_{k}^{\alpha uv}V_{k+1}^{(1-\alpha)uv}V_{0}^{v}$ were being used between these two vertices with $\alpha$ ranging from 0 to 1.  
  
\item[II.] The multiplicity $\tilde{c}_{\bullet, i}$ does not jump at row $k$ or at row $k+1$.   Recall that we have $(s(k,i),t(k,i)) = (s(k+1,i),t(k+1,i))$.  In the profile, spaces of the form $V_{s(k,i)}^{(u-W)v}V_{t(k,i)}^{W v}V_{0}^{v}$ are being used between these two vertices.  In the virtual profile, we are calculating as if spaces of the form $V_{s(k,i)}^{\alpha uv}V_{t(k,i)}^{(1-\alpha)v} V_{0}^{v}$ were being used between these two vertices (though here the range of $\alpha$ is a subinterval strictly in the interior of $[0,1]$).

\item[III.] The multiplicity $\tilde{c}_{\bullet,i}$ does not jump at row $k$ but jumps at row $k+1$.  Recall that $t(k,i) = k+1$.  Once again, in the profile, spaces of the form $V_{s(k,i)}^{(u-W )v}V_{t(k,i)}^{W v}V_{0}^{v}$ are being used in this region.  For this reason Case III is very similar to Case II.  In the virtual profile, we are calculating as if spaces of the form $V_{s(k,i)}^{\alpha uv}V_{t(k,i)}^{(1-\alpha)v} V_{0}^{v}$ were being used in this region, with $\alpha$ beginning at a value strictly smaller than 1 and decreasing to 0.    
  
\item[IV.] The multiplicity $\tilde{c}_{\bullet,i}$ jumps at row $k$ but not at row $k+1$.  By the definition of $s$ we have $s(k,i) = k$, and in the profile spaces of the form $V_{s(k,i)}^{(u-\wkwi )v}V_{t(k,i)}^{\wkwi v}V_{0}^{v}$ are being used in this region.  In the virtual profile, we are calculating as if spaces of the form $V_{s(k,i)}^{\alpha uv}V_{t(k,i)}^{(1-\alpha)v} V_{0}^{v}$  were being used in this region, with $\alpha$ starting at 1 and ending at a value strictly greater than 0.  
\end{enumerate}

{\it Computing $A_{k,i}^{\vir}$, Case I.}  By Definition \ref{virtual profile definition} we have  $\tilde{f}_{i}(k+1) = \tilde{c}_{k+1,i}$ and $\tilde{f}_{i}(k) = \tilde{c}_{k,i}$.  Thus 
\begin{eqnarray} A_{k,i}^{\vir} & = &  \frac{1}{2}(\tilde{r}_{k+1}+\tilde{r}_{k})(\tilde{f}_{i}(k+1)-\tilde{f}_{i}(k)) \nonumber \\
& = &\frac{1}{2}(uvr_{k+1} + vr_{0} +uvr_{k}+vr_{0})(uvc_{k+1,i} - uv c_{k,i}) \nonumber \\ 
\label{AkivirCaseI} & = & u^2 v^2 (\frac{1}{2}(r_{k+1}+r_{k})(c_{k+1,i} - c_{k,i} + uv^2(r_{0}(c_{k+1,i} - c_{k,i})).
\end{eqnarray}

{\it Cases II, III, and IV.}  In Case II we have
\[ \tilde{f}_{i}(k+1)  = \frac{\tilde{r}_{k+1}-\tilde{r}_{t(k+1,i)}}{\tilde{r}_{s(k+1,i)}-\tilde{r}_{t(k+1,i)}} \tilde{c}_{s(k+1,i)} + (1-\frac{\tilde{r}_{k+1}-\tilde{r}_{t(k+1,i)}}{\tilde{r}_{s(k+1,i)}-\tilde{r}_{t(k+1,i)}}) \tilde{c}_{t(k+1,i)}
\]
and
\[ \tilde{f}_{i}(k)  = \frac{\tilde{r}_{k}-\tilde{r}_{t(k,i)}}{\tilde{r}_{s(k,i)}-\tilde{r}_{t(k,i)}} \tilde{c}_{s(k,i)} + (1-\frac{\tilde{r}_{k}-\tilde{r}_{t(k,i)}}{\tilde{r}_{s(k,i)}-\tilde{r}_{t(k,i)}}) \tilde{c}_{t(k,i)},
\]
and $(s(k,i),t(k,i)) = (s(k+1,i),t(k+1,i))$.  Thus
\begin{eqnarray} A_{k,i}^{\vir} & = &  \frac{1}{2}(\tilde{r}_{k+1}+\tilde{r}_{k})(\tilde{f}_{i}(k+1)-\tilde{f}_{i}(k)) \nonumber \\
& = & \frac{1}{2}(\tilde{r}_{k+1}+\tilde{r}_{k})(\frac{\tilde{r}_{k}-\tilde{r}_{k+1}}{\tilde{r}_{s(k,i)}-\tilde{r}_{t(k,i)}} (\tilde{c}_{t(k,i)} - \tilde{c}_{s(k,i)})) \nonumber \\
& = & \frac{1}{2}(uvr_{k+1}+uvr_{k}+2vr_{0})(uv \frac{r_{k}-r_{k+1}}{r_{s(k,i)}-r_{t(k,i)}}(c_{t(k,i)} - c_{s(k,i)}) ) \nonumber \\
\label{AkivirCaseII} & = & u^2 v^2 \left( \frac{1}{2}(r_{k+1}+r_{k})(c_{t(k,i)} - c_{s(k,i)}) \frac{r_{k}-r_{k+1}}{r_{s(k,i)}-r_{t(k,i)}} \right) \nonumber \\
&&\hspace{.125in}  + uv^2 \left( r_{0} (c_{t(k,i)} - c_{s(k,i)}) \frac{r_{k}-r_{k+1}}{r_{s(k,i)}-r_{t(k,i)}} \right)
\end{eqnarray}
By a similar calculation, and using some of the information presented in paragraphs III and IV above, we derive the same formula in Case III and Case IV.

\subsection{Computing bounds for $A_{k,i}$} \label{Computing bounds for A_k,i}
We have defined $A_{k,i}^{\vir}$ but have not yet defined a corresponding quantity $A_{k,i}$.  We do this now.   Let $A_{k,i}$ denote the following sum:
\begin{equation} \label{Akidef} A_{k,i} := \sum_{w=0}^{u-1} ((u-w)r_{k} + wr_{k+1}+r_{0}) ( \tilde{x}(k,w+1, i) -  \tilde{x}(k,w,i) ).
\end{equation}
In pictures, $\sum_{i=1}^{q}A_{k,i}$ is the area under the profile between $\tilde{x}(k,0)$ and $\tilde{x}(k+1,0)$.


We wish to bound $A_{k,i}$.  We split into Cases I-IV as in Section \ref{computeAkivir}.

{\it Case I.}  Again, using Definition \ref{ProfileDefinition} we have $\tilde{x}(k,w+1,i) = v(u-(w+1)) c_{k,i} + v(w+1)c_{k+1,i}$ and $\tilde{x}(k,w,i) = v(u-w) c_{k,i} + vwc_{k+1,i}$, so $ \tilde{x}(k,w+1,i) -  \tilde{x}(k,w,i) = c_{k+1,i} - c_{k,i}$.  We have:
\begin{eqnarray}A_{k,i} & = &\sum_{w=0}^{u-1} ((u-w)r_{k} + wr_{k+1}+r_{0}) ( \tilde{x}(k,w+1,i) -  \tilde{x}(k,w,i) ) \nonumber \\
& = & \sum_{w=0}^{u-1} ((u-w)r_{k} + wr_{k+1}+r_{0}) (c_{k+1,i} - c_{k,i}) \nonumber \\
&= & \label{AkiCaseI} u^2 v^2 (\frac{1}{2}(r_{k+1}+r_{k})(c_{k+1,i} - c_{k,i})) + uv^2((r_{0}+  \frac{1}{2}(r_{k+1}+r_{k}))(c_{k+1,i} - c_{k,i})).
\end{eqnarray}

{\it Cases II, III, and IV.}  The calculation is long; fortunately, we can treat Cases II, III, and IV. {\it  Also, from here to the end of Section \ref{Computing bounds for A_k,i}, we will suppress the subscripts k,i as much as possible, as they do not change.}  We will reintroduce them at the end of this subsection in line (\ref{AkiCaseII}).

Recall that in Definition \ref{ProfileDefinition}, for each $w$, we defined $W = \wkwi$ to be the smallest integer such that the space $(V_{s}^{u-W} V_{t}^{W} V_{0})^{v}$ has weight less than or equal to $v(u-w)r_{k} + vw r_{k+1} + vr_{0}$.  We use this property to get an expression for $W$ in Case II or Case III:
\[ v(u-W)r_{s} + vW r_{t} + vr_{0} \leq v(u-w)r_{k} + vw r_{k+1} + vr_{0}
\]
\begin{eqnarray}
\Leftrightarrow W & \geq & \frac{u(r_{s} - r_{k}) + w(r_{k} - r_{k+1}) }{r_{s}- r_{t} } \nonumber \\
\label{Wusefulformula} \Rightarrow  W(w)  = \wkwi & = & \left\lceil \frac{u(r_{s} - r_{k}) + w(r_{k} - r_{k+1}) }{r_{s}- r_{t} } \right\rceil
\end{eqnarray}
It is useful to write 
\begin{eqnarray} \label{zetadef} \zeta = \zeta_{k,i} & := & \frac{ r_{k} - r_{k+1} }{ r_{s}- r_{t} } \\
\label{xidef} \xi = \xi_{k,i} & := & \frac{ r_{s} - r_{k} }{ r_{s}- r_{t} }.
\end{eqnarray}
Then 
\begin{equation} \label{omegaexp} W = \lceil u \xi +w \zeta \rceil.
\end{equation}
Also, since $s < k < t$, we have $0 \leq \zeta < 1$ and $0 \leq \xi < 1$.

Proceeding, we have: 
\begin{eqnarray*}
\tilde{x}(k,w,i) & = & v(u-W(w)) c_{s} + v W(w) c_{t} \\
\tilde{x}(k,w+1,i) &  = & v(u-W(w+1)) c_{s} + v W(w+1) c_{t} \\
\Rightarrow \tilde{x}(k,w+1,i) - \tilde{x}(k,w,i)  & = & v (c_{t} -c_{s})(W(w+1)- W(w)). 
\end{eqnarray*}
Putting this into (\ref{Akidef}) we have:
\begin{eqnarray} A_{k,i} & = & \sum_{w=0}^{u-1} v \left( (u-w)r_{k} + w r_{k+1} + r_{0} \right) v \left( (c_{t} -  c_{s} )(W(w+1) - W(w)) \right) \nonumber \\
& = & v^2  (c_{t} -  c_{s} ) \left( \sum_{w=0}^{u-1} (ur_{k} + r_{0} -w(r_{k}-r_{k+1})  )(W(w+1) - W(w) ) \right) \nonumber \\
& = & v^2  (c_{t} -  c_{s} ) \left( (ur_{k} + r_{0}) \sum_{w=0}^{u-1} (W(w+1) - W(w) ) \right. \nonumber \\
&& \label{wkwiline} \hspace{1in} \left. - (r_{k}-r_{k+1}) \sum_{w=0}^{u-1} w(W(w+1) - W(w) ) \right).
\end{eqnarray}

\subsubsection{Calculating pieces of (\ref{wkwiline}) }
Before we continue computing $A_{k,i}$ it is helpful to work out the sums appearing in (\ref{wkwiline}).  We begin with the first sum, $\sum_{w=0}^{u-1} (W(w+1) - W(w) )$.  Let 
\begin{equation}\langle y \rangle := y - \lfloor y \rfloor,
\end{equation}
so that $\langle y \rangle$ denotes the fractional part of $y$.  Then:
\begin{eqnarray}\lefteqn{ \sum_{w=0}^{u-1} (W(w+1) - W(w) ) } \nonumber \\
& = & \sum_{w=0}^{u-1} (\lceil u \xi +w \zeta +\zeta  \rceil - \lceil u \xi +w \zeta  \rceil ) \nonumber \\
\label{anglebracketssum} & = & \sum_{w=0}^{u-1}   (\lceil \langle u \xi\rangle +w \zeta +\zeta  \rceil - \lceil \langle u \xi\rangle +w \zeta  \rceil ).
\end{eqnarray}
Now, imagining the summation as a dynamic process, the sum in line (\ref{anglebracketssum}) increases by one every time the first summand passes an integer and the second summand hasn't caught up yet.  This happens $\lfloor u \zeta + \langle u \xi\rangle \rfloor$ times, so we have 
\begin{equation} \label{firstsum} \sum_{w=0}^{u-1} (W(w+1) - W(w))  = \lfloor u \zeta + \langle u \xi\rangle \rfloor.
\end{equation}
It is helpful to have a nicer expression for $\lfloor u \zeta + \langle u \xi\rangle \rfloor$.  We write 
\begin{displaymath} \lfloor u \zeta + \langle u \xi\rangle \rfloor = u \zeta + \langle u \xi\rangle - \langle u \zeta + \langle u \xi\rangle\rangle
\end{displaymath}
and define
\begin{equation} \label{etadef} \eta := \langle u \xi\rangle - \langle u \zeta + \langle u \xi\rangle\rangle
\end{equation}
so that 
\begin{equation} \label{uzeta+eta} \lfloor u \zeta + \langle u \xi\rangle \rfloor = u \zeta + \eta.
\end{equation}
Note that $-1 < \eta < 1$.

We also compute $\sum_{w=0}^{u-1} w(W(w+1) - W(w) )$.  Simplifying as above, we have:
\begin{eqnarray*}\lefteqn{ \sum_{w=0}^{u-1} w(W(w+1) - W(w) ) } \\
& = & \sum_{w=0}^{u-1}   w (\lceil \langle u \xi\rangle +w \zeta +\zeta  \rceil - \lceil \langle u \xi\rangle +w \zeta  \rceil ).
\end{eqnarray*}
I claim 
\begin{equation} \label{lround} \sum_{w=0}^{u-1}   w (\lceil \langle u \xi\rangle +w \zeta +\zeta  \rceil - \lceil \langle u \xi\rangle +w \zeta  \rceil )  = \left\{ \begin{array}{ll} \displaystyle \sum_{\ell = 1}^{u \zeta + \eta} \left\lceil \frac{\ell - \langle u \xi\rangle}{\zeta} \right\rceil , & \zeta \neq 0 \\ 0, & \zeta=0 \end{array} \right.
\end{equation}
As before, the factor $(\lceil \langle u \xi\rangle +w \zeta +\zeta  \rceil - \lceil \langle u \xi\rangle +w \zeta  \rceil )$ is 0 except when the first summand has just passed an integer and the second summand has not caught up, and then this factor is 1.  We can describe the values of $w$ which are multiplied by nonzero coefficient: for each integer $\ell$ in the appropriate range, we have $w = \left\lceil \frac{\ell - \langle u \xi\rangle}{\zeta} \right\rceil$.

Note $\zeta$ appears in the denominator, and $\zeta$ can take the value 0.  It could be forgetten all too easily that these two things do not happen at the same time, causing concern that this summand (or later quantities) is undefined, so I will write an indicator function $\indz = \indzki$ to remind us that when $\zeta=0$, we add 0.  

\subsubsection{The main calculation resumed}
We now resume the main calculation by reprinting line (\ref{wkwiline}), and then substituting in (\ref{firstsum}), (\ref{uzeta+eta}), and (\ref{lround}):
\begin{eqnarray} A_{k,i} & =  &v^2  (c_{t} -  c_{s} ) \left( (ur_{k} + r_{0}) \sum_{w=0}^{u-1} (W(w+1) - W(w) ) \right. \nonumber \\
&&  \hspace{.5in} \left. - (r_{k}-r_{k+1}) \sum_{w=0}^{u-1} w(W(w+1) - W(w) ) \right) \nonumber \\
& = & v^2  (c_{t} -  c_{s} ) \left( (ur_{k} + r_{0}) (u \zeta + \eta) \right. \nonumber\\ 
&&  \hspace{.5in} \left. - (r_{k}-r_{k+1}) \indz \sum_{\ell = 1}^{u \zeta + \eta} \left\lceil \frac{\ell - \langle u \xi\rangle}{\zeta} \right\rceil \right) \nonumber \\
& = & v^2  (c_{t} -  c_{s} ) \left( (ur_{k} + r_{0}) (u \zeta + \eta) \right. \nonumber \\ 
&&  \hspace{.5in} \left. - (r_{k}-r_{k+1})  \indz \sum_{\ell = 1}^{u \zeta + \eta} \left( \frac{\ell}{\zeta} - \frac{\langle u \xi\rangle}{\zeta} -\left\langle \frac{\ell - \langle u \xi\rangle}{\zeta} \right\rangle +1 \right) \right) \nonumber \\
&= & v^2  (c_{t} -  c_{s} ) \left( (u^2 r_{k}\zeta  + (\eta r_{k} + \zeta r_{0})u + r_{0}\eta \right.  \nonumber \\ 
&&  \hspace{.5in} -\indz \frac{ (r_{k}-r_{k+1})}{\zeta} (\frac{1}{2}(u \zeta+ \eta)(u \zeta+ \eta+1))  \nonumber \\
&& \hspace{.5in}+ \indz \frac{ (r_{k}-r_{k+1})}{\zeta}(\langle u \xi\rangle)(u \zeta+ \eta)  \nonumber \\ 
&&  \hspace{.5in} \left. - \indz \sum_{\ell = 1}^{u \zeta + \eta} \left( 1 -\left\langle \frac{\ell - \langle u \xi\rangle}{\zeta} \right\rangle \right) \right)
\end{eqnarray}
In the last line, we have $0 \leq  1 -\langle\frac{\ell - \langle u \xi\rangle}{\zeta}\rangle $.  Since this quantity is subtracted, we obtain an upper bound for $A_{k,i}$ by replacing this by zero.  We also begin grouping terms by their $u$-degree:
\begin{eqnarray} A_{k,i} & \leq & 
v^2  (c_{t} -  c_{s} ) \left( ( r_{k}\zeta  -\indz \frac{ (r_{k}-r_{k+1})}{\zeta} \frac{1}{2} \zeta^2 ) u^2 \right. \nonumber\\
&& \hspace{.5in} +(\eta r_{k} + \zeta r_{0} + \indz \frac{ (r_{k}-r_{k+1})}{\zeta}(\langle u \xi\rangle -  \frac{1}{2} (2 \eta +1)) \zeta  )u \nonumber \\
&& \hspace{.5in} \left. +(r_{0}\eta + \indz \frac{ (r_{k}-r_{k+1})}{\zeta}(\eta \langle u \xi\rangle - \frac{1}{2}\eta^{2} - \eta )) 1 \right) \nonumber \\
 & = & v^2  (c_{t} -  c_{s} ) \left( (\frac{1}{2}(r_{k}+r_{k+1}) \zeta) u^2 \right. \nonumber \\
&& \hspace{.5in} + (\eta r_{k} + \zeta r_{0} + \indz (r_{k}-r_{k+1})(\langle u \xi\rangle - \eta + \frac{1}{2}) ) u \nonumber \\
&& \hspace{.5in} \left. + (\eta r_{0} + \indz (r_{s}- r_{t})(\eta \langle u \xi\rangle - \frac{1}{2}\eta^{2} - \eta )) 1 \right).  \nonumber
\end{eqnarray}

Finally, we restore the $k,i$ symbols which have been suppressed throughout this subsection, yielding:

\begin{eqnarray}\label{AkiCaseII} A_{k,i} & \leq  & v^2  (c_{t(k,i)} -  c_{s(k,i)} ) \left( (\frac{1}{2}(r_{k}+r_{k+1}) \zeta_{k,i}) u^2 \right. \nonumber \\
&& \hspace{.5in} + (\eta_{k,i} r_{k} + \zeta_{k,i} r_{0} + \indzki (r_{k}-r_{k+1})(\langle u \xi_{k,i} \rangle - \eta_{k,i} + \frac{1}{2}) ) u \nonumber \\
&& \hspace{.5in} \left. + (\eta_{k,i} r_{0} + \indzki (r_{s(k,i)}- r_{t(k,i)})(\eta_{k,i} \langle u \xi_{k,i} \rangle - \frac{1}{2}\eta_{k,i}^{2} - \eta_{k,i} )) 1 \right).
\end{eqnarray}

This completes our calculation of $A_{k,i}$ in Case II, III, or IV.  

\subsection{Bounding the discrepancy}

We now have all the ingredients we need to bound $\Delta$.

\[ \Delta := A - A^{\vir} \leq \sum_{k=0}^{\tilde{N}-1} \sum_{i=1}^{q} ( A_{k,i} - A_{k,i}^{\vir}).
\]

In Case I, by comparing (\ref{AkivirCaseI}) and (\ref{AkiCaseI}) we see that
\begin{eqnarray} \Delta_{k,i} := A_{k,i} - A_{k,i}^{\vir} & = & uv^2 (c_{k+1,i} - c_{k})(\frac{1}{2}(r_{k}-r_{k+1}) \nonumber \\
\label{DeltakiCaseI} & \leq & uv^2 (c_{k+1,i} - c_{k})(\frac{7}{2}) + v^{2} (c_{k+1,i} - c_{k})(3).
\end{eqnarray}
Of course this last estimate is far from sharp, but it is useful to estimate this way to match what appears in Cases II-IV.

In Cases II-IV, by comparing (\ref{AkivirCaseII}) and (\ref{AkiCaseII}) (and using the definition of $\zeta_{k,i}$ at (\ref{zetadef})) we see that 
\begin{eqnarray} \lefteqn{\Delta_{k,i} := A_{k,i} - A_{k,i}^{\vir} } \nonumber \\
 & \leq &  uv^2 ( (c_{t(k,i)} -  c_{s(k,i)} )  (\eta_{k,i} r_{k} + \indzki (r_{k}-r_{k+1})(\langle u \xi_{k,i} \rangle - \eta_{k,i} + \frac{1}{2}) )) \nonumber \\ 
\label{DeltakiCaseIIorIII} && \hspace{.5in} +  v^2 ( (c_{t(k,i)} -  c_{s(k,i)} )  (\eta_{k,i} r_{0} + \indzki (r_{s}- r_{t})(\eta_{k,i} \langle u \xi_{k,i} \rangle - \frac{1}{2}\eta_{k,i}^{2} - \eta_{k,i} )
\end{eqnarray}

Recall that the weights $r_{j}$ and the fractional parts of any quantity must be between 0 and 1, and $-1 < \eta_{k,i} < 1$.  Therefore we may make various coarse estimates:
\begin{eqnarray} \eta_{k,i} r_{k} & < & 1; \nonumber \\
\langle u \xi_{k,i} \rangle - \eta_{k,i} + \frac{1}{2} & < & 1 + 1 + \frac{1}{2} \nonumber \\
\Rightarrow \indzki (r_{k}-r_{k+1})(\langle u \xi_{k,i} \rangle - \eta_{k,i} + \frac{1}{2})  & < & \frac{5}{2}; \nonumber \\
\eta_{k,i} r_{0} & < & 1; \nonumber \\
\eta_{k,i} \langle u \xi_{k,i} \rangle - \frac{1}{2}\eta_{k,i}^{2} - \eta_{k,i} & < & 1 - 0 + 1 = 2 \nonumber \\
\Rightarrow \indzki  (r_{s}- r_{t})(\eta_{k,i} \langle u \xi_{k,i} \rangle - \frac{1}{2}\eta_{k,i}^{2} - \eta_{k,i} ) & < & 3.
\end{eqnarray}
Combining these inequalities with (\ref{DeltakiCaseIIorIII}) we obtain:
\begin{equation}
\label{DeltakiCaseII} \Delta_{k,i}  \leq  uv^2 ( (c_{t(k,i)} -  c_{s(k,i)} ) (\frac{7}{2})  +  v^2 ( (c_{t(k,i)} -  c_{s(k,i)} )(3).
\end{equation}

Next, I claim that the estimates (\ref{DeltakiCaseI}) and (\ref{DeltakiCaseII})
yield
\begin{equation}\label{Deltai} \sum_{k=0}^{\tilde{N}-1} \Delta_{k,i} \leq uv^2 (\frac{7}{2} c_{\bar{N},i}) + v^2 (3 c_{\bar{N},i}).
\end{equation}
Refer back to the definition of $s$ and $t$ in Section \ref{subscriptsI}.  Equation (\ref{Deltai}) follows because the pairs $k,k+1$ from Case I and the pairs $(s,t)$ from Case II, III, and IV  
fit together in such a way that when the estimates  (\ref{DeltakiCaseI}) and (\ref{DeltakiCaseII}) are summed over $k$, the sum telescopes.  

Finally, using the estimates obtained in (\ref{Deltai}), we obtain
\begin{equation} \label{Deltabound} \Delta \leq \sum_{i=1}^{q} \sum_{k=0}^{\tilde{N}-1} \Delta_{k,i} \leq uv^2 (\frac{7}{2}d) + v^2  (3d).
\end{equation}

Observe that $\Delta$ is of order $uv^2$ and not of order $u^2 v^2$.

\section{Bounding $T^{\vir}$}
\label{boundingTvir}

The reader is strongly encouraged to review the subscript notations introduced in Section \ref{subscriptsI}, especially the definitions of $j(i,\ell)$ and $k(i, \ell)$, before proceeding.

\subsection{Setting up a comparison}
  
Recall that in line (\ref{Tbound}) we obtained the following bound on $T^{\vir}$:

\begin{eqnarray} T^{\vir} & \leq &\sum_{k=0}^{\tilde{N}-1} \frac{1}{2}(\tilde{f}(k+1)-\tilde{f}(k))(\tilde{r}_{k+1}+\tilde{r}_{k}) +  ((d-d_{\bar{N}})uv + dv -g+1)vr_{0} \nonumber \\
&& + (u+1)^{2} v^{2} \gamma \sum^{q} B_{i} r_{j(i,0)} \nonumber \\
&\label{nextTbound} = &\sum_{k=0}^{\tilde{N}-1} \frac{1}{2}(\tilde{f}(k+1)-\tilde{f}(k))(\tilde{r}_{k+1}+\tilde{r}_{k}) +  ((d-d_{\bar{N}})uv + dv -g+1)vr_{0} \\
&& + u^{2} v^{2} \gamma \sum^{q} B_{i} r_{j(i,0)} + \sum^{q} \gamma  B_{i} r_{j(i,0)} (2uv^2+v^2) \nonumber
\end{eqnarray}

Everything in this sum is in terms of $k$ (it is, after all, the weight of a basis of $H^{0}(C,\mathcal{O}(m))$).  Almost the only bound available is that the weights sum to 1: $\sum_{j=0}^{\bar{N}} z_{j} r_{j} = 1$.  Our goal in this subsection is to rewrite (\ref{nextTbound}) in a form that makes it easy to compare to $\sum z_{j} r_{j}$. 
 
We focus on the first term of (\ref{nextTbound}): 
\begin{eqnarray} \sum_{k=0}^{\tilde{N}-1} \frac{1}{2}(\tilde{f}(k+1)-\tilde{f}(k))(\tilde{r}_{k+1}+\tilde{r}_{k}) & = & \sum_{k=0}^{\tilde{N}-1} \frac{1}{2}( \sum_{i=1}^{q} \tilde{f}_{i}(k+1)-\sum_{i=1}^{q}\tilde{f}_{i} (k))(\tilde{r}_{k+1}+\tilde{r}_{k}) \nonumber \\ 
\label{Tsubinto} & = & \sum_{k=0}^{\tilde{N}-1} \sum_{i=1}^{q} \frac{1}{2}(  \tilde{f}_{i}(k+1)- \tilde{f}_{i} (k))(\tilde{r}_{k+1}+\tilde{r}_{k})
\end{eqnarray}

Let $A_{k,i}^{\vir}$ denote the area of the region described in Definition \ref{virtual profile definition}.  Then we have:

\begin{eqnarray}  \sum_{k=0}^{\tilde{N}-1} \sum_{i=1}^{q} \frac{1}{2}(  \tilde{f}_{i}(k+1)- \tilde{f}_{i} (k))(\tilde{r}_{k+1}+\tilde{r}_{k}) & = & \sum_{k=0}^{\tilde{N}-1} \sum_{i=1}^{q} A_{k,i}^{\vir} \nonumber \\
& = & \sum_{i=1}^{q} \sum_{k=0}^{\tilde{N}-1} A_{k,i}^{\vir} \label{Aivirdef}
\end{eqnarray}
where in the last line we have changed the order of summation.  Let $A_{i}^{\vir} = \sum_{k=0}^{\tilde{N}-1} A_{k,i}^{\vir}$.    Observe that, for a fixed $i$, it may not be necessary to partition this region into $\tilde{N}$ vertical trapezoids to compute the area $A_{i}^{\vir}$; a partition corresponding to the domains of definition of the piecewise linear function $f_{i}$, which may be coarser than that given by the full set of $k$'s, will do.

Recall that $k(i,\bullet)$ indexes the rows $k$ where the multiplicity $\tilde{c}_{\bullet,i}$ jumps.    Then we may compute:
\begin{eqnarray}  A_{i}^{\vir} & = & \sum_{k=0}^{\tilde{N}-1}  \frac{1}{2}(  \tilde{f}_{i}(k+1)- \tilde{f}_{i} (k))(\tilde{r}_{k+1}+\tilde{r}_{k}) \nonumber \\ 
& = & \sum_{ \ell=0 }^{\Ki-1} \frac{1}{2}(\tilde{c}_{k(i,\ell+1)}-\tilde{c}_{k(i,\ell)})(\tilde{r}_{k(i,\ell+1)}+\tilde{r}_{k(i,\ell)})  \nonumber \\ 
\label{Aivirintermed} & = & u^{2}v^{2} \left(\sum_{\ell=0 }^{\Ki-1} \frac{1}{2}(c_{j(i,\ell+1)}-c_{j(i,\ell)})(r_{j(i,\ell+1)}+r_{j(i,\ell)}) \right) + uv^2(c_{\bar{N},i} r_{0}) 
\end{eqnarray}
We develop the coefficient of the $u^2 v^2$ term of (\ref{Aivirintermed}):
\begin{eqnarray}
 \lefteqn{\left(\sum_{ \ell=0 }^{\Ki-1} \frac{1}{2}(c_{j(i,\ell+1)}-c_{j(i,\ell)})(r_{j(i,\ell+1)}+r_{j(i,\ell)}) \right)} & & \nonumber \\ 
& = & \left(    \sum_{ \ell=1 }^{\Ki}  \frac{1}{2} (c_{j(i,\ell)}-c_{j(i,\ell-1)})r_{j(i,\ell)} +   \sum_{  \ell=0 }^{\Ki-1}  \frac{1}{2}(c_{j(i,\ell+1)}-c_{j(i,\ell)})r_{j(i,\ell)} \right) \nonumber \\
\label{u2v2Aivir} & = &  \left( \sum_{\ell=1 }^{\Ki-1}  \frac{1}{2}(c_{j(i,\ell+1)}-c_{j(i,\ell-1)})r_{j(i,\ell)} +  \frac{1}{2} c_{j(i,1)} r_{j(i,0)}  \right). 
\end{eqnarray}
Once again, $c_{j(i,1)}$ is the first nonzero multiplicity of $Q_{i}$ in a base locus in $V_{\bullet}$, and $r_{j(i,0)}$ is the least weight of a section not vanishing at $Q_{i}$.    Putting (\ref{u2v2Aivir}), (\ref{Aivirintermed}), and (\ref{Aivirdef}) into (\ref{Tsubinto}),  we have:

\begin{eqnarray}  T^{\vir} & \leq & u^2 v^2 \sum_{i=1}^{q}  \left( \sum_{ \ell=1 }^{\Ki-1}  \frac{1}{2} (c_{j(i,\ell+1)}-c_{j(i,\ell-1)})r_{j(i,\ell)}+  (\frac{1}{2}c_{j(i,1)} + \gamma B_{i}) r_{j(i,0)}  \right)  \label{beforeswitching} \\
&& \hspace{.25in}+ \left( \sum_{i=1}^{q} \gamma B_{i} r_{j(i,0)} \right) (2uv^2 + v^2) + ((d-d_{\bar{N}})uv + dv -g+1)vr_{0} + uv^2(\sum_{i=1}^{q} c_{\bar{N},i} r_{0}) \nonumber 
\end{eqnarray}

It is convenient to define $\mathcal{I}_{j}$ to be the set of $i$'s where the multiplicity jumps at row $j$, and not for the first or last time:
\begin{equation} \mathcal{I}_{j} := \{ i \,  | \,  \exists \,  \ell \neq 0, \Ki \,   \mathrm{s.t. } \, j=j(i,\ell) \}. 
\end{equation}

We switch the order of summations in (\ref{beforeswitching}) to obtain:
\begin{eqnarray}
\label{LTbound} T^{\vir} & = & u^2 v^2 \sum_{j=0}^{\bar{N}} \left( \sum_{\mathcal{I}_{j} } \frac{1}{2} (c_{j(i,\ell+1)}-c_{j(i,\ell-1)})  +  \sum_{ i: \,  j=j(i,0)}  (\frac{1}{2}c_{j(i,1)} + \gamma B_{i})  \right) r_{j}  \\
&& \hspace{.25in}+ \left( \sum_{i=1}^{q} \gamma B_{i} r_{j(i,0)} \right) (2uv^2 + v^2) + (duv + dv -g+1)vr_{0} \nonumber 
\end{eqnarray}
which is of the form we desired.

\subsection{Comparing}
The next lemma gives a bound for the coefficient of $u^2v^2$ in (\ref{LTbound}).
\begin{lemma} \label{creep}
\begin{displaymath}
\sum_{j=0}^{\bar{N}}  \left( \sum_{\mathcal{I}_{j} }  \frac{1}{2}(c_{j(i,\ell+1)}-c_{j(i,\ell-1)})  +  \sum_{ i: \,  j=j(i,0)}  (\frac{1}{2}c_{j(i,1)} + \gamma B_{i})   \right) r_{j} \leq  \sum_{j=0}^{\bar{N}} Z_{j}r_{j},
\end{displaymath}
where
\begin{displaymath} Z_{j} := \left\{ \begin{array}{cl} z_{j}, & j < j_{\RR} \\ 
\displaystyle z_{j}+ (\sum_{\tau=0}^{j} z_{\tau} - (\NA - g)), & j = j_{\RR} \\
2 z_{j}, & j \geq j_{\Cliff}
\end{array} \right.
\end{displaymath}
\end{lemma}

{\it Idea of proof (Wall Street version).}  Think of $j$ as being time in days, the $Z_{j}$'s as daily income, and the coefficient of $r_{j}$ on the left hand side as daily losses.  We will show that every time you have a losing day, you have enough in the bank to see you through.

{\it Idea of proof (algebraic geometry version).}  The $Z_{j}$'s defined above bound the change in degree of the base loci from $V_{j}$ to $V_{j+1}$.  The only way there can be a jump larger than this is if $d_{j}$ lags behind the maximum allowable degree for this codimension.  In this case, we are using more small weights and fewer large weights than we conceivably could, so the weight of the resulting basis will not be maximal.

{\it Proof.}  We may rewrite the desired inequality as

\[
\sum_{j=0}^{\bar{N}}  \left( Z_{j} -  \sum_{ i: \,  j=j(i,0)}  (\frac{1}{2}c_{j(i,1)} + \gamma B_{i})  - \sum_{\mathcal{I}_{j} }  \frac{1}{2}(c_{j(i,\ell+1)}-c_{j(i,\ell-1)})  \right) r_{j} \geq 0.
\]

We work successively on each index $j$ where 
\[ Z_{j}  -  \sum_{ i: \,  j=j(i,0) } (\frac{1}{2}c_{j(i,1)} + \gamma B_{i}) -  \sum_{\mathcal{I}_{j} }  \frac{1}{2}(c_{j(i,\ell+1)}-c_{j(i,\ell-1)}) <0.
\]

If there are no such $j$, we are done.  So suppose there is at least one such index, and let the set of these be indexed $j_{e}$ beginning with $e=1$.  By the definition of $j_{1}$ we have 
\[ Z_{j}  -  \sum_{ i: \,  j=j(i,0)}  (\frac{1}{2}c_{j(i,1)} + \gamma B_{i}) -  \sum_{\mathcal{I}_{j} }  \frac{1}{2}(c_{j(i,\ell+1)}-c_{j(i,\ell-1)}) > 0
\] 
for all $j < j_{1}$, so 
\begin{eqnarray*} \lefteqn{\sum_{j=0}^{j_{1}-1} \left( Z_{j} -   \sum_{i: \,  j=j(i,0)}  (\frac{1}{2}c_{j(i,1)} + \gamma B_{i})  - \sum_{\mathcal{I}_{j} }  \frac{1}{2}(c_{j(i,\ell+1)}-c_{j(i,\ell-1)}) \right) r_{j} } && \\ & \geq &  \sum_{j=0}^{j_{1}-1} \left( Z_{j} -     \sum_{i: \,  j=j(i,0)}(\frac{1}{2}c_{j(i,1)} + \gamma B_{i})  - \sum_{\mathcal{I}_{j} }  \frac{1}{2}(c_{j(i,\ell+1)}-c_{j(i,\ell-1)}) \right) r_{j_{1}}
\end{eqnarray*}
and
\[ \sum_{j=0}^{j_{1}-1} \left( Z_{j} -   \sum_{i: \,  j=j(i,0)} (\frac{1}{2}c_{j(i,1)} + \gamma B_{i})  - \sum_{\mathcal{I}_{j} }  \frac{1}{2}(c_{j(i,\ell+1)}-c_{j(i,\ell-1)}) \right) \geq 0. 
\]

We wish to establish that 
\begin{eqnarray*} \lefteqn{ \sum_{j=0}^{j_{1}} \left( Z_{j}  -    \sum_{i: \,  j=j(i,0)} (\frac{1}{2}c_{j(i,1)} + \gamma B_{i})  - \sum_{\mathcal{I}_{j} }  \frac{1}{2}(c_{j(i,\ell+1)}-c_{j(i,\ell-1)}) \right) r_{j}} && \\ & \geq&  \sum_{j=0}^{j_{1}} \left( Z_{j} -    \sum_{i: \,  j=j(i,0)} (\frac{1}{2}c_{j(i,1)} + \gamma B_{i})  - \sum_{\mathcal{I}_{j} }  \frac{1}{2}(c_{j(i,\ell+1)}-c_{j(i,\ell-1)}) \right) r_{j_{1}}
\end{eqnarray*}
(which is easy) and that
\begin{displaymath} \sum_{j=0}^{j_{1}} \left( Z_{j}  -   \sum_{i: \,  j=j(i,0)} (\frac{1}{2}c_{j(i,1)} + \gamma B_{i})  - \sum_{\mathcal{I}_{j}}  \frac{1}{2}(c_{j(i,\ell+1)}-c_{j(i,\ell-1)}) \right)  \geq  0.
\end{displaymath}
We rewrite this last inequality as
\begin{equation}
\label{telescope} \left(\sum_{j=0}^{j_{1}}  Z_{j} \right) - \sum_{j=0}^{j_{1}}  \left(  \sum_{i: \,  j=j(i,0)} (\frac{1}{2}c_{j(i,1)} + \gamma B_{i}) + \sum_{\mathcal{I}_{j}}  \frac{1}{2}(c_{j(i,\ell+1)}-c_{j(i,\ell-1)})\right) \geq 0.
\end{equation}

We study the second sum in (\ref{telescope}) above.  Each $i$ falls into exactly one of the following cases:

{\it Case 0.} If $c_{\bullet,i}$ does not jump before or at $j_{1}$---that is, $j(i,0)> j_{1}$---then this $i$ does not contribute.  

{\it Case 1.} If $c_{\bullet,i}$ jumps exactly once before or at $j_{1}$---that is, $j(i,0) \leq j_{1} < j(i,1)$---then this $i$ contributes  \[ \frac{1}{2}c_{j(i,1)} + \gamma B_{i} \leq \frac{1}{2}c_{j_{1}+1,i} + \frac{1}{2} \leq c_{j_{1}+1,i},
\]
since $c_{j(i,1)} = c_{j_{1}+1,i}$ and $\gamma B_{i} \leq \frac{1}{2}$ and $c_{j_{1}+1,i} \geq 1$.  

{\it Case 2.} If $c_{\bullet,i}$ jumps exactly twice before or at $j_{1}$---that is, $j(i,1) \leq j_{1} < j(i,2)$---then the contribution to the second term is 
\[ \frac{1}{2}c_{j(i,1)} + \gamma B_{i} + \frac{1}{2}c_{j(i,2)} \leq c_{j_{1}+1,i}.
\]
This follows because $c_{j(i,2)}= c_{j_{1}+1,i}$ and $c_{j(i,2)} \geq c_{j(i,1)} + 1$.  
  
{\it Case 3.} If $c_{\bullet,i}$ jumps three or more times before or at $j_{1}$, then some telescoping occurs, and the contribution is 
\[\frac{1}{2}c_{j(i,1)} + \gamma B_{i} +\frac{1}{2}c_{t(j_{1},i)}+ \frac{1}{2}c_{s(j_{1},i)} - \frac{1}{2}c_{j(i,1)} \leq c_{j_{1}+1,i}.  
\]
Here I am abusing notation a little (according to Section \ref{subscriptsI} the first argument of $s(\bullet,i)$ or $t(\bullet,i)$ is supposed to be a $k$, not a $j$).  Here $s(j_{1},i)$ denotes the largest index less than or equal to $j_{1}$ where $c_{\bullet,i}$ jumps, and $t(j_{1},i)$ denotes the smallest index strictly greater than $j_{1}$ index where  $c_{\bullet,i}$ jumps.  Thus, $c_{t(j_{1},i)} = c_{j_{1}+1,i}$ and $c_{s(j_{1},i)} \leq c_{j_{1},i}$.

To summarize, in each case, we see that the contribution is no more than $c_{j_{1}+1,i}$.  

If $j_{1} < j_{\RR}$, so that $j_{1}+1$ is in the Riemann-Roch region, then by (\ref{RR region}) we have

\begin{displaymath} 
\sum_{i=1}^{q} c_{j_{1}+1,i} \leq \sum_{j=0}^{j_{1}} z_{j}, 
\end{displaymath}
so the left hand side of (\ref{telescope}) is indeed nonnegative:
\begin{eqnarray*}  \left(\sum_{j=0}^{j_{1}}  Z_{j} \right) - \sum_{j=0}^{j_{1}}  \left(  \sum_{i: \,  j=j(i,0)}(\frac{1}{2}c_{j(i,1)} + \gamma B_{i}) + \sum_{ \mathcal{I}_{j} }  \frac{1}{2}(c_{j(i,\ell+1)}-c_{j(i,\ell-1)}) \right) &&   \\
 & & \hspace{-3in}  \geq \left(\sum_{j=0}^{j_{1}}  z_{j} \right) - \left( \sum_{j=0}^{j_{1}} z_{j} \right)  = 0.
\end{eqnarray*}

We have thus dealt with the first index, if it falls inside the Riemann-Roch region.  We may repeat the argument at each $j_{e}$ in the Riemann-Roch successively, stopping when either the $j_{e}$'s are exhausted or we reach the Clifford region.  At each step we need to show two things in order to proceed to the next step: first, 

\begin{eqnarray*} \lefteqn{ \sum_{j=0}^{j_{e}} \left( Z_{j}  -    \sum_{i: \,  j=j(i,0)} (\frac{1}{2}c_{j(i,1)} + \gamma B_{i})  - \sum_{\mathcal{I}_{j} }  \frac{1}{2}(c_{j(i,\ell+1)}-c_{j(i,\ell-1)}) \right) r_{j}} && \\ & \geq&  \sum_{j=0}^{j_{e}} \left( Z_{j} -    \sum_{i: \,  j=j(i,0)} (\frac{1}{2}c_{j(i,1)} + \gamma B_{i})  - \sum_{\mathcal{I}_{j} }  \frac{1}{2}(c_{j(i,\ell+1)}-c_{j(i,\ell-1)}) \right) r_{j_{1}}
\end{eqnarray*}
(which is always easy to check), and second,
\begin{displaymath} \sum_{j=0}^{j_{e}} \left( Z_{j}  -   \sum_{i: \,  j=j(i,0)} (\frac{1}{2}c_{j(i,1)} + \gamma B_{i})  - \sum_{\mathcal{I}_{j} }  \frac{1}{2}(c_{j(i,\ell+1)}-c_{j(i,\ell-1)}) \right)  \geq  0.
\end{displaymath}

Next suppose that $j_{e} = j_{\RR}$, so $j_{e}+1 = j_{\Cliff}$.   Then by (\ref{Clifford region}) we have
\begin{displaymath} 
\sum_{i=1}^{q} c_{j_{e}+1,i} \leq \sum_{j=0}^{j_{e}} z_{j} + \sum_{j=0}^{j_{e}} z_{j} - (\NA-g)
\end{displaymath}
\begin{eqnarray*}  \left(\sum_{j=0}^{j_{e}}  Z_{j} \right) - \sum_{j=0}^{j_{e}}  \left(  \sum_{i: \,  j=j(i,0)} (\frac{1}{2}c_{j(i,1)} + \gamma B_{i}) + \sum_{\mathcal{I}_{j} }  \frac{1}{2}(c_{j(i,\ell+1)}-c_{j(i,\ell-1)})\right)&&   \\
 & & \hspace{-5in}  \geq \left(\sum_{j=0}^{j_{e}}  z_{j} + \sum_{j=0}^{j_{e}} z_{j} - (\NA-g)  \right) - \left( \sum_{j=0}^{j_{e}} z_{j} + \sum_{j=0}^{j_{e}} z_{j} - (\NA-g)  \right)  = 0.
\end{eqnarray*}

Finally suppose that some $j_{e}+1$ falls within the Clifford region.  Then by (\ref{Clifford region}) we have 
\begin{displaymath} 
\sum_{i=1}^{q} c_{j_{e}+1,i} \leq \sum_{j=0}^{j_{e}} z_{j} + \sum_{j=0}^{j_{e}} z_{j} - (\NA-g).
\end{displaymath}

Using the definitions given in the statement of the lemma, we compute
\begin{eqnarray*} \sum_{j=0}^{j_{e}} Z_{j} & = & \sum_{j=0}^{j_{\RR}-1} Z_{j} + Z_{j_{\RR}}  + \sum_{j = j_{\Cliff}}^{j_{e}} Z_{j} \\
& = & \sum_{j=0}^{j_{\RR}-1} z_{j} + \left(z_{j_{\RR}} + \sum_{j=0}^{j_{\RR}} z_{j} + (\NA-g) \right) + 2z_{j_{\Cliff}} + \cdots + 2z_{j_{e}} \\
& = & 2 \sum_{j=0}^{j_{e}} z_{j} - (\NA -g)
\end{eqnarray*}

and once again the left hand side of (\ref{telescope}) is nonnegative:
\begin{eqnarray*}   \left(\sum_{j=0}^{j_{e}}  Z_{j} \right) - \sum_{j=0}^{j_{e}}  \left(  \sum_{i: \,  j=j(i,0)} (\frac{1}{2}c_{j(i,1)} + \gamma B_{i}) + \sum_{\mathcal{I}_{j} }  \frac{1}{2}(c_{j(i,\ell+1)}-c_{j(i,\ell-1)})\right) &&   \\
 & & \hspace{-5in}  \geq \left(2 \sum_{j=0}^{j_{e}}  z_{j} -(\NA-g) \right) - \left( \sum_{j=0}^{j_{e}} z_{j}  + \sum_{j=0}^{j_{e}} z_{j} - (\NA-g) \right)  = 0.
\end{eqnarray*}

Again, proceed to the next $j_{e}$ until the set of these has been exhausted. $\hfill \Box$

\vspace{.25in}

Ideally, we would now show that the bound obtained in Lemma \ref{creep} is smaller than what is required in the numerical criterion.  Unfortunately, this is not always true.  Lemma \ref{creep} is sufficient for most, but not all, sets of linearizing weights $\mathcal{B}$.   Below I have listed five cases which exhaust all possibilities.  This partitioning may look strange, but it is in order of difficulty of proof.  In Cases A-C, I can prove asymptotic stability of smooth curves.   In Cases D and E, I cannot prove stability, so I will ultimately impose hypotheses to ensure that these cannot occur.

  Choose any sufficiently small value $\epsilon > 0$.  (The size of $\epsilon$ allowed will become clear in Cases B and C below, and the role of $\epsilon$ will become clear in the proof of Theorem \ref{stabilitytheorem}.)  Then we consider the following five cases:  
\begin{equation} \label{CaseICaseII}
\begin{array}{lccccccl}
\mbox{{\it Case A.}} & n & \geq & 1 & \mbox{and}  & \gamma b & \geq & \frac{g-1}{N} + \epsilon (N+1). \\
\mbox{{\it Case B. }}& n & \geq  & 1 & \mbox{and} & \gamma b & < & \frac{g-1}{N} + \epsilon (N+1)  <  \frac{1}{2}\\
\mbox{{\it Case C. }}& n & = & 0 & \mbox{and} & N & \geq & 2g-2  \\
\mbox{{\it Case D. }}& n & = & 0 & \mbox{and} & N & < & 2g-2  \\
\mbox{{\it Case E. }}& n & \geq & 1 & \mbox{and} & \gamma b & < & \frac{g-1}{N} + \epsilon (N+1)  \geq  \frac{1}{2} \\
\end{array}
\end{equation}

Let us proceed first with Case A:  To apply Lemma \ref{creep} to our problem, we need to bound $\sum Z_{j}r_{j} $.  Let $r_{\NA-g+1}, \ldots, r_{\NA-1}, r_{\NA}=0$ be the last $g$ weights (that is, ignore the index $j$ and list the smallest weights as many times as indicated by their multiplicities). 
   Then we have 
\begin{eqnarray*}\sum Z_{j}r_{j}  & \leq & \sum z_{j}r_{j} +  r_{N-g+1} + \cdots + r_{N} \\
&  \leq & 1 + r_{N-g+1} + \cdots + r_{N}  
\end{eqnarray*}
Now we bound $r_{N-g+1} + \cdots + r_{N}$:

\begin{lemma} \label{Clifftailbound} $r_{N-g+1} + \cdots + r_{N} \leq \frac{g-1}{N}$.
\end{lemma}
{\it Proof.}  Recall that $r_{N}= 0$, so we may omit it from all the following sums.  We argue similarly to \cite{Morr} Theorem 4.1.  We wish to maximize $r_{N-g+1} + \cdots + r_{N-1}$, which is linear in the $r$'s, subject to the constraints $\sum_{j=0}^{\bar{N}-1} z_{j}r_{j} =1$ and that the $r$'s are decreasing.  In the affine hyperplane in $(N-1)$-dimensional $r$-space determined by the equation $\sum_{j=0}^{\bar{N}-1} z_{j}r_{j} =1$, the condition that the $r$'s are decreasing defines an $(N-1)$-simplex.  The vertices of this simplex correspond to sequences of the following form:
\[ r_{0} = \cdots = r_{h} > r_{h+1} = \cdots = r_{N-1} = 0.
\]
The function must take its maximum at (at least) one of these vertices, and it is easy to check that the maximum occurs when 
\[ r_{0} = \cdots = r_{N-1} > 0,
\]
or $r_{j} = \frac{1}{N}$ for all $j$, yielding a maximum value of $\frac{g-1}{N}$.  $\hfill \Box$

\vspace{.5in}

Also, the defining hypothesis of Case A at line (\ref{CaseICaseII}) may be written as follows. 
\begin{eqnarray}
\gamma b & \geq & \frac{g-1}{N} + \epsilon (N+1) \nonumber \\
\Leftrightarrow \frac{g-1}{N} &  \leq & \frac{g-1+ \gamma b}{N+1} - \epsilon \nonumber
\end{eqnarray}
Therefore, as a trivial extension of Lemma \ref{Clifftailbound}, we have:
\begin{equation} \label{desiredrbound} r_{N-g+1} + \cdots + r_{N} \leq \frac{g-1+ \gamma b}{N+1} - \epsilon 
\end{equation}
We combine (\ref{desiredrbound}) with the bound found in (\ref{LTbound}) to obtain:
\begin{equation} \label{Tvirbound}  T^{\vir}  \leq   \left( 1 + \frac{g-1+ \gamma b}{N+1} - \epsilon  \right) u^2 v^2 +  \left( \sum_{i=1}^{q} \gamma B_{i} r_{j(i,0)} \right) (2uv^2 + v^2) + (duv + dv -g+1)vr_{0}.
\end{equation}
Note that the leading coefficient $ 1 + \frac{g-1+ \gamma b}{N+1} - \epsilon $ is less than the leading coefficient $1+ \frac{g-1+ \gamma b}{N+1}$ of the numerical criterion (\ref{numcrit}) by $\epsilon$.  This completes our discussion of Case A.

Next we turn to Cases B and C, defined in line (\ref{CaseICaseII}).  In these cases, the bound given in Lemma \ref{Clifftailbound} is too large to use with the numerical criterion.  Fortunately, if we examine the proof of Lemma \ref{creep} closely, we can improve the bound there a little bit.

\begin{lemma} \label{betterthancreep} \begin{enumerate}
\item Suppose a sufficiently small $\epsilon >0$ has been chosen and $n \geq 1$ and  $\gamma b <  \frac{g-1}{N} + \epsilon (N+1) < \frac{1}{2}$, so that we are in Case B.  Then 
\[
\sum_{j=0}^{\bar{N}}  \left( \sum_{\mathcal{I}_{j} }  \frac{1}{2}(c_{j(i,\ell+1)}-c_{j(i,\ell-1)})  +  \sum_{ i: \,  j=j(i,0)}  (\frac{1}{2}c_{j(i,1)} + \gamma B_{i})   \right) r_{j} \leq  \sum_{j=0}^{\bar{N}} Z_{j}r_{j} - \left( \frac{1}{2} - \gamma b \right) r_{N-1},
\]
where the $Z_{j}$ are as in Lemma \ref{creep}, and 
\[
r_{N-1} = \left\{ \begin{array}{ll} 0, & z_{\bar{N}} >1
\\ r_{\bar{N}-1}, & z_{\bar{N}} = 1. 
\end{array} \right.
\]  
\item Suppose $n=0$.  Then 
\[
\sum_{j=0}^{\bar{N}}  \left( \sum_{\mathcal{I}_{j} }  \frac{1}{2}(c_{j(i,\ell+1)}-c_{j(i,\ell-1)})  +  \sum_{ i: \,  j=j(i,0)}  (\frac{1}{2}c_{j(i,1)} + \gamma B_{i})   \right) r_{j} \leq  \sum_{j=0}^{\bar{N}} Z_{j}r_{j} -  \frac{1}{2}r_{N-1},
\]
where the $Z_{j}$ are as in Lemma \ref{creep}, and 
\[
r_{N-1} = \left\{ \begin{array}{ll} 0, & z_{\bar{N}} >1
\\ r_{\bar{N}-1}, & z_{\bar{N}} = 1. 
\end{array} \right.
\]  
\end{enumerate}
\end{lemma}
{\it Proof.}  Note this is a trivial extension of Lemma \ref{creep} if $z_{\bar{N}} > 1$, as then $r_{N-1} = 0$.  So suppose $z_{\bar{N}} = 1$; then $\sum_{j=0}^{\bar{N} - 1} z_{j} = N-1$.  By the proof of Lemma \ref{creep} we know that
\begin{eqnarray*} \lefteqn{ \sum_{j=0}^{\bar{N}-1} \left( Z_{j}  -    \sum_{i: \,  j=j(i,0)} (\frac{1}{2}c_{j(i,1)} + \gamma B_{i})  - \sum_{\mathcal{I}_{j}}  \frac{1}{2}(c_{j(i,\ell+1)}-c_{j(i,\ell-1)}) \right) r_{j}} && \\ & \geq&  \sum_{j=0}^{\bar{N}-1} \left( Z_{j} -    \sum_{i: \,  j=j(i,0)} (\frac{1}{2}c_{j(i,1)} + \gamma B_{i})  - \sum_{\mathcal{I}_{j} }  \frac{1}{2}(c_{j(i,\ell+1)}-c_{j(i,\ell-1)}) \right) r_{\bar{N}-1}
\end{eqnarray*}
and 
\begin{displaymath} \sum_{j=0}^{\bar{N}-1} \left( Z_{j}  -   \sum_{i: \,  j=j(i,0)} (\frac{1}{2}c_{j(i,1)} + \gamma B_{i})  - \sum_{\mathcal{I}_{j}}  \frac{1}{2}(c_{j(i,\ell+1)}-c_{j(i,\ell-1)}) \right)  \geq  0.
\end{displaymath}
So if 
\begin{equation} \label{Zminusc} \sum_{j=0}^{\bar{N}-1} Z_{j} - \sum_{i=1}^{q} c_{\bar{N},i} \geq \frac{1}{2} - \gamma b,
\end{equation}
then we are done.  Note the left hand side of (\ref{Zminusc}) is a nonnegative integer.  So suppose the left hand side of (\ref{Zminusc}) is zero; we will explain how to improve the estimates used in the proof of Lemma \ref{creep} by at least $\frac{1}{2} - \gamma b$.  

First, if $n=0$, there are no marked points, and $B_{i} = 0$ for all $i$.  Since we estimated $\gamma B_{i} \leq \frac{1}{2}$, we have the improvement we need.

So suppose $n \geq 1$.  If there is at least one point $Q_{i}$ appearing in a base locus in $V_{\bullet}$ which is not one of the marked points $P_{i}$, then similarly since $B_{i} = 0$ and we always estimated $\gamma B_{i} \leq \frac{1}{2}$, we have the improvement we need.  So we may suppose that every $Q_{i}$ is a $P_{j}$ (hence $q<n$).

If there are no points $Q_{i}$---that is, the base locus of $V_{\bar{N}}$ is empty---then the weight $vr_{0}$ space has codimension 0 in $H^{0}(C,\mathcal{O}(m))$, and we can easily show $T^{\vir}$ is smaller than what is required by the numerical criterion.

So suppose there is at least one point $Q_{1}$ in the base locus of $V_{\bar{N}}$.  But now, on the one hand we have by hypothesis that $\gamma B_{i} \leq \gamma b <  \frac{g-1}{N+1} + \epsilon (N+1) \leq \frac{1}{2}$; but in the proof of Lemma \ref{creep} we only estimated $\gamma B_{i} \leq \frac{1}{2}$; so we see that we may improve our estimate by at least the desired amount.  
$\hfill \Box$

We proceed with Case B.  We may argue just as we did in Lemma \ref{Clifftailbound} to get
\begin{equation} \label{ClifftailboundB} r_{N-g+1} + \cdots + r_{N} - (\frac{1}{2} - \gamma b)r_{N-1} \leq  \left( g-1 - (\frac{1}{2}  - \gamma b) \right)\frac{1}{N}
\end{equation}
Combining (\ref{ClifftailboundB}) with (\ref{LTbound}), we obtain:
\begin{equation} \label{TvirboundII}  T^{\vir}  \leq   \left( 1 + \frac{g-\frac{3}{2} + \gamma b}{N}  \right) u^2 v^2 +  \left( \sum_{i=1}^{q} \gamma B_{i} r_{j(i,0)} \right) (2uv^2 + v^2) + (duv + dv -g+1)vr_{0}.
\end{equation}
We desire that the leading coefficient should be smaller than what is required by the numerical criterion by $\epsilon$.  That is, we want:
\begin{eqnarray} \frac{g-\frac{3}{2} + \gamma b}{N} & \leq & \frac{g-1+ \gamma b}{N+1} - \epsilon \nonumber \\
\label{epsilonlineB} \Leftrightarrow \epsilon & \leq & \frac{1}{2N(N+1)} (N - 2g + 3 - 2\gamma b).
\end{eqnarray}
The right hand side of (\ref{epsilonlineB}) is positive because the hypotheses of Case B imply that $N \geq 2g-1$, and we also have  $\gamma b < \frac{1}{2}$.  Thus, when $\epsilon$ is sufficiently small (depending on $N$, $\nu$, and $\mathcal{B}$) then (\ref{epsilonlineB}) is satisfied.

Next we consider Case C.  Lemma \ref{betterthancreep}.2 covers this situation, and we may argue just as we did in Lemma \ref{Clifftailbound} to get
\begin{equation} r_{N-g+1} + \cdots + r_{N} - \frac{1}{2}r_{N-1} \leq  \left( g - \frac{3}{2}  \right)\frac{1}{N}
\end{equation}
Then, we want to arrange that 
\begin{eqnarray} \frac{g-\frac{3}{2} }{N} & \leq & \frac{g-1}{N+1} - \epsilon \nonumber \\
\label{epsilonlineC} \Leftrightarrow  \epsilon & \leq & \frac{1}{2N(N+1)} (N-2g+3).
\end{eqnarray}
Since $N \geq 2g-2$ in Case C, (\ref{epsilonlineC}) is satisfied for all $\epsilon$ sufficiently small.  

This completes our discussion of Cases B and C.

Unfortunately, in Cases D and E, I know of no way to improve the bound of Lemma \ref{creep} in order to get the leading coefficient of $T^{\vir}$ small enough to use with the numerical criterion in this case!  Therefore, at present I am forced to make the following hypotheses to ensure that Cases D and E do not occur:
\begin{enumerate}
\item If $n=0$, then $N \geq 2g-2$.
\item If $n \geq 1$ and $g \geq 2$ then either $\gamma b  \geq  \frac{g-1}{N} + \epsilon (N+1)$ or else $\gamma b  <  \frac{g-1}{N} + \epsilon (N+1)  <  \frac{1}{2}$.
\end{enumerate}
Note that for $n \geq 1$ and $g=0$ or $g=1$ and $b >0$, we always have  $\gamma b  \geq  \frac{g-1}{N} + \epsilon (N+1)$, so this hypothesis does not impose any restriction on $d$ or $N$ in these cases; we only need the linear system embedding the curve to be complete.

\section{GIT stability of smooth pointed curves}
\label{putTogether}

\subsection{The stability theorem}
We are ready to prove the main result:

\begin{theorem} \label{stabilitytheorem} Let $\gamma = \nu/(2 \nu-1)$.   Choose any $\epsilon >0$ which is sufficiently small depending on $d$, $g$, and $n$.  If $n=0$ assume $N \geq 2g-2$.   If $n \geq 1$ and $g \geq 2$ then suppose $\gamma b  \geq  \frac{g-1}{N} + \epsilon (N+1)$ or else $\gamma b  <  \frac{g-1}{N} + \epsilon (N+1)  <  \frac{1}{2}$.    Consider a point in the incidence locus $I$ parametrizing a smooth pointed curve $(C, \{ P_{i} \})$ embedded in $\Pro^{\NA}$ by any (i.e. not necessarily pluricanonical) complete linear system of degree $\dA$.  Assume also that the points $P_{i}$ are distinct. 
 
If $n \geq 1$, suppose each $b_{i} \in \mathcal{B}$ satisfies $\gamma b_{i} < \frac{1}{2}$ (this may not be covered by the previous assumptions).  Let $m = (u+1)v$.   Then for certain large values of $m$, the point of $I$ parametrizing $(C, \{ P_{i} \}, C \subset \Pro^{N}) $  is GIT stable for the $SL(N+1)$-action with the linearization specified by $m_{i}' = \gamma b_{i} m^{2}$ for each $i$.  More precisely, there exist:
\begin{enumerate}
\item a positive integer $u_{0}$  depending on $d$, $g$, $n$, and $\mathcal{B}$, but not on the curve $C$, the points $P_{i}$, or the embedding $C \subset \Pro^{N}$
\item a function $v_{0}(u)$ whose domain is all integers greater than $u_{0}$, and which depends on $u$, $d$, $g$, $\mathcal{B}$ and $\epsilon$  but not on the curve $C$, the points $P_{i}$, or the embedding $C \subset \Pro^{N}$
\end{enumerate}
such that for any integers $u \geq u_{0}$ and $v \geq v_{0}(u)$, the point of $I$ parametrizing $(C, \{ P_{i} \}, C \subset \Pro^{N}) $  is GIT stable for the $SL(N+1)$-action with the linearization specified by $m_{i}' = \gamma b_{i} m^{2}$ for each $i$.
\end{theorem}
{\it Proof.}  By (\ref{Tvirbound}) and (\ref{Deltabound}) we have
\begin{eqnarray} T & = & T^{\vir} + \Delta  \nonumber \\
& \leq &  \left( 1 + \frac{g-1+ \gamma b }{N+1} - \epsilon \right) u^2 v^2 +  \left( \sum_{i=1}^{q} \gamma b_{i} r_{j(i,0)} \right) (2uv^2 + v^2) \nonumber \\
&& \hspace{.5in}+ (duv + dv -g+1)vr_{0}   + \frac{7}{2}d uv^2  + 3d v^2 \nonumber \\
& = & \left( 1 + \frac{g-1+ \gamma b}{N+1} - \epsilon \right) u^2 v^2 + \left(2 \sum_{i=1}^{n} \gamma b_{i} r_{j(i,0)} + \frac{7}{2}d \right) uv^2 + \left(2 \sum_{i=1}^{n} \gamma b_{i} r_{j(i,0)} + 3d \right) v^2 \nonumber \\
& \leq & \left( 1 + \frac{g-1 + \gamma b}{N+1} - \epsilon \right) u^2 v^2 + \left(2 \gamma b + \frac{7}{2}d \right) uv^2 + \left(2 \gamma b + 3d \right) v^2 \nonumber \\
\label{lastTbound} & \leq & \left( 1 + \frac{g-1 + \gamma b}{N+1} - \epsilon \right) u^2 v^2 + \left(n + \frac{7}{2}d \right) uv^2 + \left(n + 3d \right) v^2
\end{eqnarray}

Note that this bound depends on $d$, $g$, and $n$.  Therefore, in the important special case when $d = \nu(2g-2+a)$, it also depends on $\nu$ and $a$.   But we emphasize that in every case, this bound does \emph{not} depend on the particular curve $C$, the points $P_{i}$, the embedding $C \subset \Pro^{N}$, or the 1-PS $\lambda$.

Recall the bound required in the numerical criterion: 
\begin{equation} \label{finalNumCrit} \left( 1 + \frac{g-1+ \gamma b}{\NA+1} \right) m^2  - \frac{g-1}{\NA+1} m = \left( 1 + \frac{g-1+ \gamma b}{\NA+1} \right) (u^2v^2 + 2uv^2 +v^2) - \frac{g-1}{\NA+1} (uv+v).
\end{equation}

We want to show that (\ref{lastTbound}) is less than (\ref{finalNumCrit}), or equivalently that 
\begin{eqnarray} 0 & \leq & \left( \left( \frac{g-1+ \gamma b}{\NA+1} - ( \frac{g-1+ \gamma b}{N+1} - \epsilon) \right) u^2 + \left( 2 +  \frac{2g-2+ 2\gamma b}{\NA+1} - 2 \gamma b - \frac{7}{2}d \right) u \right. \nonumber \\
&& \left. \hspace{.5in}+ \left( 1 + \frac{g-1+ \gamma b}{\NA+1} - 2 \gamma b - \frac{7}{2}d \right) \right) v^2  - \left( \frac{g-1}{N+1} (u+1) \right) v. 
\end{eqnarray}

But the coefficient of $u^{2}$ in the coefficient of $v^{2}$ is $\epsilon > 0$.  So for all sufficiently large $u$, the polynomial 
\begin{displaymath} \epsilon u^2 + \left( 2 +  \frac{2g-2+ 2\gamma b}{\NA+1} - 2 \gamma b - \frac{7}{2}d \right) u + \left( 1 + \frac{g-1+ \gamma b}{\NA+1} - 2 \gamma b - \frac{7}{2}d \right)
\end{displaymath}
is positive; but then for all sufficiently large $v$, the polynomial 
\begin{eqnarray*} \left( \epsilon u^2 + \left( 2 +  \frac{2g-2+ 2\gamma b}{\NA+1} - 2 \gamma b - \frac{7}{2}d \right) u + \left( 1 + \frac{g-1+ \gamma b}{\NA+1} - 2 \gamma b - \frac{7}{2}d \right) \right) v^2  && \\
  \hspace{.5in} \mbox{} - \left( \frac{g-1}{N+1} (u+1) \right) v &&
\end{eqnarray*}
is positive, too.  Once again, we emphasize that the size of $u$ required depends on $d$, $g$, $\mathcal{B}$, and  $\epsilon$ but not on the particular curve $C$, the points $P_{i}$, the embedding $C \subset \Pro^N$, or the 1-PS $\lambda$. Similarly the size of $v$ required depends on $d$, $g$, $\mathcal{B}$, $\epsilon$ and $u$  but not on the particular curve $C$, the points $P_{i}$, the embedding $C \subset \Pro^N$, or the 1-PS $\lambda$.  $\hfill$ $\Box$

{\it Remark.}  Theorem \ref{stabilitytheorem} as stated does not establish stability for all large values of $m$, only for some large values of $m$.  Similarly, Gieseker's stability proof (\cite{G}, Theorem 1.0.0) only establishes stability for some, not all, large values of $m$.  In both cases it seems possible that one may be able to use variation of GIT arguments to conclude stability for all sufficiently large values of $m$, but I have not checked this.

\subsection{Application to the construction of moduli spaces}
\label{moduli section}

My motivation for studying this problem was to give GIT constructions of moduli spaces of weighted pointed stable curves.  We describe the parameter spaces and linearizations for this application now.  

Let $(C, P_{1},\ldots, P_{n}, \mathcal{A})$ be a weighted pointed stable curve with $n$ marked points.  Write $a := \sum a_{i}$, and assume that $2g-2+a > 0$.  Then for $\canonicalPower$ sufficiently large, $(\omega_{C}(\sum a_{i}P_{i}))^{\otimes \canonicalPower} =:  \mathcal{O}_{C}(1)$ is a very ample line bundle.  Write 
\begin{eqnarray*} V_{\canonicalPower,\mathcal{A}} & = & H^{0}(C, (\omega_{C}(\sum a_{i}P_{i}))^{\otimes \canonicalPower}) = H^{0}(C,\mathcal{O}_{C}(1)) \\
\dA & = & \deg \mathcal{O}_{C}(1) = \canonicalPower (2g-2+a) \\  
\NA +1 & = & \dim V_{\canonicalPower,\mathcal{A}} = \canonicalPower (2g-2+a) -g+1 \\
P(t) & = & h^{0}(C, \mathcal{O}_{C}(t)) = \dA t - g+1.
\end{eqnarray*}
Then $(C,  P_{1},\ldots,P_{n}, \mathcal{A})$ is represented by a point (in fact, many) inside the incidence locus $\incidenceLocus \subset  \Hilb(\Pro(V_{\canonicalPower,\mathcal{A}}),P(t)) \times \prod^{n} \Pro(V_{\canonicalPower,\mathcal{A}})$ where the points in the second factor land on the curve in the first factor.  In fact, $(C,  P_{1},\ldots,P_{n}, \mathcal{A})$ lies in a locally closed subscheme of $\incidenceLocus$ corresponding to weighted pointed curves embedded by $(\omega_{C}(\sum a_{i}P_{i}))^{\otimes \canonicalPower}$.

 It is very important to note that $d$, $N$, and $P(t)$ all depend on $g$, $n$,  $\mathcal{A}$ and $\nu$.  So, even if $g$ and $n$ are held constant, if $\mathcal{A}$ or $\nu$ varies, one is moving between loci in different Hilbert schemes---that is, one is using different parameter spaces---and this is not variation of GIT in the sense of Thaddeus and Dolgachev and Hu.  On the other hand, if $g$, $\nu$, and $\mathcal{A}$ are held constant and only $\mathcal{B}$ varies, this is VGIT in the sense of Thaddeus and Dolgachev and Hu.

I claim the following theorem, although the proof is not completely written down yet:

\begin{theorem} \label{modulitheorem} Suppose $g$, $n$, $d$, $\nu$, $\mathcal{A}$, and $\mathcal{B}$ fit the setup of this paper and satisfy the hypotheses of Theorem \ref{stabilitytheorem}.  Let $\gamma = \nu / (2 \nu -1)$.   Suppose $\nu \geq 5$ and $d= \nu (2g-2+a)$,  and let $J$ be the locus in $I$ where $\mathcal{O}(1) \cong (\omega( \sum a_{i} P_{i}))^{\nu}$. Then:
\begin{itemize}
\item  If $\mathcal{A} = \mathcal{B}$ and $ b_{i} \leq 1$, then $J \dblq SL(N+1) \cong \MgA$.
\item In particular, if $\mathcal{A} = \mathcal{B}$ and $\frac{1}{2} + \epsilon < b_{i} < \frac{1}{2 \gamma}$ for each $i = 1,\ldots,n$, then $J \dblq SL(N+1) \cong \Mgn$. 
\end{itemize}  
\end{theorem}

How much of Theorem \ref{modulitheorem} has been checked?  I believe all that is needed is extremely minor changes to the Potential Stability Theorem of \cite{BS}.  It should still say that nothing ``bad'' can be GIT stable; the argument is very long, so I have not checked all of it, but it is also extremely robust, and I am very confident that it will work.   One can easily write down the ``Basic Inequality'' when there are weighted marked points.  I have done this, and checked that the condition on points colliding agrees exactly with the definition of $\MgA$, and that the argument that $J^{ss}$ is closed inside $I^{ss}$ still goes through.  It then follows that all weighted pointed stable curves are GIT stable, justifying the title of this paper and completing the proof of Theorem \ref{modulitheorem}.

If $g$, $\nu$ and $\mathcal{A}$ are held fixed and the set of linearizing weights $\mathcal{B}$ is allowed to vary sufficiently far from $\mathcal{A}$, the quotient may undergo a flip.  Identifying these quotients is a project I am currently working on.

\section{Additional remarks (Director's cut)}
In the course of my research I have learned a little bit more about this problem than just what appears in this paper.  In particular, I relate my proof to Gieseker's in the unpointed case, and this leads to a conjecture about the worst 1-PS.  Next, I mention two suggestions for improving the main result, one that I expect would not work, and one that probably would.

\subsection{Comparison to Gieseker and Morrison's results, and the worst 1-PS}
We may interpret Gieseker's proof (\cite{G}, Theorem 1.0.0) as the $n=0$, $q=1$ case of Theorem \ref{stabilitytheorem}.  This easily leads to a coarse upper bound for $T$.  The bound so obtained is not quite as good as the bound given in \cite{Morr}, Section 4 and used in Gieseker's proof.  However, after running the proof here, one can perform their analysis on top of that, and the resulting bounds for the leading coefficient would then agree.

Kempf and Rousseau showed that when $x$ is GIT-unstable, there is a ``worst 1-PS'' destabilizing $x$.  This suggests the following strategy for proving stability: suppose for purposes of contradiction that $x$ is unstable, then find the worst 1-PS, then show that it is actually not destabilizing.  Morrison and I have never gotten this strategy to work in our situation (we can't find the worst 1-PS, for the same reason that we can't compute the absolute weight filtration discussed in Section \ref{profilessection}).

However, we can describe the 1-PS for which it is most difficult to prove stability using our methods: it is the 1-PS for which there is only one point $Q_{1} = P_{i}$ in the base locus of $V_{\bar{N}}$, where $b_{i}$ is the largest value in $\mathcal{B}$, every stage of the filtration is a complete sublinear series of $H^{0}(C,\mathcal{O}(1))$, and the weights are linearly decreasing (hence, uniquely determined by the conditions that they decrease to zero and sum to 1).

Of course, just because it is hard for us to show that this 1-PS is stable does not mean it is actually the worst 1-PS, but it certainly is a candidate.  I believe it would be an interesting to show either that this is the worst 1-PS, or exhibit another 1-PS which is worse.  In the meantime, I mention this 1-PS for its value  as a heuristic test for GIT stability for parameter spaces and linearizations where this is currently unknown, and for testing putative stability proofs.

\subsection{Can we improve these results if we use a more complicated filtration than $\tilde{V}_{\bullet}$ as scaffolding?}
{\it Q:}  We only take the span of  ``three-layer'' spaces $V_{s}^{\alpha uv}V_{t}^{(1-\alpha)uv}V_{0}^{v}$.  Could we get any further improvement by defining a filtration using spaces of the form $V_{s}^{\alpha uv}V_{t}^{\beta uv}V_{w}^{(1-\alpha - \beta)uv}V_{0}^{v}$?

{\it A:}  There may be room for improvement of our results, but when $m$ is large, adding more layers will not buy you anything.  We never really asked what is the best way to produce a basis.  We always began with a space of the form $V_{k}^{\alpha uv}V_{k+1}^{(1 - \alpha) uv V_{0}^{v}}$ having weight $\alpha r_{k} uv + (1-\alpha)r_{k+1}uv + vr_{0}$ and asked the question: for what choice of $\beta_{j}$ for $j=0$ to $N$ will \\ $\codim \mbox{span} (V_{k}^{\alpha uv}V_{k+1}^{(1 - \alpha) uv}, V_{0}^{\beta_{0}uv}V_{1}^{\beta_{1}uv} \cdots V_{N}^{\beta_{N}uv}V_{0}^{v})$ be minimized?  

There are constraints.  First, $\sum_{j=0}^{N} \beta_{k} = 1$.  Also, the weight of the second space in the span should be less than or equal to that of the first, so 

\[ \beta_{0}r_{0} + \cdots + \beta_{N}r_{N} \leq \alpha r_{k} + (1-\alpha)r_{k+1}.  
\]

These conditions give a polytope in $\beta$-space.  Minimizing the multiplicity of each $P_{i}$ means minimizing the linear function 
\[ f(\beta_{0},\ldots,\beta_{N}) = c_{1,i}\beta_{1} + c_{2,i}\beta_{2} + \cdots+ c_{N,1}\beta_{N} 
\]
over this polytope.  The minimum must occur on the boundary, specifically at one (or more) of the vertices of the polytope, and these are precisely the ``three-layer'' spaces.  

The argument just given should be approximately true when $m$ is very large and divisible (so that all the exponents are integers), but it could break down badly for small $m$.  So, for small $m$ stability, we might want to consider filtrations which are much more complicated than those used in this paper.

\subsection{Lower convex envelopes might give better bounds for $T$} 
\label{lowerenvelopesarebetter}
Recall from Section \ref{xtilde section} that in the definition of $\tilde{X}_{\bullet}$, we do not minimize the multiplicity of each $Q_{i}$.  In fact, it is not hard to find the minima; instead of using the functions $s(k,i)$ and $t(k,i)$, defined as `` `previous' and `next' among values where the multiplicity of $Q_{i}$ jumps,'' we should instead use $\sigma (k,i)$ and $\tau (k,i)$, defined as  `` `previous' and `next' among values where the multiplicity of $Q_{i}$ jumps which lie on the lower envelope of these.''  That is, there are $q$ lower envelopes to keep track of.

It is possible that if one defines a filtration $\tilde{Y}_{\bullet}$ using lower envelopes like this, one might be able to prove stability under a weaker hypotheses than those used in this paper.  In particular I believe that this might yield a proof of asymptotic stability of canonically embedded smooth nonhyperelliptic curves.  The obstacle is the proof of the analogue of Lemma \ref{creep}.  I can't figure out how to get this to work if you use lower envelopes instead of just the next value; instead of relating everything to $c_{j_{1}+1}$ one would need to work with much later $c$'s, and I don't see how to do this.






\begin{thebibliography}{MMMMM}

\setlength{\baselineskip}{12pt} 

\bibitem[AG]{AG}Alexeev, V. and G. M. Guy.  ``Moduli of weighted stable maps and their gravitational descendants.''  math.AG/0607683.

\bibitem[BM]{BM}Bayer, A. and Y. Manin.  ``Stability conditions, wall-crossing and weighted Gromov-Witten invariants.'' math.AG/0607580.

\bibitem[BS]{BS}Baldwin, E. and D. Swinarski.  ``A geometric invariant theory construction of moduli spaces of stable maps.''  arXiv:0706.1381.

\bibitem[Gies]{G}Gieseker, D.  {\it Lectures on Moduli of Curves.}  Tata Institute Lecture Notes, Springer, 1982.  

\bibitem[G2]{GCime}Gieseker, D.  ``Geometric invariant theory and applications to moduli problems.'' 45--73, LNM {\bf 996}, Springer, 1983.  

\bibitem[Gotz]{Gotz}Gotzmann, G.   ``Eine Bedingung f\"{u}r die Flachheit unda das hilbertpolynom eines graduierten Ringes.''  {\it Math. Z.} {\bf 158} (1978), 61--70. 

\bibitem[HM]{HM}Harris, J. and I. Morrison. {\it Moduli of Curves.} Graduate Texts in Mathematics \textbf{ 107}, Springer, 1998.

\bibitem[Hass]{Hass}Hassett, B.  ``Moduli spaces of weighted pointed stable curves.''  {\it Adv. Math.} {\bf 173} no. 2 (2003), 316--352.

\bibitem[Morr]{Morr}Morrison, I.   ``Projective Stability of Ruled Surfaces.''  {\it Inv. Math.} {\bf 56} (1980), 269--304.

\bibitem[MM]{MM}Musta\c{t}\u{a}, A. and A. Musta\c{t}\u{a}.  ``Intermediate moduli spaces of stable maps.''  {\it Invent. Math.} {\bf 167} no. 1 (2007), 47--90.

\bibitem[GIT]{GIT}Mumford, D., Fogarty, J. and F.C. Kirwan.  {\it Geometric Invariant Theory.}  Third Edition.  Springer, 1994.  

\bibitem[Mum]{StProjVar}Mumford, D.  ``Stability of Projective Varieties.''  {\it Enseignement Math. (2)} {\bf 23} (1977), no. 1-2, 39--110.


\end{thebibliography}
\end{document}